\documentclass[11pt]{amsart}
\usepackage[margin=1in]{geometry}

\usepackage{mmacells}
\mmaDefineMathReplacement[≠]{!=}{\neq}

\usepackage{scalefnt}

\usepackage{graphicx}
\usepackage{wrapfig}

\usepackage{menukeys} 
\usepackage[utf8]{inputenc}

\usepackage{caption}

\usepackage{tikz}

\newenvironment{myquotation}%
{\list{}{\leftmargin=0.4in\rightmargin=0.4in}\item[]}%
{\endlist}

\usepackage{enumitem}
\setlist[description]{leftmargin=2\parindent,labelindent=\parindent}

\theoremstyle{definition}

\theoremstyle{remark}

\numberwithin{equation}{section}

\begin{document}

\title[Encouraging student creativity in mathematics]{Encouraging student creativity in mathematics\\ through 3D design and 3D printing}

\author{Christopher R.\ H.\ Hanusa}
\address{Department of Mathematics \\ Queens College (CUNY) \\ 65-30 Kissena Blvd. \\ Flushing, NY 11367\\ United States}
\email{chanusa@qc.cuny.edu}
\thanks{The author was supported in part by NSF Grant DUE-1928565.}






\begin{abstract}
This is a case study of teaching 3D design and 3D printing in a project-based computing course for undergraduate math majors. This article discusses content organization, implementation, project grading, and includes a personal reflection. There is an emphasis on lessons learned and how to encourage student creativity and artistic expression. An appendix details 3D design techniques in Mathematica.
\end{abstract}

\maketitle

\section*{Introduction}
As mathematicians we recognize that the research process can be messy and non-linear, that productive failure is a key part of achieving a new result, that creativity and pushing boundaries is necessary to our craft, and that there is beauty in the elegance of mathematical results and visualizations. Now ask yourself: How many of our courses let 

\begin{wrapfigure}[16]{R}{0.34\textwidth}
	\centering
		\includegraphics[angle=270,origin=c,height=0.3\textwidth]{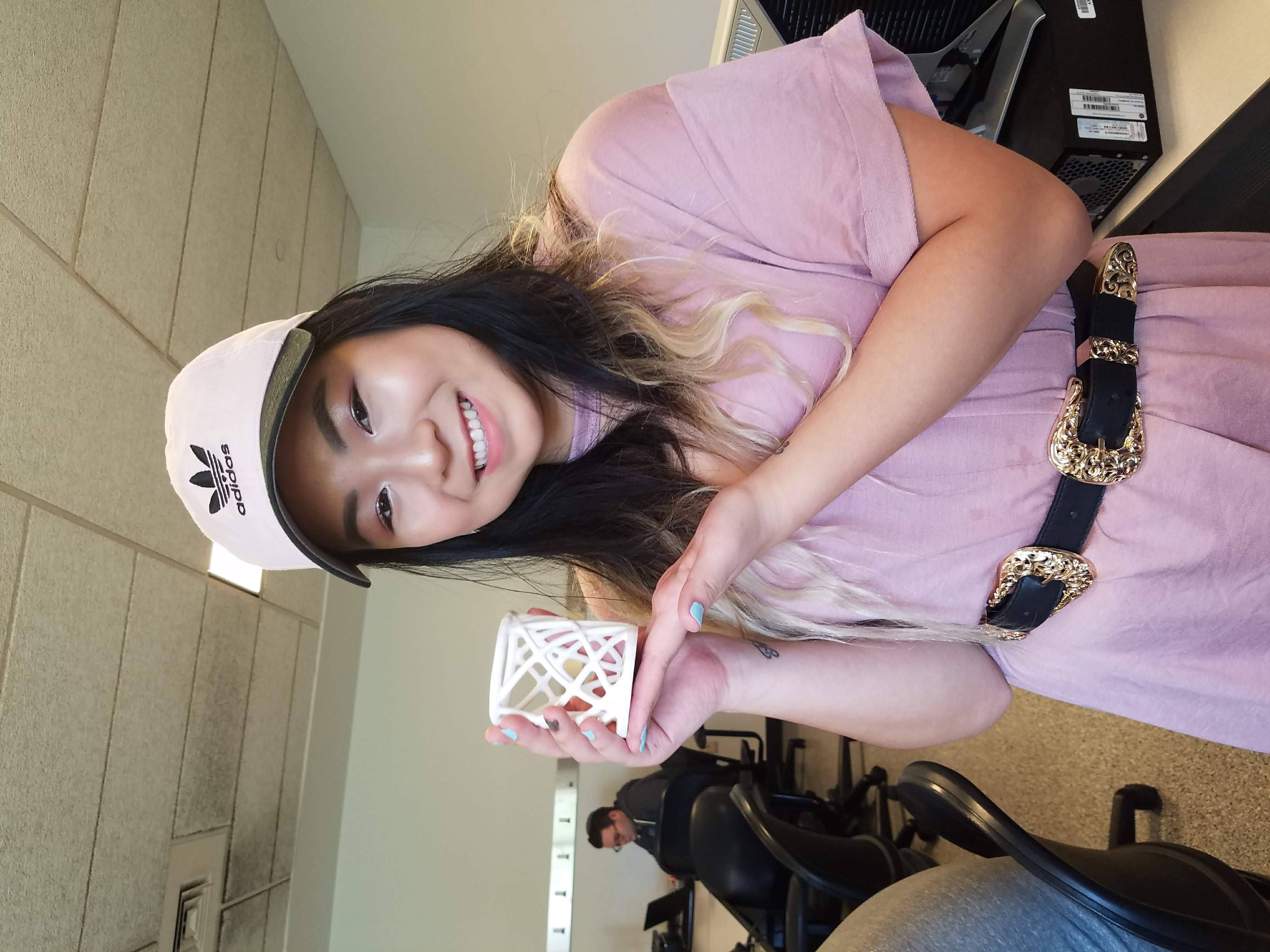}
		
\begin{minipage}{0.33\textwidth}
	\caption{Jenny Xu with her 3D printed ceramic pot.}\end{minipage}
\end{wrapfigure}

\noindent
students experience these aspects that we find so fundamental in mathematics?

This article describes the project-based Mathematical Computing course at Queens College that aims to give students this experience.  In this course, students learn the basics of programming in Wolfram Mathematica and actively apply their skills to design mathematical sculpture and create an interactive app. Section~1 describes the structure of the course, including learning objectives and ensuring student success. Section~2 details the philosophy behind and the implementation of the 3D design and 3D printing module in the course, including how students learn about design thinking and the elements of art. Section~3 discusses the written lab report and the way the project deliverables are graded. Section~4 is a personal reflection on the course. Finally, Appendix~\ref{sec:3DDesign} summarizes key techniques for 3D design in Mathematica.

\section{Course Structure}

\subsection{Student Population}
Queens College is located in its namesake borough of Queens in New York City and is one of the 25 campuses of the City University of New York. The student body of Queens College of approximately 18,000 mirrors the population of Queens, the most ethnically diverse county in the United States. Among freshmen, 44\% are immigrants born outside the mainland US, and 34\% are first generation students. Queens College serves as an engine of upward mobility and is proud to offer a high-quality education to students at a fraction of the cost of private institutions.

The Queens College Department of Mathematics has 31 full-time faculty members and many part-time instructors. There are approximately 450 math majors in five specialties: pure mathematics, applied mathematics, data science and statistics, secondary education, and elementary education. The largest cohorts are applied mathematics and secondary education.

The course Mathematical Computing (MATH 250) is offered once per year with 75-minute classes twice a week during a 14-week semester. The course has an enrollment between 15 and 25 students, most of whom are juniors or seniors. This course fulfills the computing requirement for applied mathematics majors and serves as a math elective for all other math majors. The prerequisite for the course is either multivariable calculus or linear algebra, to ensure that students have prior experience in multi-dimensional thinking and visualization. The course has no computer programming prerequisite in order to serve the largest possible audience. At the same time, a handful of enrolled students do have extensive programming experience as they are majoring or minoring in computer science. 

\medskip
\begin{figure}[h]
	\includegraphics[width=0.8\textwidth]{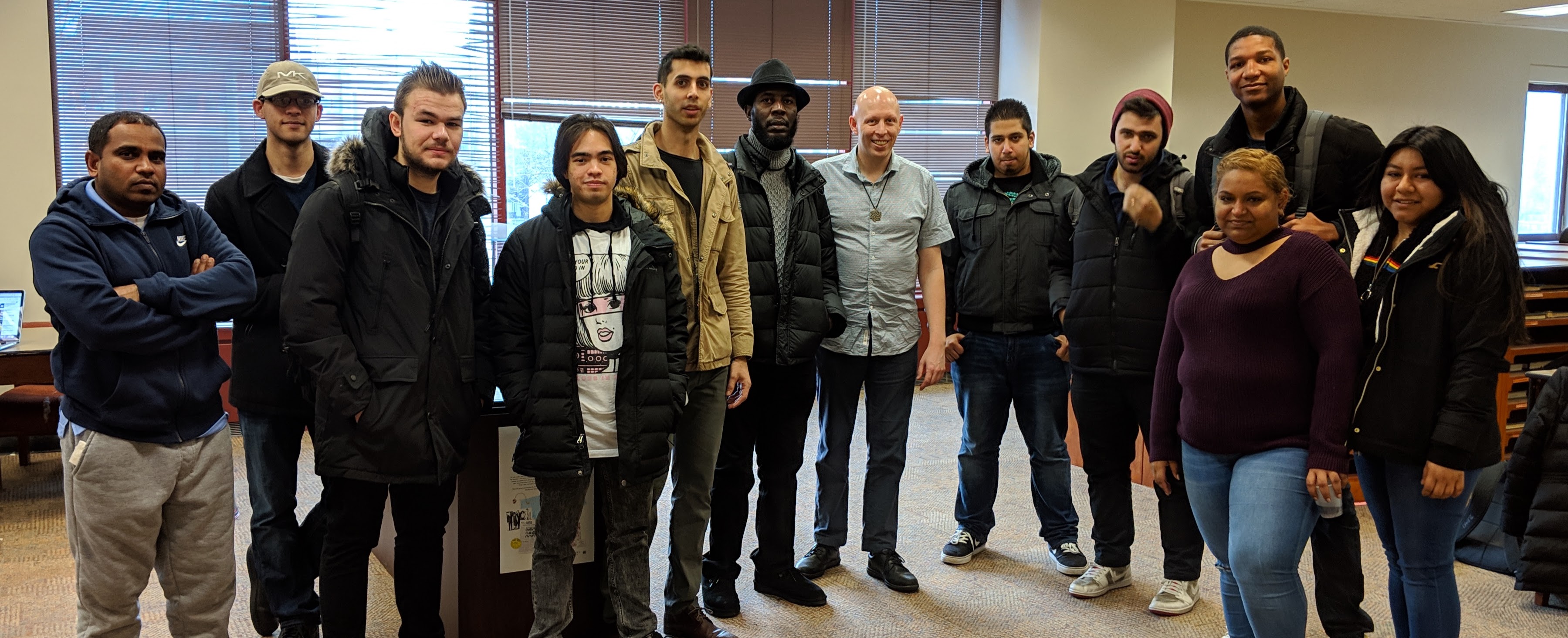}
	\caption{The Fall 2018 Mathematical Computing class.}
\end{figure}

\subsection{Course Materials}

The current structure of the course was initiated in Spring 2015. There is no textbook for this course; the instructor provides freely available Mathematica tutorials and accompanying video lectures on his webpage.\footnote{\url{http://qc.edu/~chanusa/courses/250/21/}} The students contribute to an online discussion platform during the semester for pre-class activities; past platforms have included Blackboard, Microsoft Teams, and Campuswire. 

The City University of New York has a university-wide site license for Wolfram Mathematica, which means all Queens College campus computers have Mathematica installed, and all students may install Mathematica on their home machines, so software access and cost are not barriers for our students.  

Students have access to the Queens College Makerspace, a new space located in the library that houses approximately eight 3D printers, a 3D scanner, a laser cutter, an embroidery machine, a CNC machine, pen plotters, electronic supplies, and a variety of workshop tools. Once students have completed the orientation, they may freely use any of the materials and machines.

Students are required to pay to print their final model through a third-party 3D printing company such as Shapeways. This is a cost of approximately \$20-\$60 dollars, depending on the size of the model and the material they choose, and is their only expense for the course.

\subsection{Learning Objectives}
\label{sec:objectives}

Our department believes that every math major should gain some programming experience. For the vast majority of students, this class is their first exposure to any mathematical software. So the most fundamental learning objective is that students learn the basics of programming in Mathematica, developing fluency with the key data structures of lists and functions. Students learn how to design, code, run, test, and debug computer programs. Along the way, students are developing good programming techniques, including sectioning and documenting their code. Students also gain skills to become confident and self-sufficient learners through the use of Mathematica's extensive documentation and internet searches.

Another important objective is that students develop mathematical, programming, and problem-solving skills. For example, students gain a deeper understanding of three-dimensional geometry, including coordinate systems, multivariable functions, and three-dimensional objects. From a programming perspective, students learn and apply functional programming, which is new even for experienced programmers. Building on this knowledge, students apply basic problem-solving skills including analyzing problems, modeling a problem as a system of objects, creating algorithms, and implementing solutions in Mathematica. 

A final objective is for students to use mathematics and programming for experimentation, as creative tools, and collaboratively. They apply the design process and communicate the decisions made therein, including ideation, artistic principles, prototyping, and revisions. Students advance teamwork skills by collaborating with classmates, discussing and solving problems in a group setting, and practicing giving and receiving constructive feedback. These ideas will be developed further in Section 2.

\subsection{A Day in the Life}

As a project-based course, class time is balanced between content acquisition and project work. The semester breaks down into three modules, each of which culminate with an individual project. 

At the beginning of each module, classes are centered around learning the course content. Before each class, students work through a tutorial notebook with an accompanying video lecture, created by the author. Students are encouraged to ask questions about the content through the online discussion board. (A tutorial structure works well given the disparate levels of programming background because the tutorials can be consumed at each student's individual speed.) During class time, students work together in groups on comprehension and challenge questions that reinforce key ideas from the day's tutorial and lecture. During this time the instructor travels from group to group, addressing any questions and, if necessary, gently nudging students in the correct direction.  Each class ends with a debriefing of the different ways to approach solving the questions.

As each module progresses, class time turns to project work. Students work on their projects outside class and bring their progress to class for instructor and peer consultation and feedback. As the module nears its end, the instructor spends much of the class period troubleshooting and debugging code and ensuring students are on track to complete their project by the deadlines. Once the students finish their final draft, they complete a peer review with one or more classmates and revise their work before submitting their work for grading. 

\medskip
\begin{figure}[h]
	\includegraphics[width=0.8\textwidth]{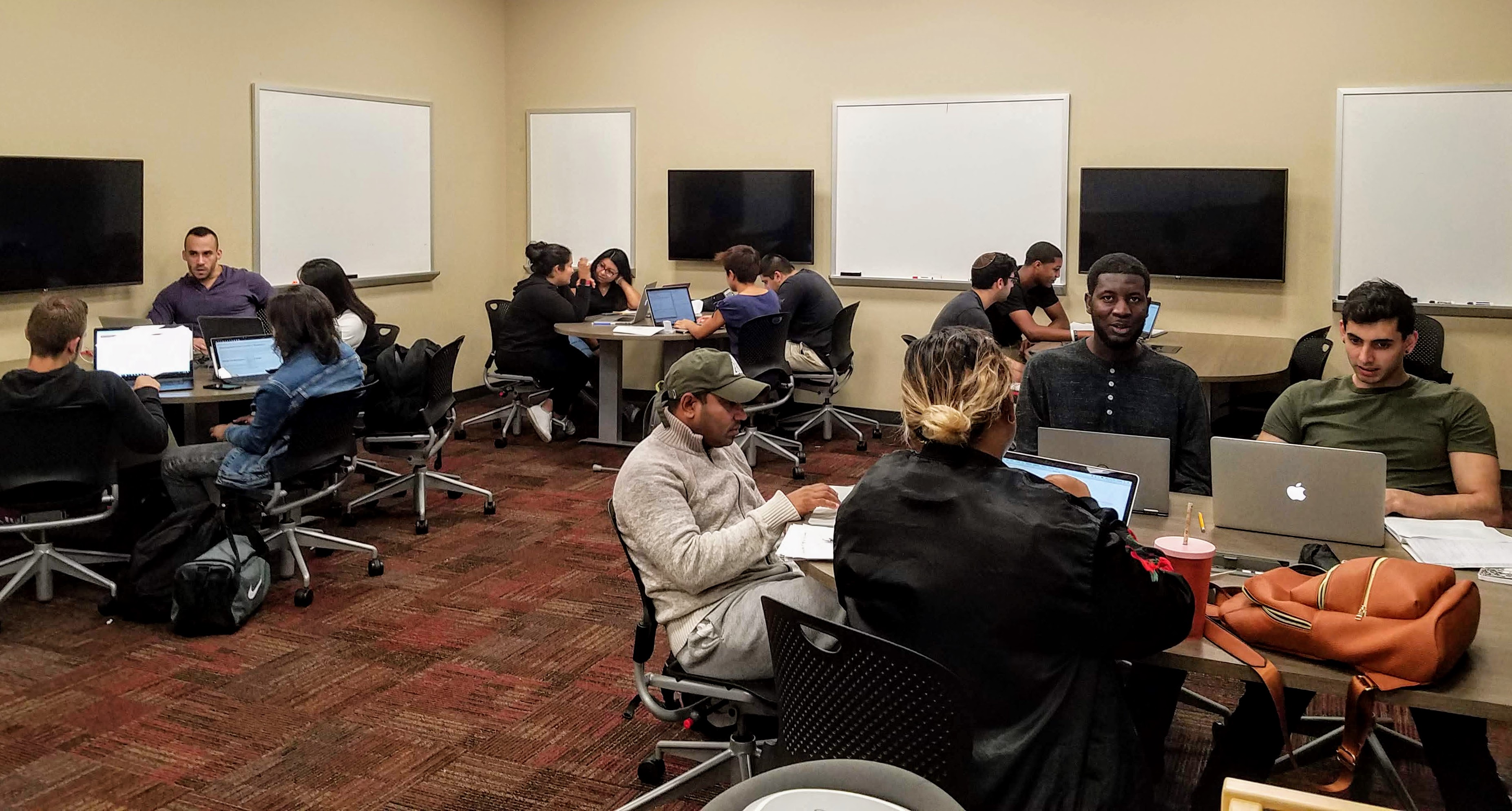}
	\caption{Students working in groups on their projects.}
\end{figure}

\subsection{Course Content Overview}

The first module of the class teaches students the basics of working with Mathematica, exposes students to the wide capabilities of Mathematica, and helps students develop fluency with lists, list manipulation, and functions. The first project involves creating a tutorial centered around a coherent set of Mathematica commands that each student chooses. The instructor encourages students to choose a topic related to the mathematics they have seen in previous courses or from one of the many subject domains to which Mathematica applies.  Past student projects include visualizing graphs of multivariable functions, matrix operations, statistical methods, image processing, machine learning, finance, geography, and sound.

The second module of the class is centered around two- and three-dimensional graphics, 3D design, and the 3D printing process. Students design and 3D print a mathematical sculpture. This is the focus of the remainder of this article. 

In the third module of the class students learn about and apply the interactive capabilities of Mathematica. The third project involves creating an interactive program similar to an app you might use on your phone. Students determine and design the entire interface and program the behavior of the interface depending on the user input. The resulting programs always vary widely and use a broad spectrum of Mathematica's visual, audio, and interactive capabilities.  Past student projects include an interactive map, a piano simulator, a fractal explorer, trivia guessing games, two-player games like dots and boxes or tic-tac-toe, and action games like pong, snake, or boulder dash.

\subsection{Ensuring Student Success}

The most important part of ensuring student success is making sure that each student is working on an individually fruitful and feasible project. It must be fruitful in that the student develops programming and problem solving skills and feels a sense of accomplishment from their work. It must be feasible in that each project must be achievable in the given timeframe with the student's skill set. To that end, the instructor must help each student develop a well-defined project goal of the right scope.  This is described in detail in Section~\ref{sec:design}.

During this project work phase, the instructor is simultaneously overseeing 15 to 25 different research projects. For this reason, it is important that the students be able to function mostly independently or with help from their peers. It has been crucial to provide a large amount of scaffolding around every project, including outlining transparent expectations and setting deadlines to keep the students on track. The course webpage provides details about how far along students should be each class period and this is reinforced in class by the instructor. This signals to students who are behind that they need to put in more time to the project and/or come to office hours for extra help. 

Because the students have such diverse prior knowledge and backgrounds, the fruitfulness and feasibility of a project will vary from student to student and the grading system must accommodate this disparity. In particular, if project expectations were uniform across all students, this might encourage advanced students to put in minimal work and discourage beginning students. Details about project grading are discussed in depth in Section~\ref{sec:grading}.

\section{The 3D Design Project}

\subsection{Overview}
\label{sec:projectoverview}

In the second module of the course, students are tasked with designing and 3D printing a mathematical sculpture. Each student artwork is expected to originate from some mathematical concept. The artwork may visualize an idea from mathematics, or instead use mathematical and computational techniques to recreate a specific object or abstract form of interest. 

In contrast to more established computer-aided design programs (such as Solidworks,  AutoCAD, or Fusion 360), 3D design in Mathematica relies on the user building their model from base mathematical concepts. To be able to do this, students must refresh their mathematical knowledge of three-dimensional coordinate systems, curves, surfaces, and geometric objects.

The author has created Mathematica tutorials on 2D graphics objects, 3D graphics objects, parametric curves and surfaces, transformations thereof, and how to make each of these types of objects 3D printable. (See Appendix~\ref{sec:3DDesign} for more information.) These ideas are interwoven with discussions of the elements of art, design thinking, and the 3D printing process.

\subsection{Creative Expression and Artistic Intentionality}
\label{sec:creativity}

The author cares deeply that students get a chance to use mathematics creatively and learn to think artistically. A few days into this second module, students are prompted to explore examples of mathematical art online, including works by Bathsheba Grossman\footnote{\url{http://www.shapeways.com/shops/bathsheba}}, Henry Segerman\footnote{\url{http://www.shapeways.com/shops/henryseg}}, Laura Taalman\footnote{\url{https://mathgrrl.com/designs/}}, and various artists exhibiting at the Bridges Mathematical Art Conference\footnote{\url{http://gallery.bridgesmathart.org/exhibitions/}}. They are asked to choose a few inspiring artworks, and share and annotate them with a few sentences about why they chose the pieces on the online discussion board before class. Here are a few examples of past student comments:

\smallskip
\begin{myquotation}
	``I find the following pieces of artwork interesting and inspiring because I can imagine the pieces being utilized as jewelry as a statement piece. The 3D prints in my option are very unique and striking visually.''
\end{myquotation}

\smallskip
\begin{myquotation}
	``This gif unlocked memories I forgot I had. I remembered learning about flowers having symmetry, patterns, and some even follow the Fibonacci sequence. They are mathematical and stunning!''
\end{myquotation}

\smallskip
\begin{myquotation}
	 ``[This sculpture] shows the beauty of the randomness and reorganization of water. This piece is not created using a 3D printer, however, I can see something like this being designed and created on a 3D printer.''
\end{myquotation}

\begin{wrapfigure}[13]{R}{0.31\textwidth}
	\centering
	\includegraphics[origin=c,width=0.26\textwidth]{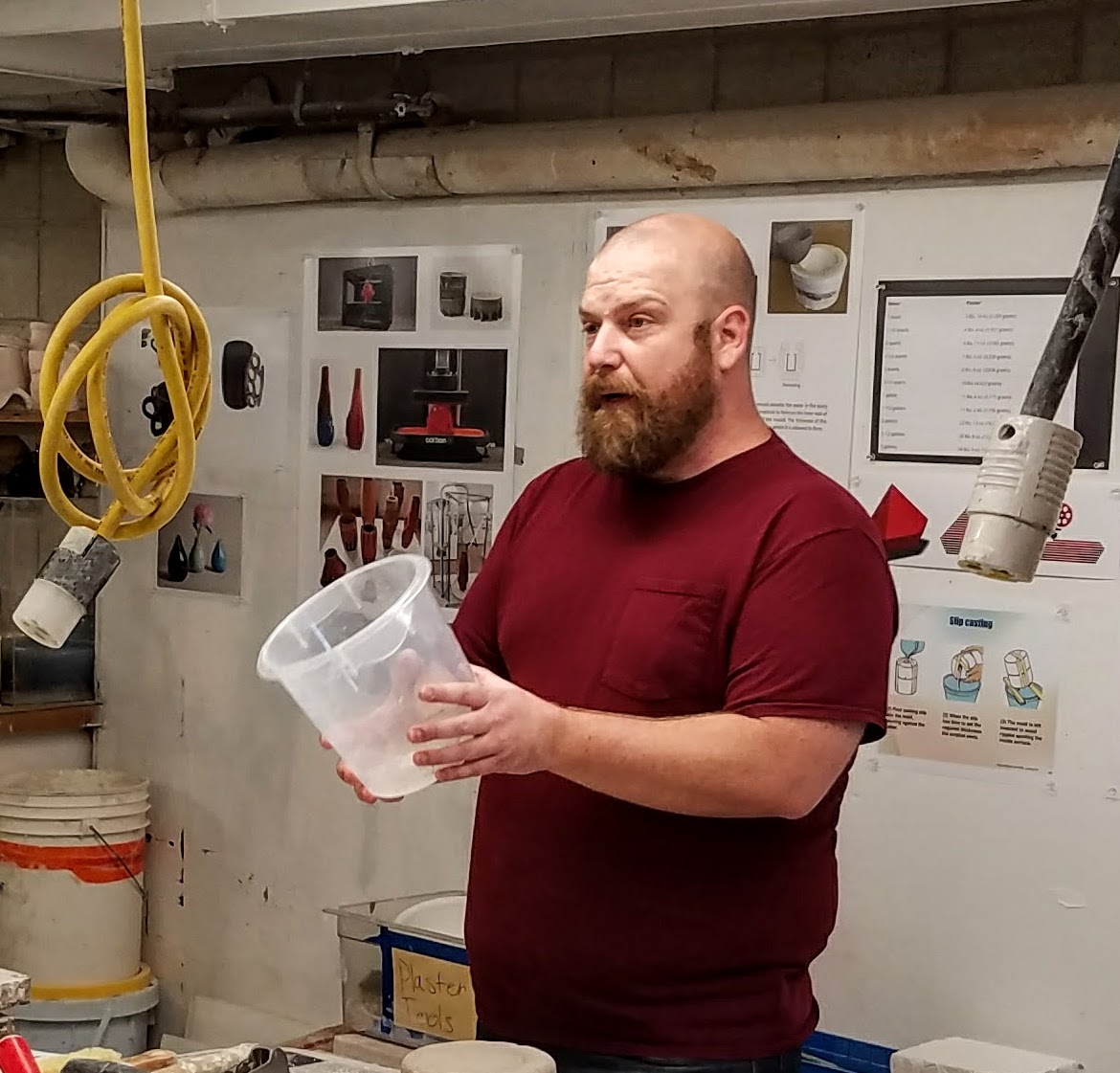}
	\vspace{-.2in}
	\begin{minipage}{0.29\textwidth}\caption{Professor of Art Matthew Greco.}\end{minipage}
\end{wrapfigure}

The next class session the students are in for a treat. A colleague from the art department comes to the classroom and gives a presentation about 3D printing and the artistic elements of sculpture, including the concepts of positive and negative space, symmetry vs. asymmetry, proportion, repetition, contrast, harmony, and movement. After the structured presentation, the class discusses the artistic elements that are present in the mathematical art the students had shared the night before. Because it is often the first time many students have thought about these ideas, they are always engaged and inspired by this discussion.

This is supplemented by handouts on the Elements of Art\footnote{
\url{https://www.getty.edu/education/teachers/building_lessons/formal_analysis.html}}
and the Principles of Design\footnote{
\url{https://www.getty.edu/education/teachers/building_lessons/formal_analysis2.html}} from the Getty Museum's education initiative. Once students are given these insights to the world of art and design, they can use this lexicon to express their goals in the creative process. Do they want to aim for a piece that has distinctive negative space? Are they looking for a piece that looks lightweight or instead has gravitas? Through this exposure, students gain the ability to express key aspects of inspiring artwork, become more intentional about how their sculpture develops, and enrich their final reflections.

\subsection{Design Thinking}
\label{sec:design}

``3D printing'' sounds like it could be simple$\hdots$ until you try it. The process involves many messy and complex steps between the initial seed of an idea and the final physical piece. Students confront the full 3D printing experience in this class, including working through productive failure at many levels. 

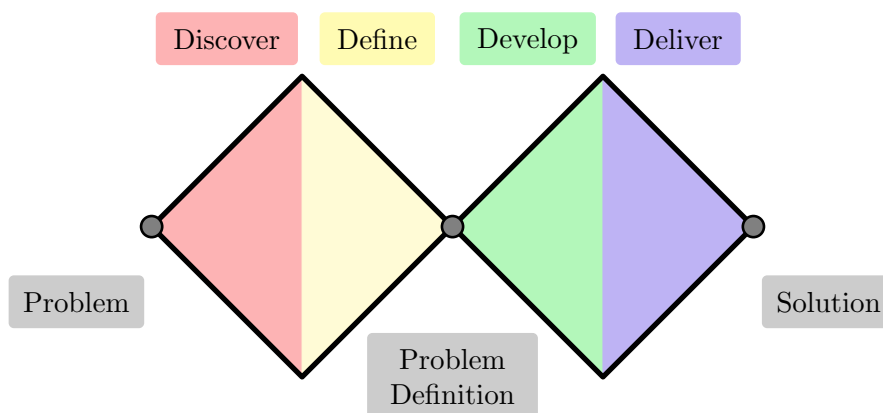
\begin{figure}[ht!]
\begin{center}
\begin{tikzpicture}
	\filldraw[fill=purple!10!red!30, draw=purple!10!red!30, line width=.1pt] (0,0) -- (2,2) -- (2,-2) -- (0,0);
	\filldraw[fill=yellow!20, draw=yellow!20, line width=.1pt] (4,0) -- (2,2) -- (2,-2) -- (4,0);
	\filldraw[fill=blue!10!green!30, draw=blue!10!green!30, line width=.1pt] (4,0) -- (6,2) -- (6,-2) -- (4,0);
	\filldraw[fill=purple!20!blue!30, draw=purple!20!blue!30, line width=.1pt] (8,0) -- (6,2) -- (6,-2) -- (8,0);
	\draw[line width=2pt, line join=round] (0,0) -- (2,2) -- (4,0) -- (6,2) -- (8,0) -- (6,-2) -- (4,0) -- (2,-2) -- (0,0);
	\filldraw [fill=gray, draw=black, line width=1pt] (0,0) circle [radius=4pt]
					 (4,0) circle [radius=4pt]
					 (8,0) circle [radius=4pt];
	\draw (1,2.5) node [fill=purple!10!red!30,rounded corners=2pt,minimum size=20pt]{\phantom{i}Discover\phantom{i}};
	\draw (3,2.5) node [fill=orange!10!yellow!30,rounded corners=2pt,minimum size=20pt]{\phantom{i}Define\phantom{i}};
	\draw (5,2.5) node [fill=blue!10!green!30,rounded corners=2pt,minimum size=20pt]{\phantom{i}\raisebox{-.03in}{Develop}\phantom{i}};
	\draw (7,2.5) node [fill=purple!20!blue!30,rounded corners=2pt,minimum size=20pt]{\phantom{i}Deliver\phantom{i}};
	\filldraw (-1,-1) node [fill=black!20,rounded corners=2pt,minimum size=20pt]{\phantom{\,}Problem\phantom{\,}};
	\filldraw (4,-2) node [fill=black!20,rounded corners=2pt,minimum size=20pt]{\begin{tabular}{c}\raisebox{-.015in}{Problem}\\\raisebox{-.015in}{Definition}\end{tabular}};
	\filldraw (9,-1) node [fill=black!20,rounded corners=2pt,minimum size=20pt]{\phantom{\,}Solution\phantom{\,}};
\end{tikzpicture}	
\end{center}
\caption{The double diamond design model}
\label{fig:doublediamond}
\end{figure}

Students are introduced to a macro vision of this adventure from start to finish using the ``double diamond'' framework in Figure~\ref{fig:doublediamond}.  This breaks down the process into four stages: Discover, Define, Develop, and Deliver. The first diamond is where a vague initial concept is made into a concrete problem specification. The second diamond is where the idealized problem specification meets practicality and feasibility to be physically realized.

\subsubsection{Discover Stage}

The process starts with the students receiving the very vague project description given in the first paragraph of Section~\ref{sec:projectoverview}. The students are primed to think artistically after the web browsing and discussion of the elements of art.  At that point, students are encouraged to explore, research, and brainstorm to create many possible project directions that touch on mathematical principles or outside interests. 

\subsubsection{Define Stage}

This brings us to the widest part of the first diamond and the next stage in the process: problem definition. These many directions created from the discover stage need to be triaged based on personal preference, time constraints, ambition, and peer feedback. 

Students propose initial project ideas over email by a specified deadline; the next class period is spent helping students hone these initial ideas to precise project goals. This requires that the instructor understands Mathematica's capabilities and can envision how the student could complete the project when they propose it. Furthermore, the instructor must informally assess each individual student's level of knowledge to make sure that the steps to completion are not too simple nor too complex.

After multiple rounds of back and forth, the result is a precise problem specification that the student can always refer back to if they need to remind themselves what they are trying to accomplish.

\subsubsection{Develop Stage}

We have reached the point in between the two diamonds and it is time for students to branch out once again. They need to take the precise idea and understand how they are going to realize it. Students are encouraged to sketch their idea on paper first and develop a plan before starting to code. They then develop the necessary coding skills and mathematical knowledge to explore the space of possibilities. They also make use of a prototyping feedback loop to learn what works and doesn't work in the 3D printing world.

This process has been enhanced by the new Queens College Makerspace, first incorporated into this course in the Fall 2021 semester. Students use the 3D printers for ``rapid prototyping'' of their initial designs, which can be printed and critiqued within a day's time.

\begin{figure}[h]
	\includegraphics[width=0.75\textwidth]{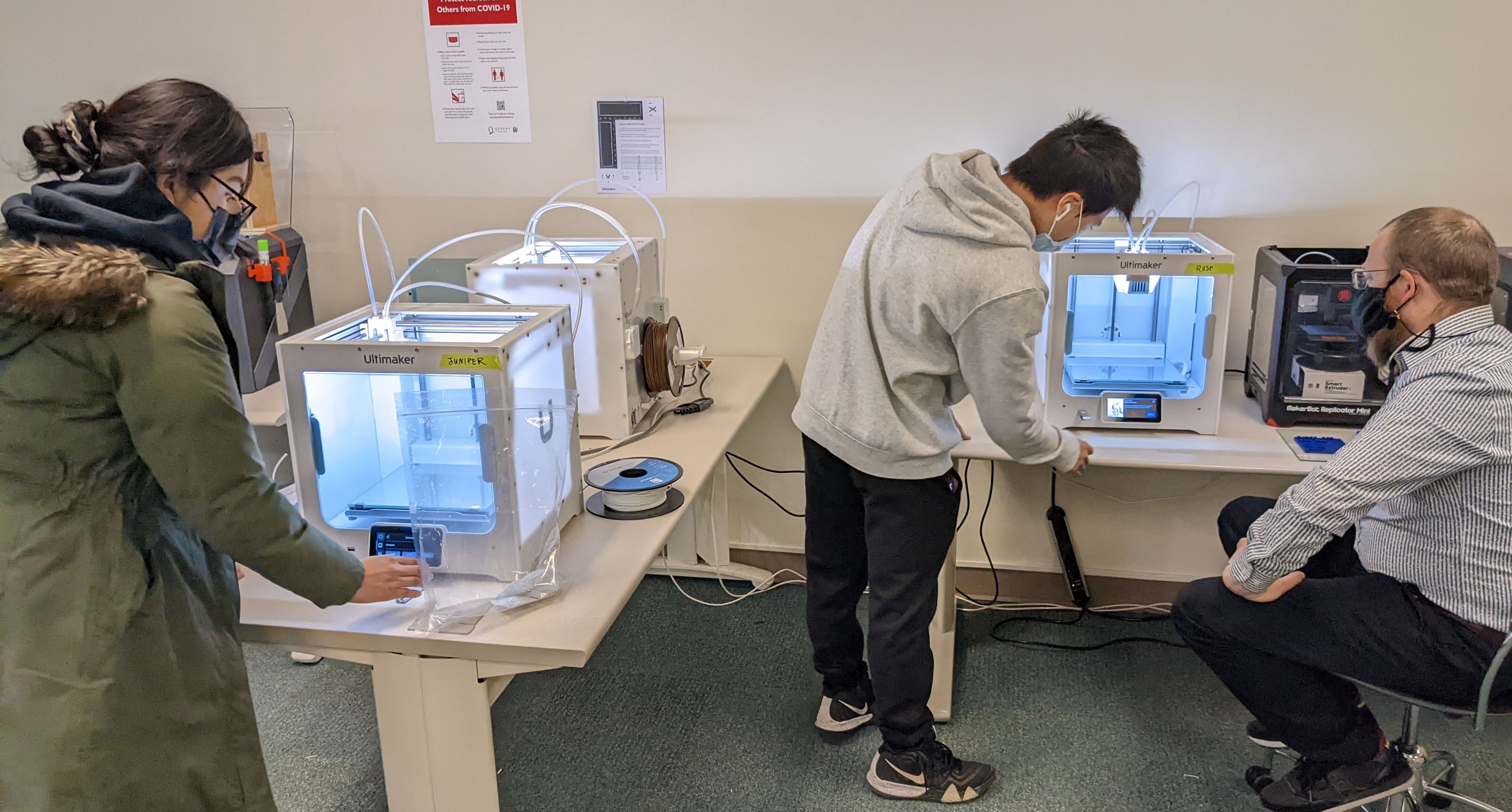}
	\caption{Students working in the Queens College Makerspace.}
\end{figure}

\subsubsection{Deliver Stage}

\begin{wrapfigure}[11]{R}{0.32\textwidth}
	\centering
	\includegraphics[origin=c,width=0.3\textwidth]{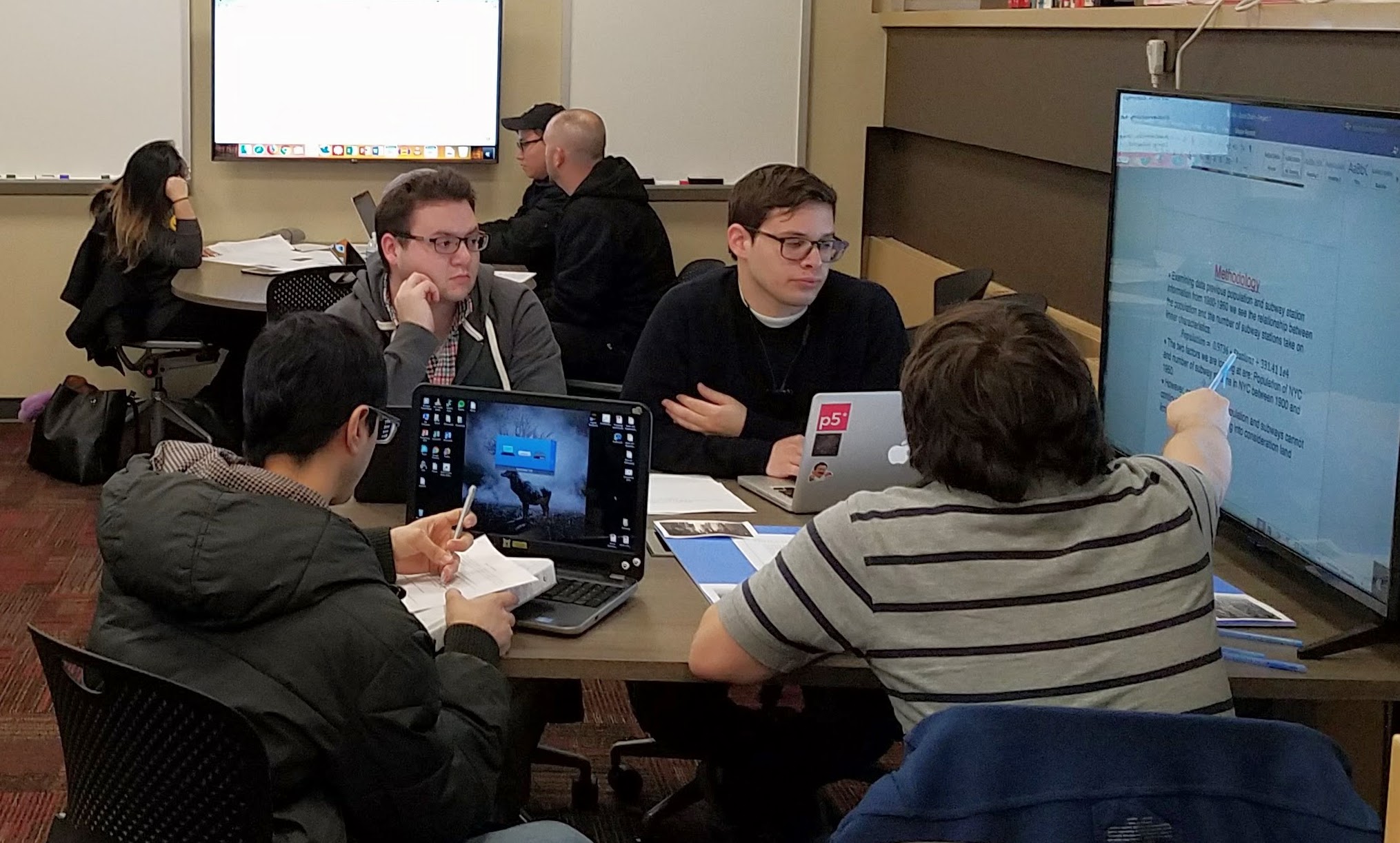}
	
	\begin{minipage}{0.3\textwidth}\caption{In class peer feedback session.}\end{minipage}
\end{wrapfigure}

The exploration, research, and prototyping that happens in the develop stage then needs to be harnessed into a final project. As the tests continue, the students give each other feedback until the finalized model is complete. A peer community is extremely beneficial during this process for giving advice about the Mathematica commands they are using and to help shape the final piece. The class period before the project is due is spent doing a formal peer review. The instructor provides pairs of students with a collection of questions that each student answers about the other student's work. The questions guide students to provide constructive feedback about whether the current model, Mathematica notebook, and report satisfy the project expectations. 

After instructor feedback on the grading criteria (see Section~\ref{sec:grading}), students revise their models one more time before sending the final, honed models to be 3D printed through Shapeways to take advantage of the variety of materials in which they can print, including high resolution nylon and resin, steel, cast metals, and full color sandstone. 

\subsection{The Culminating Experience}

At the end of the semester, the instructor organizes an art exhibition for students to display their artwork in the Queens College Library. This culminating experience of the semester gives the students an opportunity to take a victory lap and have yet another experience from the art world. The students bring their final 3D printed sculptures and we gather together for every student to give a two-minute artist talk about their piece and the highs and lows of the creation process. It is clear that throughout the process, students have developed agency over their artwork and are proud of the results. 

\begin{figure}[h]
	\includegraphics[width=0.7\textwidth]{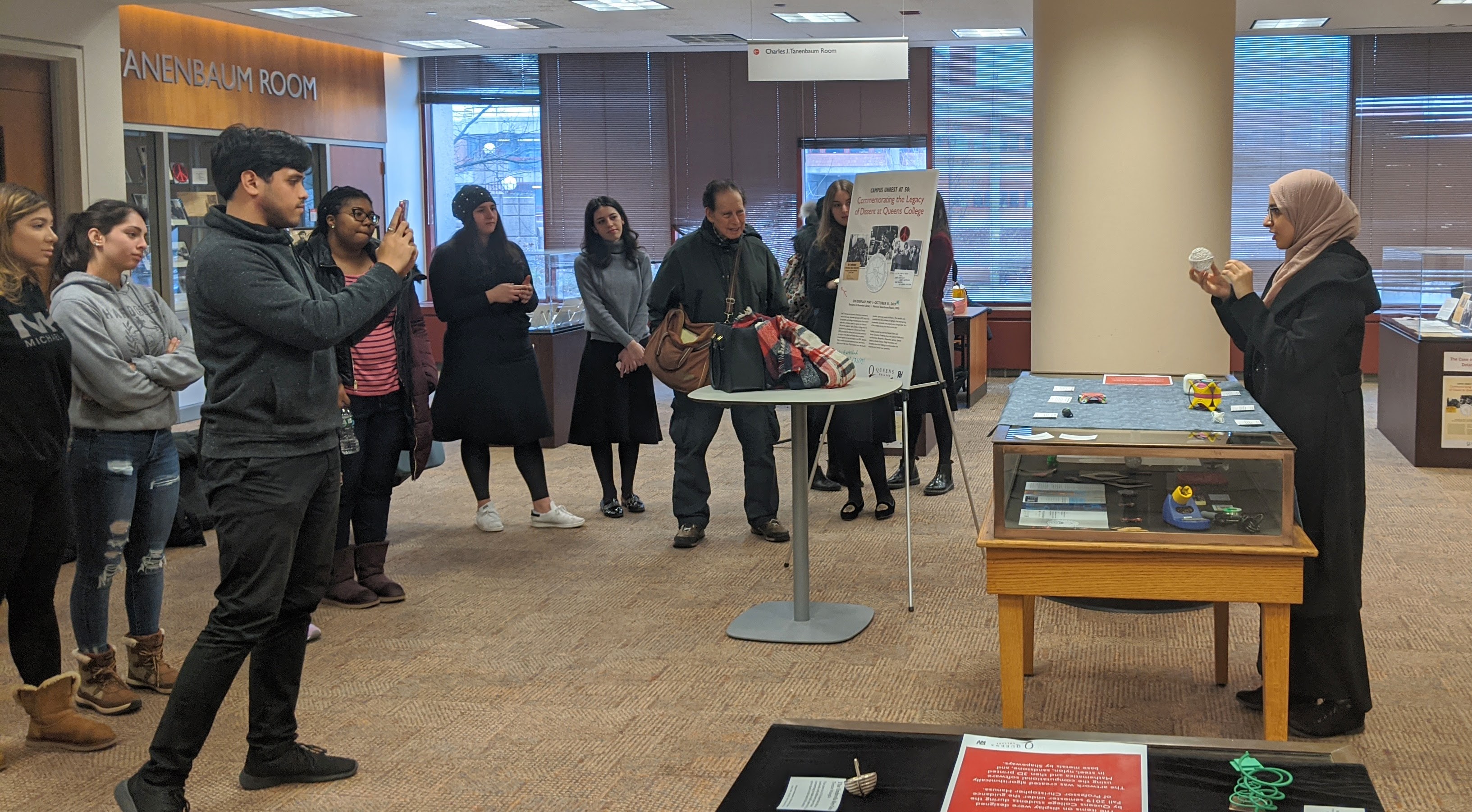}
	\caption{A student presenting work during the art exhibition.}
\end{figure}

\begin{figure}[h]
	\includegraphics[height=1.65in]{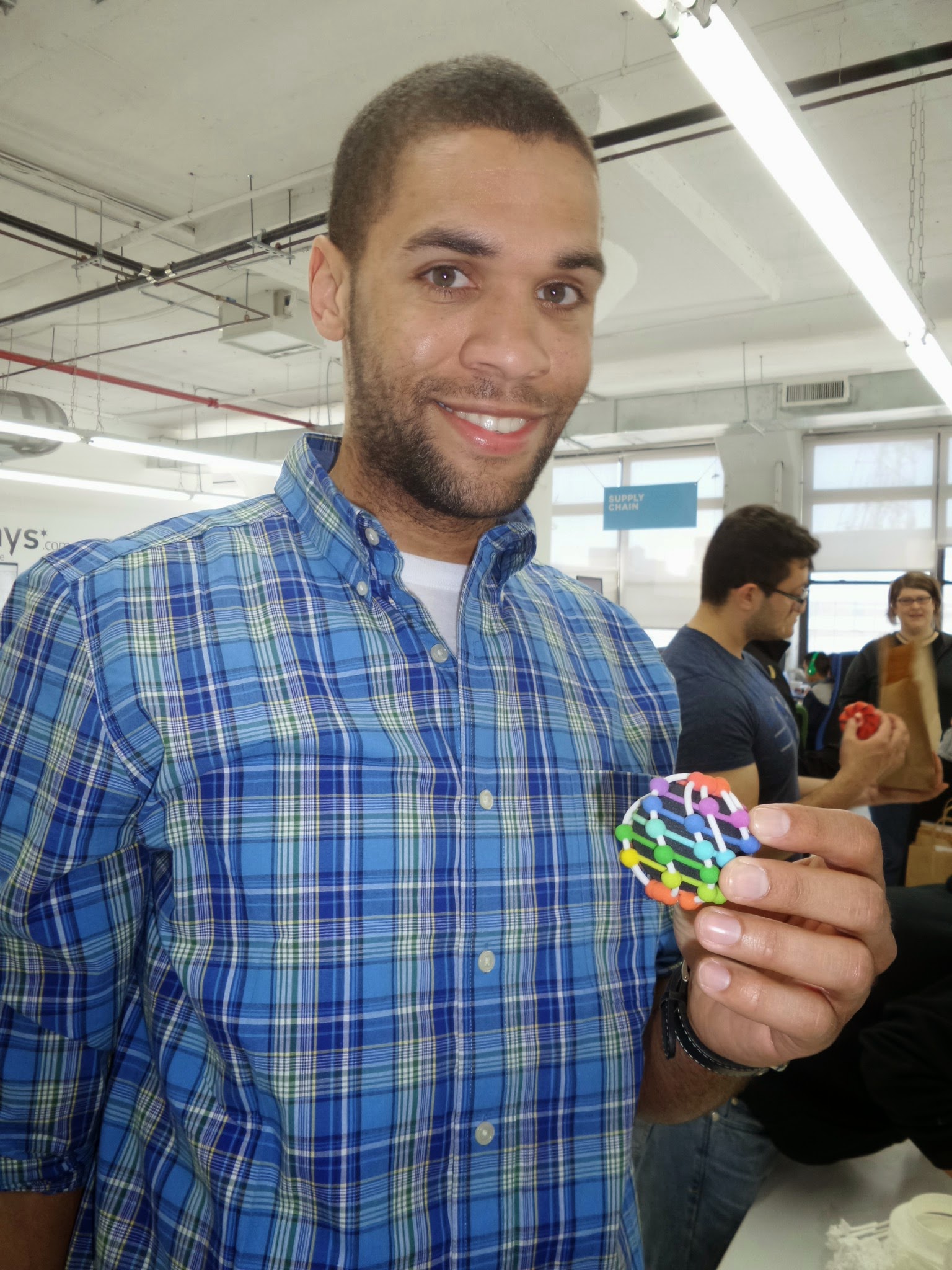}\quad
	\includegraphics[height=1.65in]{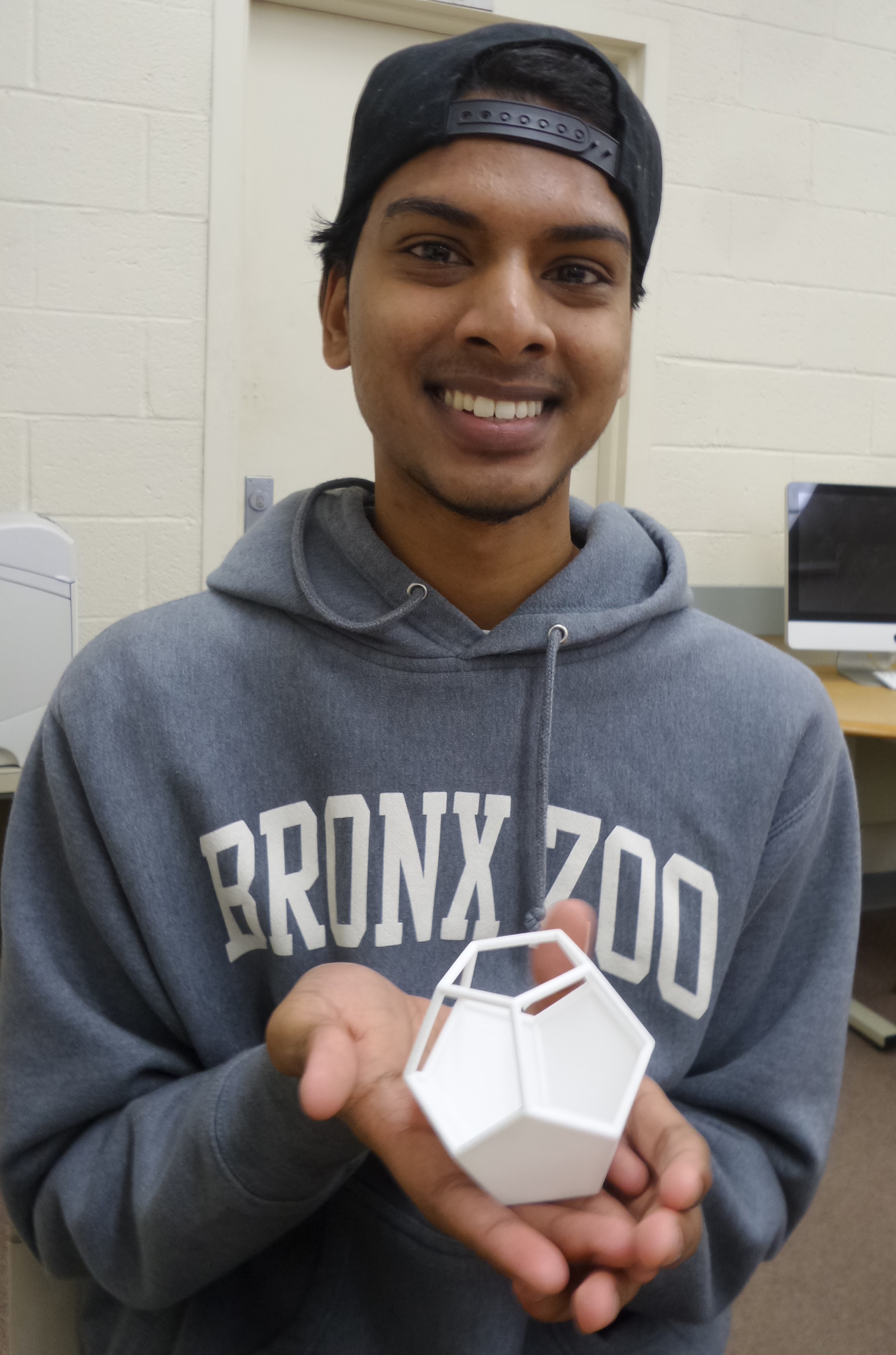}\quad
	\includegraphics[height=1.65in]{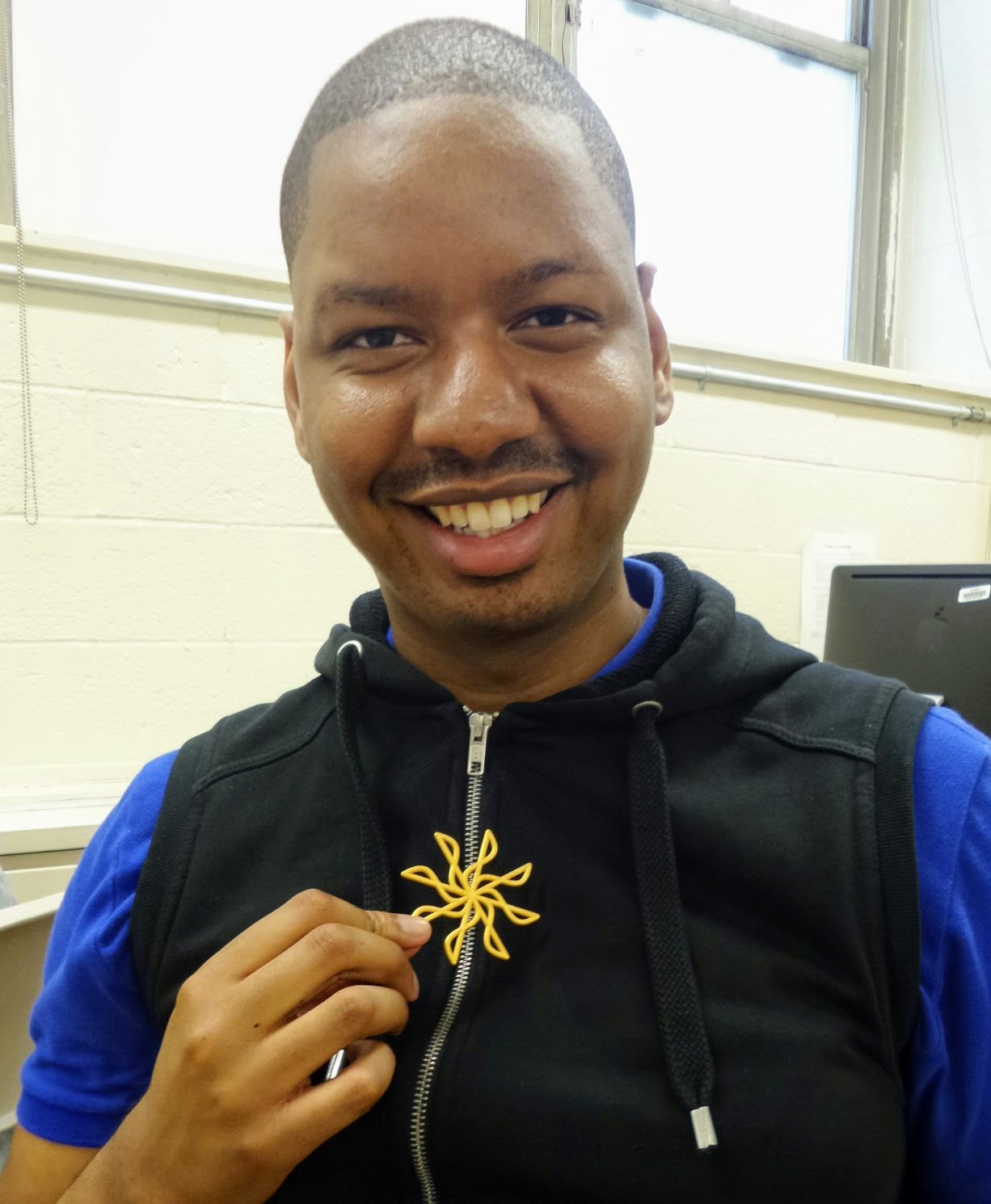}
	\caption{Students proudly showing off their final artwork, 3D printed in full color sandstone, white nylon, and gold steel.}
\end{figure}

\begin{wrapfigure}[21]{R}{0.38\textwidth}
	\centering
	\includegraphics[origin=c,width=0.34\textwidth]{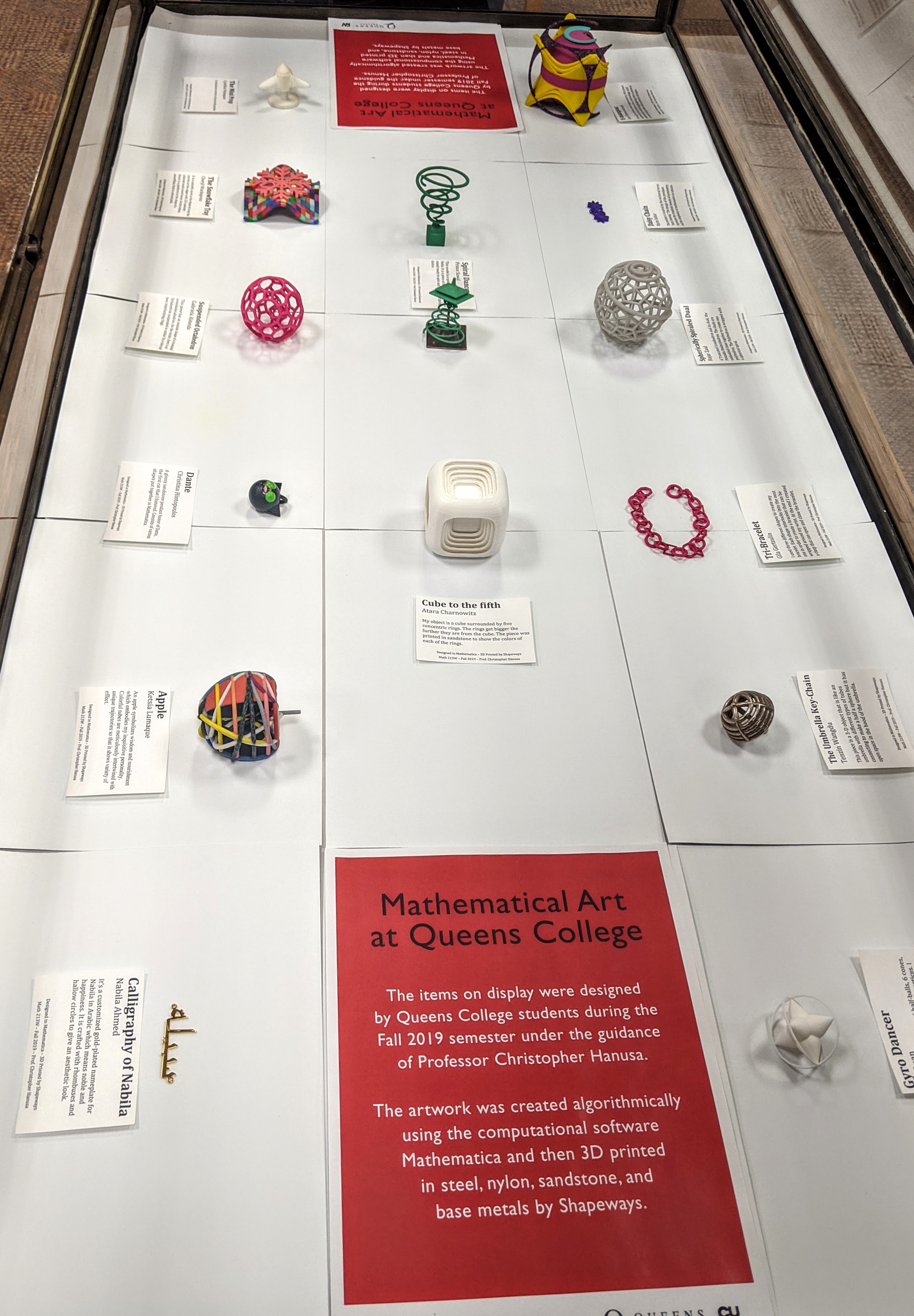}
	
	\begin{minipage}{0.35\textwidth}\caption{An installation of student art at the library.}\end{minipage}
\end{wrapfigure}

At the end of the art show, we install the artwork in an art exhibit that stays on display in the library until the end of the next semester. Before the show the instructor requests from the students the name of their artwork and provide a two-sentence description. This information is used to create professional-looking museum placards that are placed in the display next to the artwork. It is empowering to see one's work in an art exhibit alongside other people's creations. Students are always impressed that they have their work in an art exhibit and enjoy being able to show off their artwork to their friends.

The helpfulness of the Queens College Library staff has made this process a pleasure. It is wonderful to be able to hold the event in the library and for the artwork to be displayed in a public space. The exhibits have become a highlight of campus tours for prospective Queens College students.

\subsection{Project Timeline}

When creating a schedule for this module, it is important to account for the time that it takes to prepare a 3D model for printing, have it printed by the 3D printing company, and receive the model after shipping. This explains why this module is the second of three, and not the final one. Here is a representative timeline for organizing the content in a course that meets twice a week for 14 weeks:

\begin{description}

	\item[Day 9] The class wraps up Module 1. Students complete a tutorial on the basics of 2D graphics. Students sign up for a tour of the Queens College Makerspace and create a Shapeways for Education account.

	\item[Day 10] Students complete a tutorial on the basics of 3D graphics with an emphasis on functional programming. Students are given a list of sites with mathematical art and asked to share and annotate a couple examples on our online discussion board.

	\item[Day 11] Students complete a tutorial on more advanced 3D graphics. On this day Professor Greco talks about 3D printing and the artistic elements of sculpture. 

	\item[Day 12] Students complete a tutorial on making high quality 3D models and start to brainstorm possible project ideas. We discuss how 3D printing works to ensure the final models will be printable.

	\item[Day 13] Students complete a tutorial on 3D design techniques. Students finalize their project specifications. 

	\item[Day 14] Students see a number of minimal working examples to see how to realize common 3D design tasks and watch a video about advanced techniques in 3D design. Students are actively designing their sculpture inside and outside class.

	\item[Day 15] Students are in full project work mode inside and outside class. In class we introduce Ultimaker Cura and rapid prototyping, getting practice using the Makerspace to print a model.

	\item[Day 16] Students finalize their model and print a physical prototype for peer review. They complete an initial draft of their lab report (minus discussion of revision process).
	
	\item[Day 17] Students complete an in-class peer review of prototype and draft report, after which they revise their sculpture, notebook, and report.
	
	\item[Day 18] Students submit final prototype to Shapeways for 3D printing in desired material, and submit their final report to the instructor. The class transitions to thinking about Module 3.
	
\end{description}

Soon after Day 18, the instructor gives feedback on the final report and allows for revisions. Shapeways prints and ships the 3D printed model within the next few weeks. The semester culmination event and art exhibit installation take place during our assigned final exam period.

\section{Project Deliverables and Grading}

At the end of the project, students submit their physical mathematical sculpture, an organized and commented Mathematica notebook that recreates their 3D model, the digital model computer file, and a two-to-three page writeup about the project.

\subsection{The Written Lab Report}

Students are required to write and submit a two-to-three page lab report that is a reflection on the artistic and creative process and documents the student's process throughout the project.

The goal is to give students the opportunity to reflect on how they have developed in their roles as programmers, mathematicians, and artists. The report also serves to ensure that the final artwork is a product of deliberate choices and not simple happenstance or academic dishonesty.

From a programming standpoint, students are asked to share the commands and algorithmic techniques they integrate into the coding of their projects and how they advance their programming skills. They also detail how they took the initiative to learn the necessary commands and techniques in Mathematica. Students also need to describe the mathematical principles that are at work, and how those principles are evident in the finished sculpture. Further students address how the mathematical principles match their level of mathematical knowledge. 

Students are called upon to reflect upon their artistic choices: the inspiration for the project, the aesthetic they were going for, the elements of art they worked to achieve, the reason for their choice of final material, and what the viewer should appreciate about the artwork.

Students are asked to address prompts related to the creation, prototyping, and revision processes: what sort of obstacles they encountered along the way, how they experienced the first prototype they created and how that inspired them to revise the project. Students have to explain how their project evolved over time, including how discussions with their classmates and during the peer review process influenced the ways in which they changed their project. Last, students are to reflect on how they might approach the process differently if they were starting over.  

\subsection{Grading the Projects}
\label{sec:grading}

In the author's experience, instructors who exclusively grade mathematics exercises take comfort in the idea that their grading role is objective---a mathematical answer is either right or wrong. (The author disagrees with this assessment.) This can lead to discomfort about one's role in grading student artwork because the universe of possible subjects is infinite and the universe of approaches to realizing the final piece is also infinite. In effect, we are called upon to render a value judgment on student work and at the same time we recognize that our interpretation of all parts of the process is subjective. It is also important to recognize that many math students are also often uncomfortable with the freedom (and therefore the responsibility) that comes from an exercise that does not have one ``right'' answer.

The author has developed a grading method for the projects in this course that is inspired by his experience in standards-based grading, which he now uses when teaching a non-project-based class like calculus. In standards-based grading, the learning objectives are transparently provided to the students, and regular assessments and reassessments allow students to show their mastery of the material at any point during the semester. Such an approach makes it possible to assess {\it everything} that the instructor values, including so-called ``soft skills'' like the ability to discuss math with peers. This method of grading has been transformative because it changes the dynamic in the class from adversarial to collaborative between students and instructor, and grades become much more aligned with student progress towards learning objectives. 

The rubric developed for the 3D design project has been refined over the years to address the course learning objectives in Section~\ref{sec:objectives}. Students are assessed on the following standards.
\begin{description}
\item[Attentiveness] The student has made active and steady progress and met deadlines.
\item[Name and Description] The name and description of the final artwork are precise and concise.
\item[Intentionality] The vision of the artwork has been honed over time and has a story behind it.
\item[Mathiness] The object is created using mathematical content at the student's level of knowledge.
\item[Functional Techniques] The object was developed using Mathematica techniques learned in class.
\item[Notebook Organization] The notebook was well organized into sections and with text cells that provide context for the code. 
\item[Artistic Writeup] The report discusses the student's development as an artist and artistic choices.
\item[Technical Writeup] The report discusses the student's technical choices and their development as a both a mathematician and programmer.
\item[Revision Writeup] The report discusses the creation, prototyping, and revision processes.
\item[Writeup style] The report is well written and follows the format requirements.
\end{description}

After the student submits their deliverables, each standard is scored on an E-M-R-N scale (Exemplary, Meets Expectations, Revisions Needed, and Not Assessable) and the instructor provides detailed feedback about each standard. The ten scores are converted to a provisional letter grade as follows.

\medskip
\begin{tabular}{cp{4in}}
	{\bf A} &	Earn a score of M or higher on all standards and a score of E on at least eight standards.\\
	{\bf B} &	Earn a score of M or higher on all standards and a score of E on at least four standards.\\
	{\bf C} &	Earn a score of M or higher on nine standards and no N scores.\\
	{\bf D} &	Earn a score of M or higher on eight standards and at most one N score.\\
	{\bf F} &	Have fewer than eight E or M scores.\\
\end{tabular}

\medskip
Note that a student cannot receive a grade higher than C if there is some aspect of their project that needs revision. Students are given the opportunity to improve their provisional grade by revising their project deliverables and having the instructor re-score the standards. Due to time constraints and instructor mental load, students are limited to one resubmission.

Consistent with Talbert's pillars of alternative grading (see \cite{Talbert}), this grading scheme was developed to provide complete transparency in the expectations for students and how their work will be graded, and by giving students the opportunity for peer and instructor feedback throughout the process with the ability to revise work without penalty. 

Students appreciate the transparency in evaluation expectations. Since the expectations align with the class's learning objectives, a student who addresses the items in the rubric is showing that they learned what the instructor wanted them to learn. The revision process benefits both the student and the instructor. The instructor can convey a frank assessment of the student's work and how it meets or does not meet expectations. The student gets a chance to see where their work does not meet expectations and other aspects that could be improved, which reinforces that project development is a process and there is always room for improvement. The student can decide whether it is worth putting in the time to revise their work to improve their standards scores, and in a direct and transparent way, their final project grade.

\section{Course Reflection}

The students and the author enjoy this project immensely; this section includes reflections from both points of view. 

\subsection{Student Comments}

Student voices are collected at the end of the semester through a request to write ``a letter to future students'', a curated collection of which is shared with incoming students the following semester. Here are some of their words of wisdom.

\smallskip
\begin{myquotation}
``For a lot of students, this class might be one of the first instances that unravels for them what it is exactly that is so beautiful about mathematics, though things like precision are not particularly stressed, your creative freedom will know no bounds in this class. You should absolutely take advantage of this, the depth at which you're allowed to pursue a particular problem is seldom allowed in other mathematics courses, this versatility is a gift.''
\end{myquotation}

\smallskip
\begin{myquotation}
``Incoming students can expect to be as creative as they want when doing their projects.  The professor allows complete creative freedom and that is what made the project fun to complete.  It would have been nice to know before hand that I would need to spend the extra hours outside of the classroom to understand the topics.  Overall, I have nothing but positive things to say about this class and I would recommend it to anyone who is creative and who wish to show their creative prowess.''
\end{myquotation}

\smallskip
\begin{myquotation}
``The amazing thing is that the professor gives you the tool set to do these projects throughout the course and all you have to do is apply that and he even gives you the freedom to put your own flavor in to the mix. With this you get to truly see the diversity that everyone can come up with through each project.''
\end{myquotation}

\smallskip
\begin{myquotation}
	``We completed some major projects. The most exciting of them all was designing three dimensional objects to be sent to a 3D-printing company (Shapeways). I was able to make a piece of jewelry based on the intersection of sine waves and lines that I can sell if I wanted to! There's nothing like holding in your hands and feeling the math  you've been drawing on paper or seeing on the computer for your entire life. We've built tutorials. We've built games. We've built economic models. We've done so much with the robust power of Mathematica. This class inspired me so much that I want to keep playing with it and exploring it even into the Summer, long after the course is over.''
\end{myquotation}

\smallskip
\begin{myquotation}
``Do all of the assignments. Play with Mathematica in your spare time. Explore the internet---there are tons of resources, tutorials and communities of users who are already asking questions you probably are searching for the answer to. This is one of the few classes where you'll be challenged to think mathematically, critically, and creatively. This is one of the few classes where you may not know the answer and neither does the professor. This is one of the few classes where you'll learn that it's okay not to have the answer, even in a mathematics class. The professor will work hardest to help you master Mathematica and bring your dreams to live via coding. Hopefully Mathematica will take you to a new place in your journey with mathematics.''
\end{myquotation}

\smallskip
\begin{myquotation}
``Out of everything I've learned,  there's one thing that stuck out to me the most. One day I asked the professor a mathematical question which for me was hard, so I asked him for help. He thought about my question and then told me `I don't know.' As bright and intelligent as the professor is, he humanized himself in that moment. It made me feel like I have a shot at this math thing. It was more inspirational than anything. He could have given me an answer to make himself look great and leave me baffled but he didn't.''
\end{myquotation}

\subsection{Instructor Reflection}

It is always a treat to teach this class. While the project descriptions and expectations have converged over time, the course is always different and exciting to teach! Project deliverables vary greatly from student to student and semester to semester because I encourage students to draw from their lived experience and personal interests and because Mathematica's computational power extends to a vast variety of subject areas.

This course is a lot of work to coordinate and supervise. The students are given leeway to explore, experience productive failure, and take ownership of their learning. This is a messy process and often the first time the students encounter this type of learning environment, so as an instructor I need to encourage students and provide technical help if their project choices end up being a bit too ambitious. This means that in its current incarnation this course cannot realistically scale up to serve more students without additional teaching assistants or peer mentors. Relatedly, because I prioritize ensuring that the day-to-day process runs smoothly, I end up falling behind in returning project feedback promptly.

One conscious choice I made in the course construction is that there are three projects in the semester. I think the first project (the Mathematica tutorial) serves its purpose phenomenally---students learn the structure of Mathematica and its basic data types, they explore a topic of interest, they learn how to format a notebook, and they ease into the self-directed and collaborative modes that make the rest of the semester run smoothly. It also introduces the alternative grading scheme and reinforces that learning and personal growth is valued. 
Students often wish that they had more time for each project, which is impossible if there are three projects. I take the view that the second project is there to give students a taste of what 3D printing can be. If I devoted a whole semester to it, I am sure that the students' skills in 3D printing would improve as well as the quality of the final pieces. However, we would run into issues at the end of the semester with the time it takes for models to be printed and be shipped, and the students would miss out on the distinct yet complementary third project where they develop different programming skills and learn about designing interfaces with human interaction. Often the second and third projects are both rated equally positively by students in terms of preference and how much they felt they learned.

One of the most rewarding points in the class is at the end of the second project. It is such an amazing and empowering experience to hold something in your hand that had only previously existed in your mind and in the computer. When the final models arrive, students can't wait to show me and their classmates; they regularly relate how proud they were to share their sculpture with their friends and family. And I often get feedback that before the class they never expected to be able to create something of their own, yet after the class feel comfortable and confident with that process. With the opening of the Queens College Makerspace, I expect future students to take even more ownership of their learning. 

Just as I hope that my students develop their growth mindset in this course, I am continually working to improve my teaching and the curriculum. For example, each semester I modify the standards and grade conversion table as I learn more about the way I assign scores and what I expect from students. In addition, one thing I noticed that had been lacking in previous semesters was that students would proceed linearly from start to finish without exploring a wide space of possible solutions. So I added a requirement in Fall 2021 that students need to modify some parameter or aspect of their project, and reflect on which version they prefer and why. With this prompt, students thought more deeply about why they ended up with their final result and they shared more about their creative process in their reports. 

Mathematical Computing is a course destined for math majors who have a certain amount of mathematical background. I wanted to be able to offer the opportunity to be creative with mathematics and create mathematical art to non-majors. I created a project-based course for non-majors called Mathematical Design with the online graphing calculator Desmos as its medium. I look forward to further developing the curriculum for this class and publishing more about it. Stay tuned!

\section{Acknowledgments}

Many thanks to Leah Wrenn Berman, Emily Dennett, and Patrick Johnson for their advice for improving this article. I would also like to thank Matthew Greco for being a great collaborative partner at Queens College. 

\newpage
\appendix
\section{3D design in Mathematica}
\label{sec:3DDesign}

In this section the reader can learn multiple techniques to use Mathematica for 3D design. An interactive copy of this notebook is available on the Wolfram Notebook Archive \cite{HanusaMMA}. Additional resources about 3D design in Mathematica are available on the author's webpage\footnote{\url{http://qc.edu/~chanusa/mathematica}}. 

The author would like to acknowledge Mathematica's Documentation Center \cite{Mathematica} for its extensive examples that explain the Wolfram Language in great detail, the Mathematica StackExchange community\footnote{\url{https://mathematica.stackexchange.com/}} for their selfless sharing of ideas that have inspired a number of advanced explorations, and Henry Segerman's article \cite{Segerman}, which was very helpful with initial trials and tribulations in 3D modeling in Mathematica. Lastly, many thanks to Jakub Kuczmarski for his {\tt mmaCells} \LaTeX\ package\footnote{\url{https://github.com/jkuczm/mmacells}} that allowed for in-line typesetting of Mathematica code in this article.

\subsection{Basic 3D graphics and exporting}
\label{sec:basic3D}

Creating a 3D model file with Mathematica can be exceptionally easy or extremely frustrating. In many cases, creating your 3D file is but a simple call to the {\tt Export} function, where you specify the file extension {\tt .stl} for uncolored models or {\tt .wrl} for colored models.
\begin{mmaCell}{Code}
Export["file.stl", model]

\end{mmaCell}
In general, I like to append the file name to the current notebook's working directory (using the {\tt NotebookDirectory} command) to organize my files.

The most basic way to create an object in Mathematica is to specify the coordinates and parameters of the 3D primitives such as {\tt Sphere}, {\tt Cylinder}, {\tt Tube}, {\tt Tetrahedron} and wrap them all in a {\tt Graphics3D} command. Happily, starting in Mathematica version 13, these primitives behave well when they are exported. For example, here is the code that builds a model of a basic snowman using spheres, cylinders, a cone, and some rectangular prisms:
\noindent
\begin{mmaCell}[moredefined={Sphere,Cylinder,Cone}]{Code}
snowman = Graphics3D[{
  Sphere[{0, 0, 0}, 1],
  Sphere[{0, 0, 1.3}, 0.8], 
  Sphere[{0, 0, 2.4}, 0.6],
  Black, Cylinder[{{0, 0, 2.8}, {0, 0, 2.9}}, 0.7], 
  Cylinder[{{0, 0, 2.9}, {0, 0, 3.5}}, 0.5],
  Orange, Cone[{{0, -0.55, 2.4}, {0, -0.9, 2.4}}, 0.1], 
  Black, Cuboid[{0, -0.75, 1.3} - 0.1 {1, 1, 1}, 
                {0, -0.75, 1.3} + 0.1 {1, 1, 1}],
  Cuboid[{0, -0.7, 1.6} - 0.1 {1, 1, 1}, 
         {0, -0.7, 1.6} + 0.1 {1, 1, 1}],
  Cuboid[{0, -0.7, 1.0} - 0.1 {1, 1, 1}, 
         {0, -0.7, 1.0} + 0.1 {1, 1, 1}]
}] 
Export[NotebookDirectory[] <> "snowman.stl", snowman]

\end{mmaCell}
In each {\tt Sphere} command, you specify the center of the sphere and its radius. In the {\tt Cylinder} and {\tt Cone} commands you are specifying the endpoints of its axis of rotation as well as its radius. In the {\tt Cuboid} command, you specify the opposite corners of the rectangular prisms. The resulting snowman looks like this:
\mmaCellGraphics[ig={height=3in},pole2=vc,yoffset=.5ex]{Output}{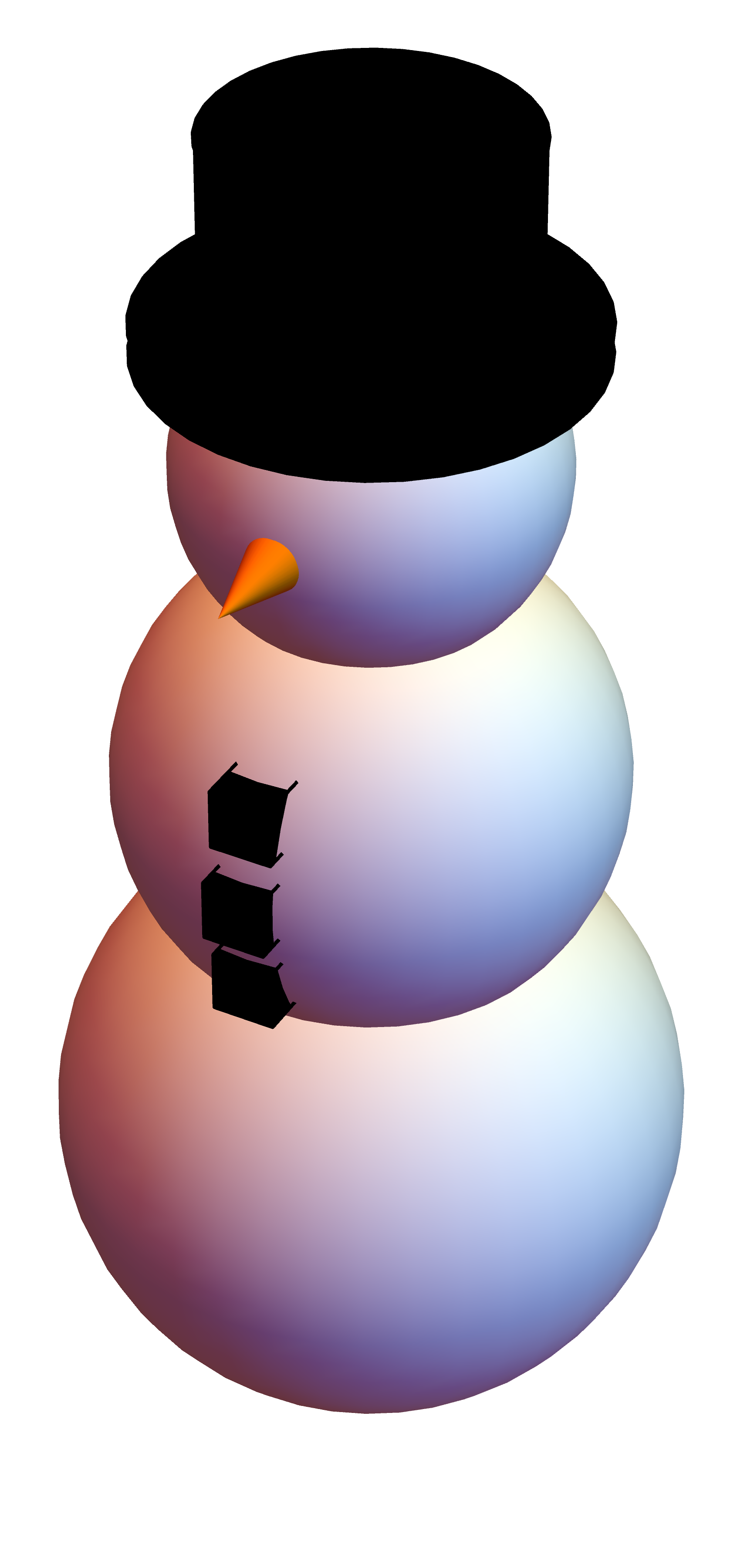}

I always suggest that you re-import your STL files to ensure that everything exported correctly. 

\begin{mmaCell}{Code}
Import[NotebookDirectory[] <> "snowman.stl"]

\end{mmaCell}

\mmaCellGraphics[ig={height=3in},pole2=vc,yoffset=.5ex]{Output}{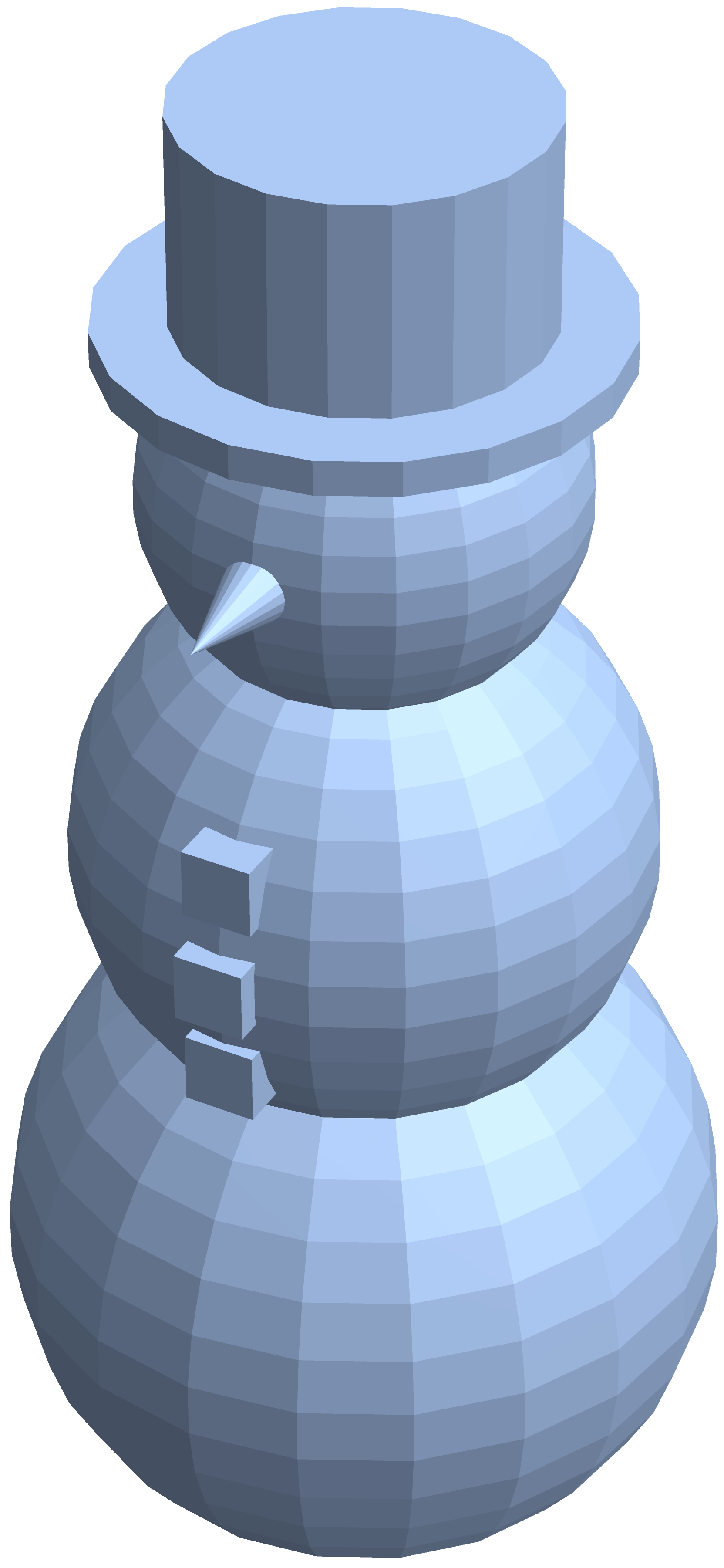}

Upon observation, we see that the imported model has no color and that the final result is very faceted because the model was exported with the default options. We can greatly improve the result with the more advanced commands shown in Section~\ref{sec:discretize}; we use these techniques to write updated code to create a high quality, full color model of the snowman in Section~\ref{sec:snowman}.

Mathematica is not currently capable of importing WRL files so any colored models should be opened in a different program or uploaded to a web service such as Sketchfab\footnote{\url{https://sketchfab.com}} to ensure that the model has exported correctly. 

\subsection{Curves}

The next level of 3D design complexity uses knowledge of multivariable calculus to create models of curves and surfaces. A 3D model of a vector curve $\mathbf{r}(t)=\langle f(t),g(t),h(t)\rangle$ is made using {\tt ParametricPlot3D} with one input variable. For example, a basic helix can be created using the code

\begin{mmaCell}{Code}
ParametricPlot3D[{t/Pi, Cos[t], Sin[t]}, {t, 0, 6 Pi}]

\end{mmaCell}

\mmaCellGraphics[ig={height=1in},pole2=vc,yoffset=.5ex]{Output}{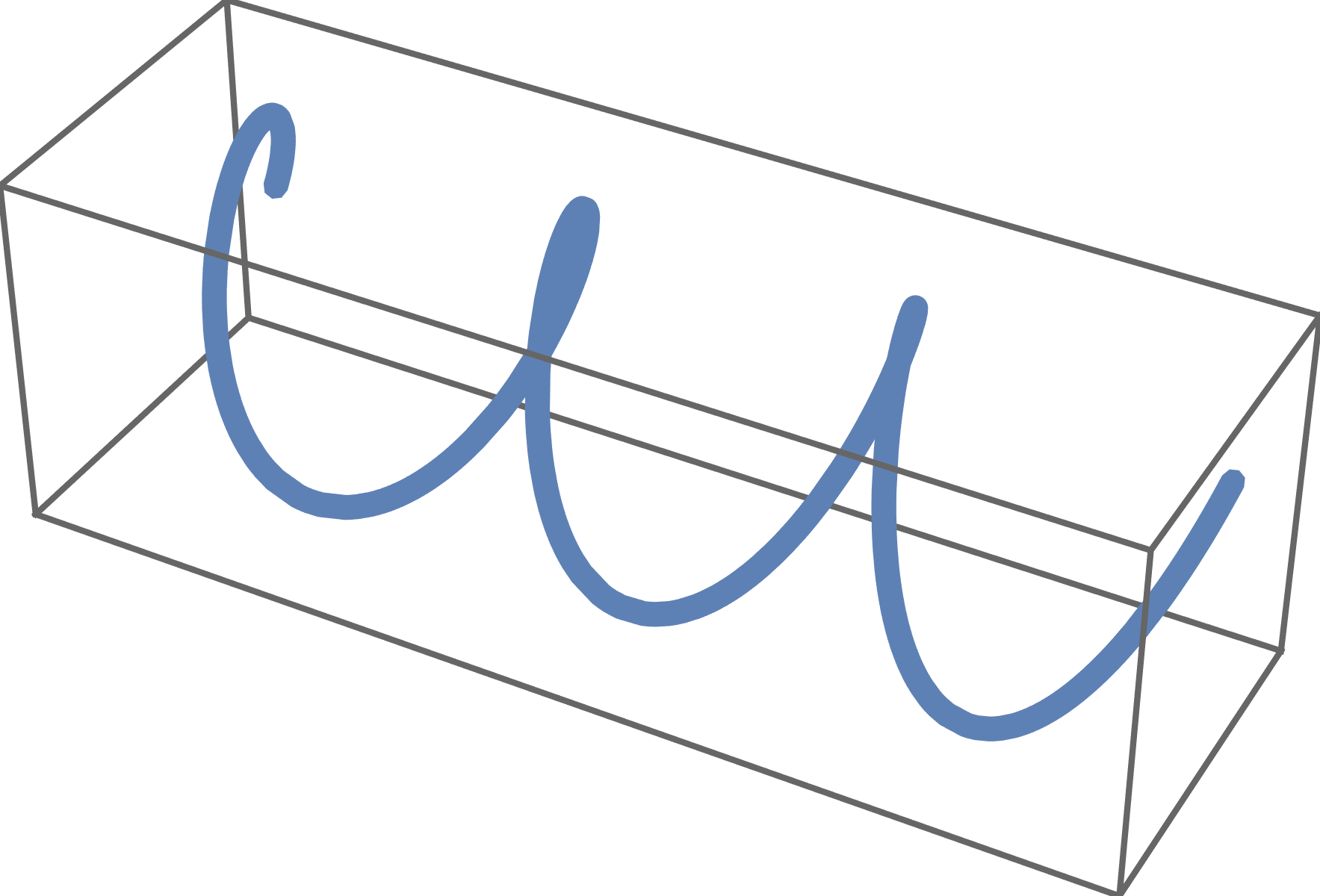}

However, this is where the precision of an idealized mathematical model runs into the reality of what can be 3D printed. A helix is a one-dimensional curve embedded in three dimensions, so we need to add some thickness to create a printable model. We can do this using a {\tt PlotStyle} option that specifies the radius of a tube that will be swept out along the curve. 
\begin{mmaCell}[moredefined={Tube}]{Code}
helix = ParametricPlot3D[{t/Pi, Cos[t], Sin[t]},
  {t, 0, 6 Pi}, PlotStyle -> {Tube[0.1]}, 
  Axes -> False, Boxed -> False, SphericalRegion -> True];
Import[Export[NotebookDirectory[] <> "helix.stl", helix]]

\end{mmaCell}

\mmaCellGraphics[ig={height=1in},pole2=vc,yoffset=.5ex]{Output}{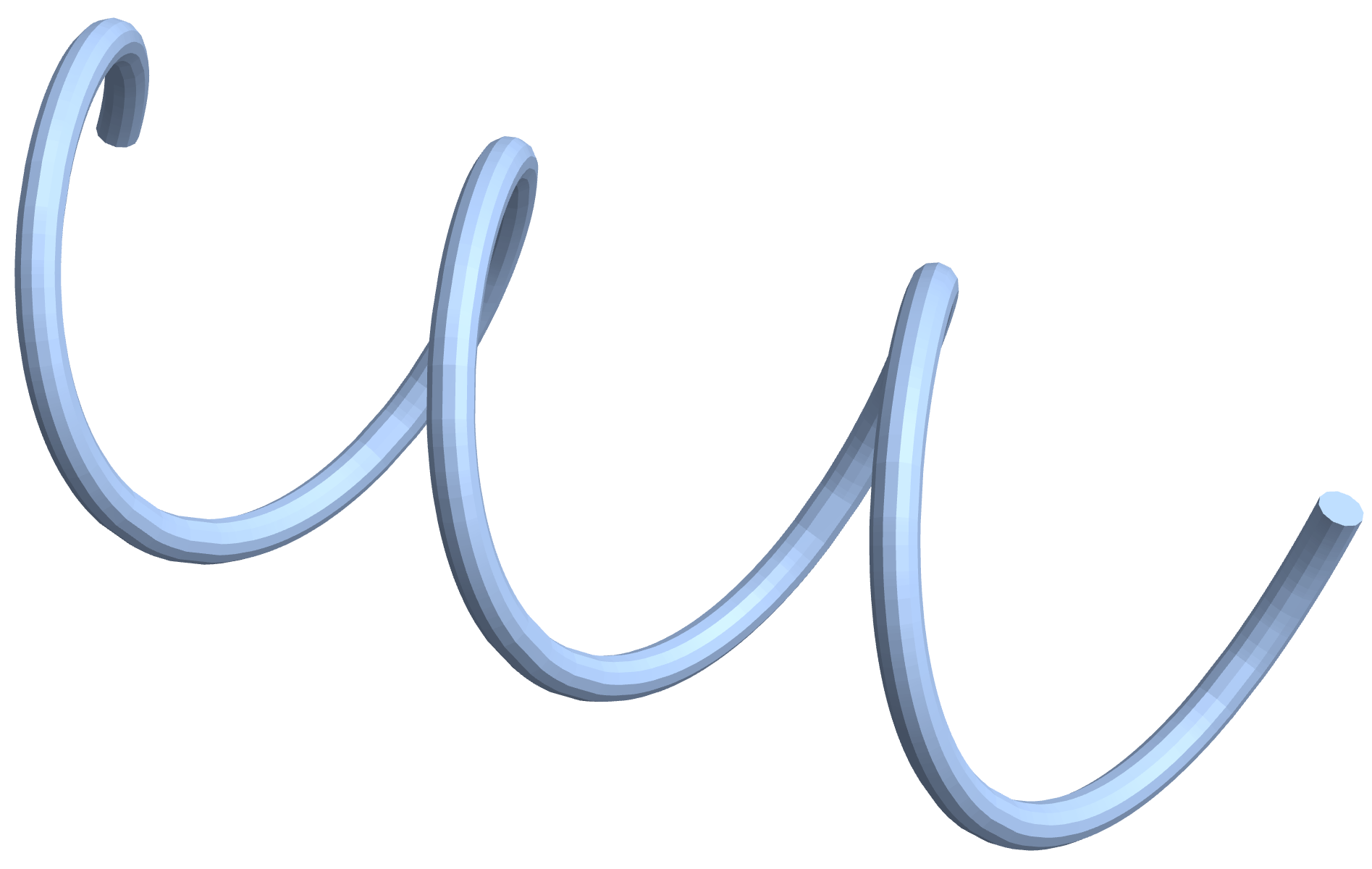}

The {\tt Axes}, {\tt Boxed}, and {\tt SphericalRegion} options improve visualization in the notebook by removing the default framing around the model and ensuring that when you rotate the model it does not resize the cell. 

It is only on re-importing the model do we realize that the result is once again very faceted. We can use the {\tt PlotPoints} option twice to improve the quality of this 3D print. We specify how many points will be sampled along the curve (here 500), as well as how many points will be sampled in the circular cross-section of the tube (here 50).
\begin{mmaCell}[moredefined={Tube}]{Code}
helixHD = ParametricPlot3D[{t/Pi, Cos[t], Sin[t]},
  {t, 0, 6 Pi}, PlotStyle -> {Tube[0.1, PlotPoints -> 50]}, 
  PlotPoints -> 500, 
  Axes -> False, Boxed -> False, SphericalRegion -> True];
Import[Export[
  NotebookDirectory[] <> "helix.HD.stl", helixHD]]
  
\end{mmaCell}
\mmaCellGraphics[ig={height=1in},pole2=vc,yoffset=.5ex]{Output}{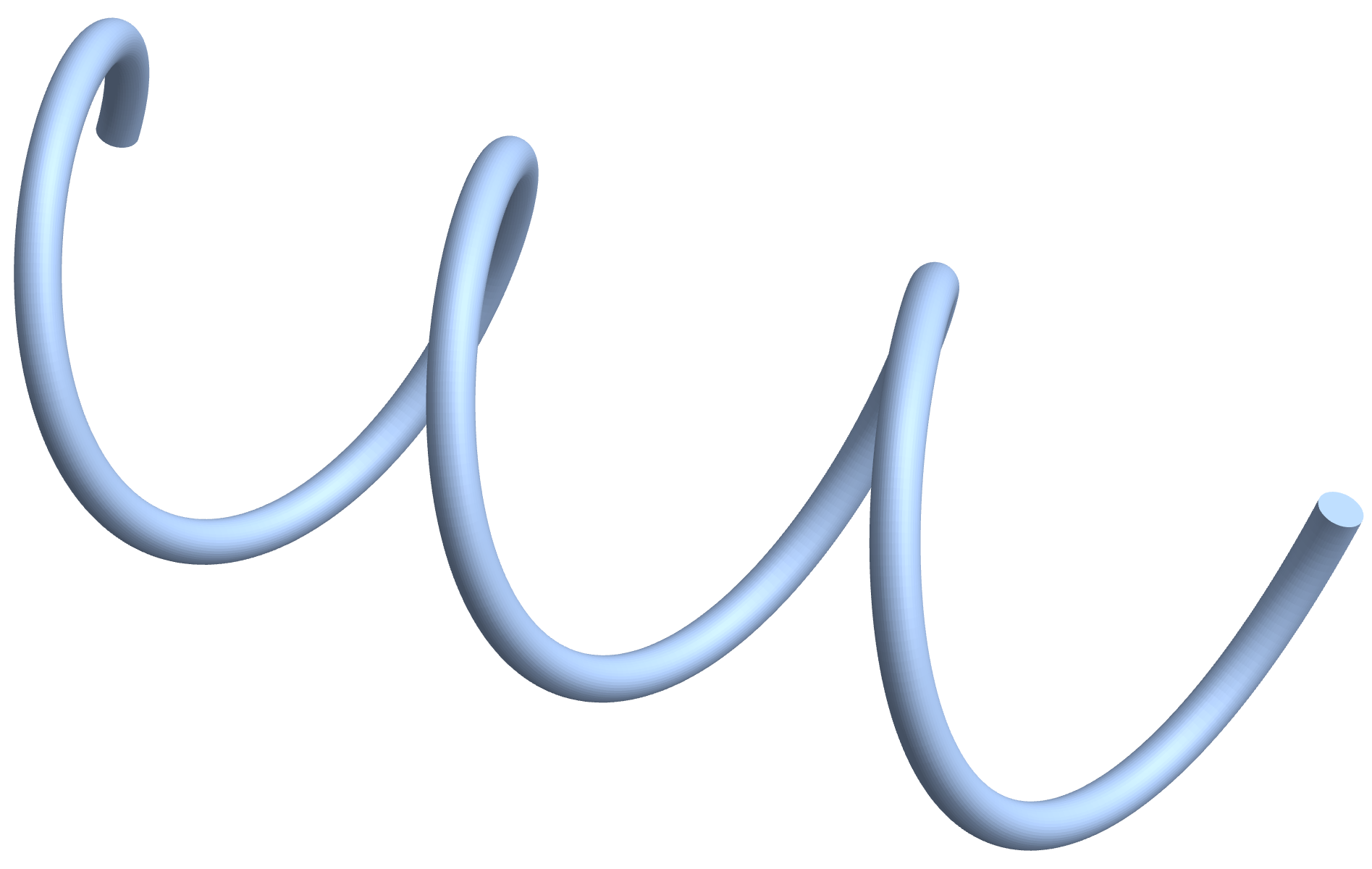}
(In your notebook the first {\tt PlotPoints} option will be {\color{red} red}.) The more points specified, the more detailed the print will be; at the same time, this increases computation time and the size of the exported file. The two plots {\tt helix} and {\tt helixHD} look the same in the notebook; it is only upon re-importing the stl models that we can see that the fidelity of the model has greatly improved.

{\tt ParametricPlot3D} also allows you to 3D print a spline specified by a set of control points. The primitive for a spline is a {\tt BSplineCurve}:
\begin{mmaCell}[moredefined={Sphere,Thick,BSplineCurve,SplineClosed}]{Code}
pts = {{0, 0, 0}, {1, 0, 0}, {0, 1, 0}, {0, 0, 1}};
Graphics3D[{
  Blue, Map[Sphere[#, 0.03] &, pts],
  Thick, Green, Line[Append[pts, pts[[1]]]],
  Red, BSplineCurve[pts, SplineClosed -> True]},
Axes -> True, AxesOrigin -> {0, 0, 0}, Boxed -> False]

\end{mmaCell}
\mmaCellGraphics[ig={height=1.5in},pole2=vc,yoffset=.5ex]{Output}{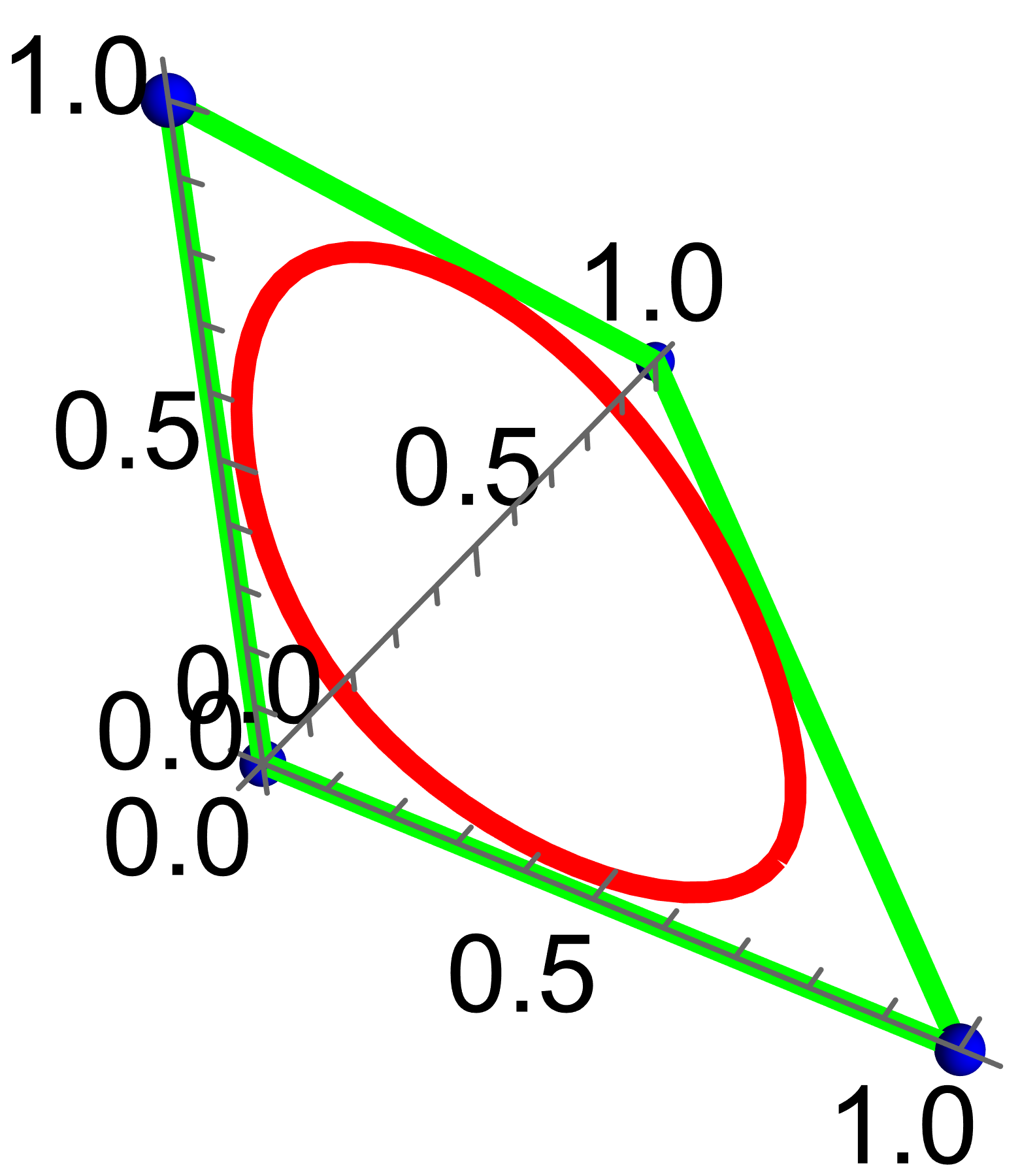}
You can 3D print this curve by first converting the primitive to a parametric function on the domain $[0,1]$ by replacing {\tt BSplineCurve} by {\tt BSplineFunction} and then using {\tt ParametricPlot3D} as above.
\begin{mmaCell}[moredefined={BSplineFunction,SplineClosed,Tube}]{Code}
ParametricPlot3D[ 
  BSplineFunction[pts, SplineClosed -> True][t], 
  {t, 0, 1}, PlotStyle -> {Tube[0.02]}, 
  PlotRange -> All, Axes -> False, Boxed -> False]
  
\end{mmaCell}
\mmaCellGraphics[ig={height=1in},pole2=vc,yoffset=.5ex]{Output}{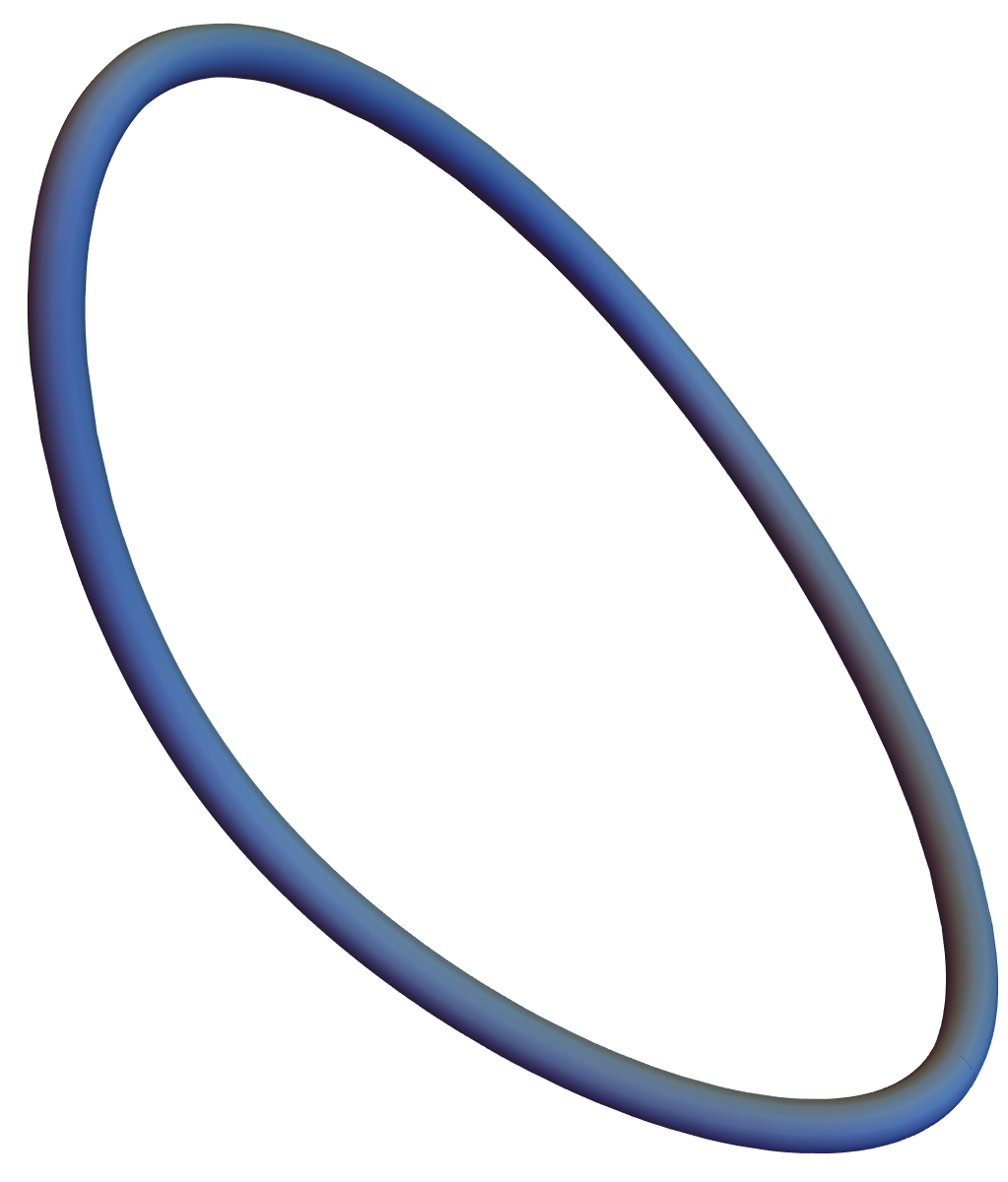}

\subsection{Surfaces}

There are four ways to create two-dimensional surface embedded in three dimensions, depending on how it is defined. 

When the surface is the graph of a function $f(x,y)$ of two variables, use {\tt Plot3D}. You need to specify the two independent variables and their domains. By default {\tt Plot3D} only outputs the 2D surface; thicken it to make it 3D printable by specifying the {\tt PlotTheme} option to be {\tt "ThickSurface"}. 

\begin{mmaCell}[moredefined={PlotTheme}]{Code}
Plot3D[Sin[x + y^2], {x, -3, 3}, {y, -2, 2}, Mesh -> None, 
  PlotPoints -> 30, PlotTheme -> "ThickSurface"]
  
\end{mmaCell}
\mmaCellGraphics[ig={height=1in},pole2=vc,yoffset=.5ex]{Output}{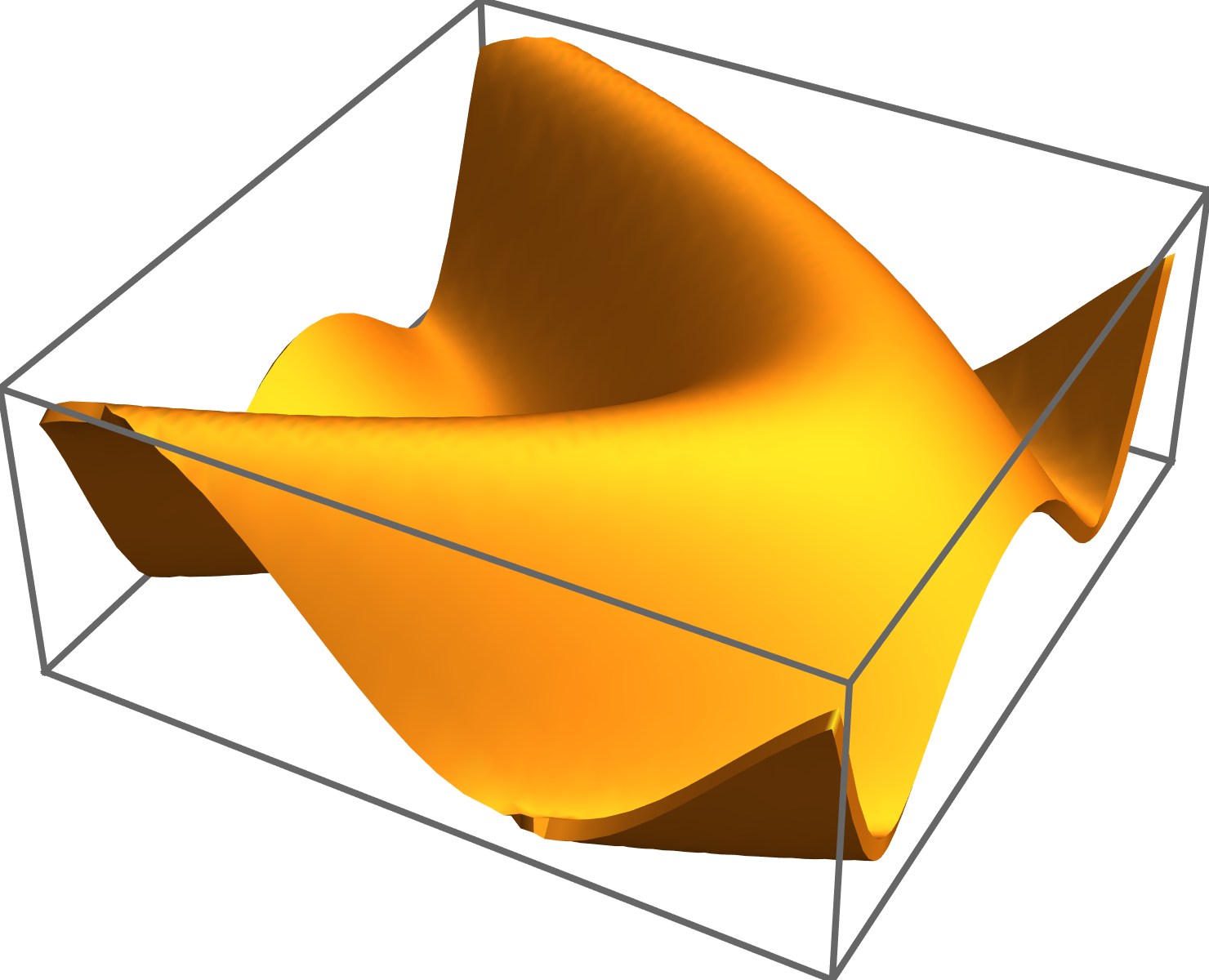}

That option does not allow you to specify a thickness; to do so instead use the hidden {\tt Extrusion} option to thicken the surface in the direction normal to the surface. (This option will show in {\color{red} red} in your notebook.) In the following example we also see the domain does not have to be a rectangle. (In your Mathematica code be sure to replace IN by $\in$ by typing \keys{esc} {\tt $\backslash$in} \keys{esc}.)

\begin{mmaCell}[morelocalconflict={Extrusion}]{Code}
Plot3D[Sin[x + y^2], {x, y} IN Disk[{0, 0}, 2], 
  Mesh -> None, PlotPoints -> 30, Extrusion -> 0.25] 
  
\end{mmaCell}
\mmaCellGraphics[ig={height=1in},pole2=vc,yoffset=.5ex]{Output}{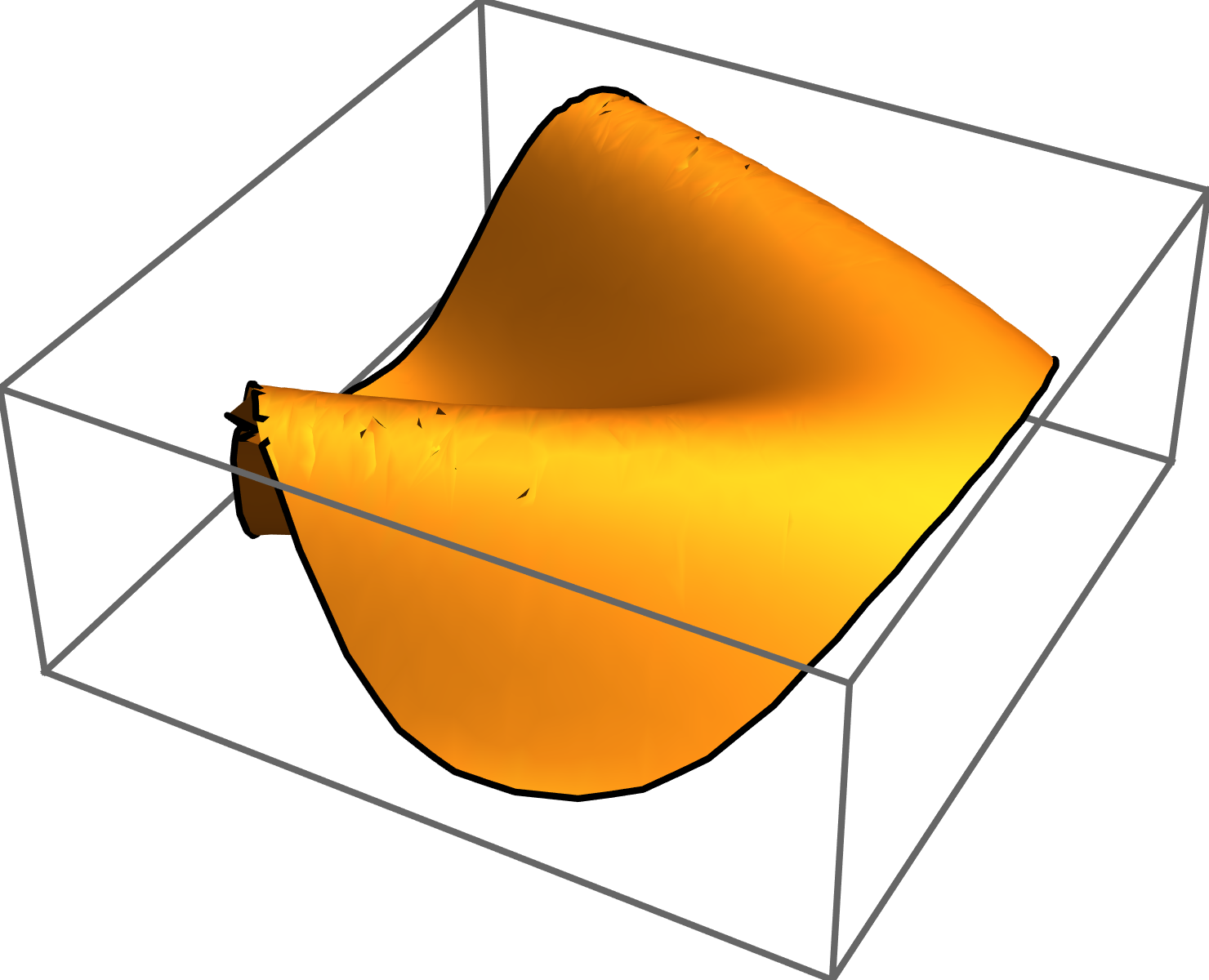}

If instead the surface is defined implicitly as the level surface of a function $f(x,y,z)=k$, use {\tt ContourPlot3D}:

\begin{mmaCell}[moredefined={ContourPlot3D},morelocalconflict={Extrusion}]{Code}
ContourPlot3D[x^3 + y^2 - z^2 == 0, 
  {x, -2, 2}, {y, -2, 2}, {z, -2, 2}, 
  Mesh -> None, PlotPoints -> 30, Extrusion -> 0.2]
  
\end{mmaCell}
\mmaCellGraphics[ig={height=1.5in},pole2=vc,yoffset=.5ex]{Output}{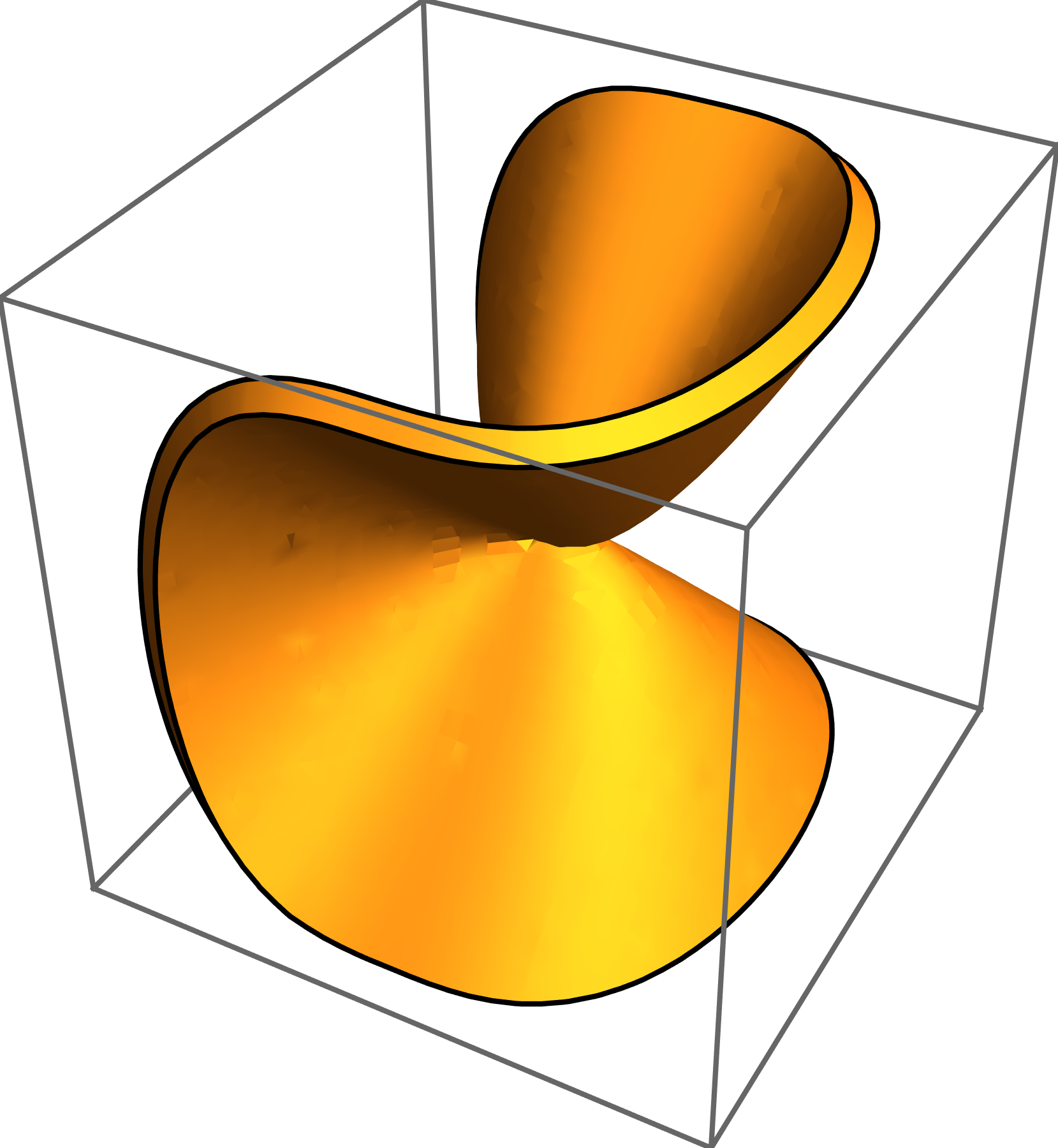}

Otherwise, if the surface is defined parametrically as a function $f(u,v)$ of two parameters, use {\tt ParametricPlot3D} and specify the domain of each parameter. The following example is a torus; we specify the number of sample points for each parameter independently.

\begin{mmaCell}{Code}
ParametricPlot3D[
  {(3 + Cos[v]) Cos[u], (3 + Cos[v]) Sin[u], Sin[v]}, 
  {u, 0, 2 Pi}, {v, 0, 2 Pi}, 
  Mesh -> None, PlotPoints -> {100, 30}]
  
\end{mmaCell}
\mmaCellGraphics[ig={height=1in},pole2=vc,yoffset=.5ex]{Output}{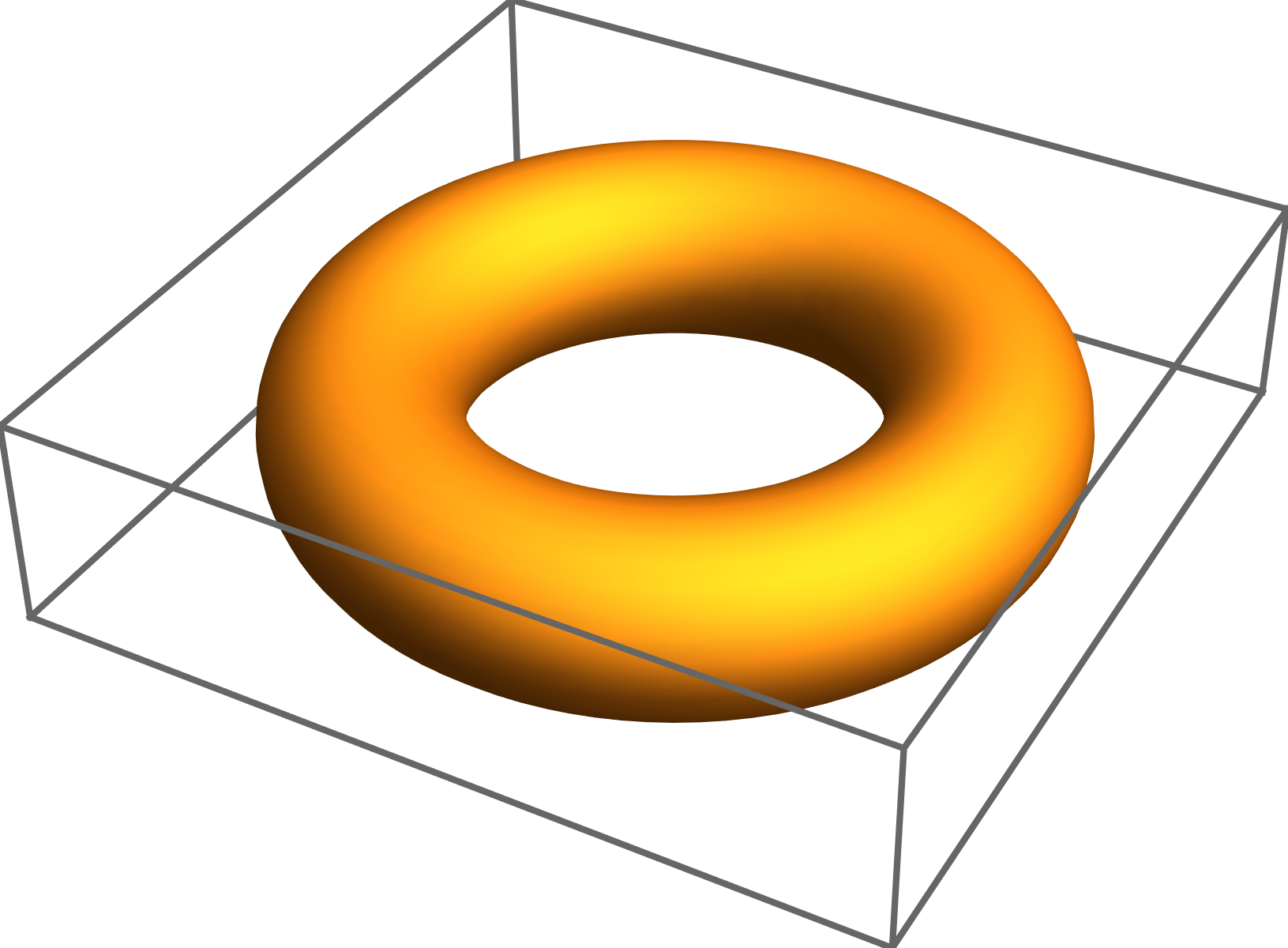}

That torus is the boundary of a 3D solid, so the file exported by Mathematica can be 3D printed. In general, a parametric surface may not be the boundary of a solid so the slicing software for 3D printing will tell you that your object is ``not manifold''. You will need to thicken those surfaces as we did before.

\begin{mmaCell}[morelocalconflict={Extrusion}]{Code}
ParametricPlot3D[
  {v Cos[u], v Sin[u], -Cos[v]}, 
  {u, 0, 2 Pi}, {v, 0, 2 Pi}, 
  Mesh -> None, PlotPoints -> {100, 30}, Extrusion -> .2]
  
\end{mmaCell}
\mmaCellGraphics[ig={height=1in},pole2=vc,yoffset=.5ex]{Output}{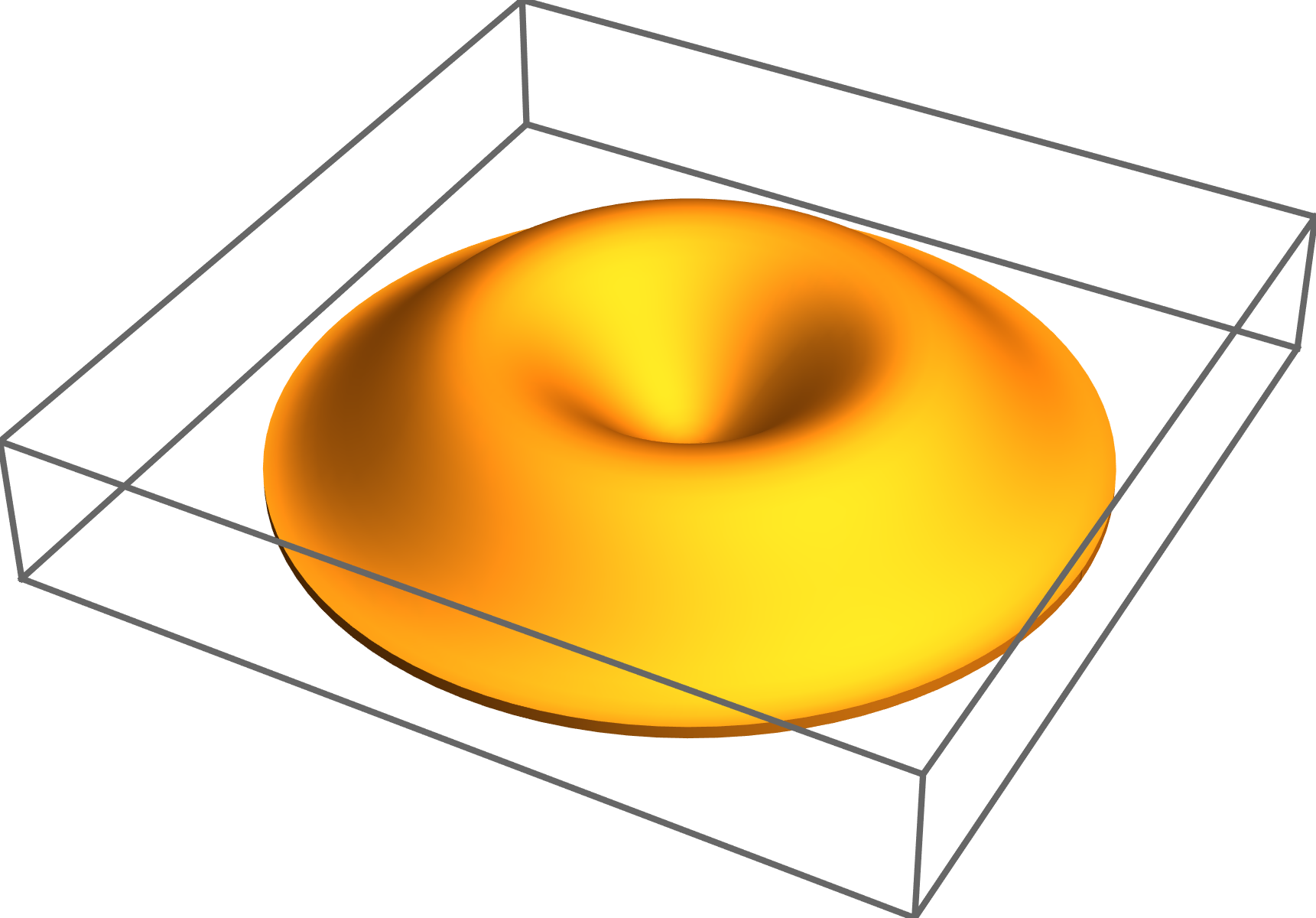}

Last, you can create a two dimensional spline surface using {\tt BSplineSurface}. Here is the example from the Wolfram Documentation Center that generates a random surface over the rectangular domain $[1,5]\times[1,5]$.
\begin{mmaCell}[moredefined={BSplineSurface,RandomReal}]{Code}
pts = Table[{i, j, RandomReal[{-1, 1}]}, {i, 5}, {j, 5}];
Graphics3D[BSplineSurface[pts]]

\end{mmaCell}
\mmaCellGraphics[ig={height=1in},pole2=vc,yoffset=.5ex]{Output}{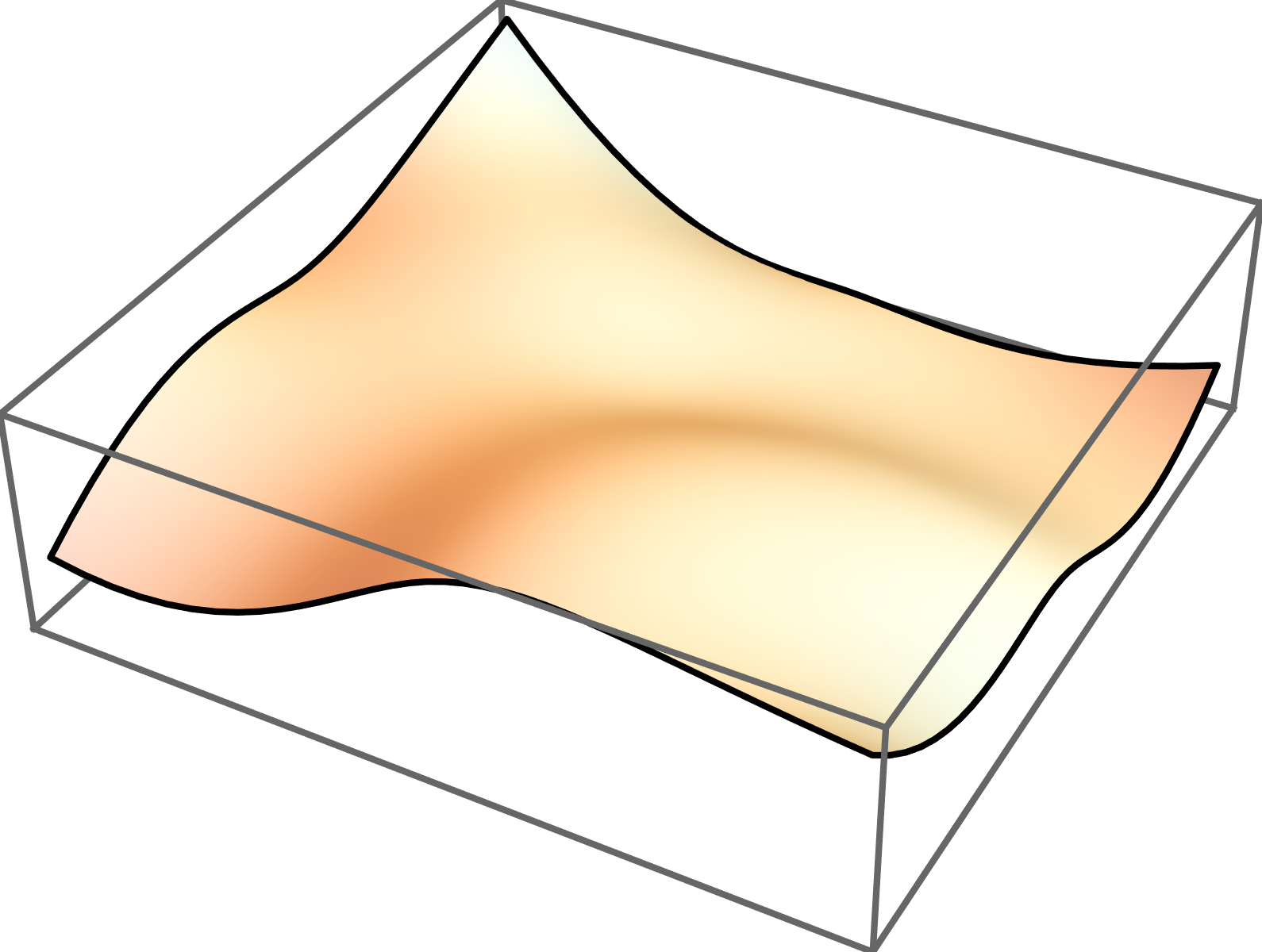}
Again we need to convert the {\tt BSplineSurface} to a {\tt BSplineFunction} to create a 3D printable model.
\begin{mmaCell}[moredefined={BSplineFunction},morelocalconflict={Extrusion}]{Code}
ParametricPlot3D[BSplineFunction[pts][u, v], 
  {u, 0, 1}, {v, 0, 1}, Mesh -> None, Extrusion -> .2]

\end{mmaCell}
\mmaCellGraphics[ig={height=1in},pole2=vc,yoffset=.5ex]{Output}{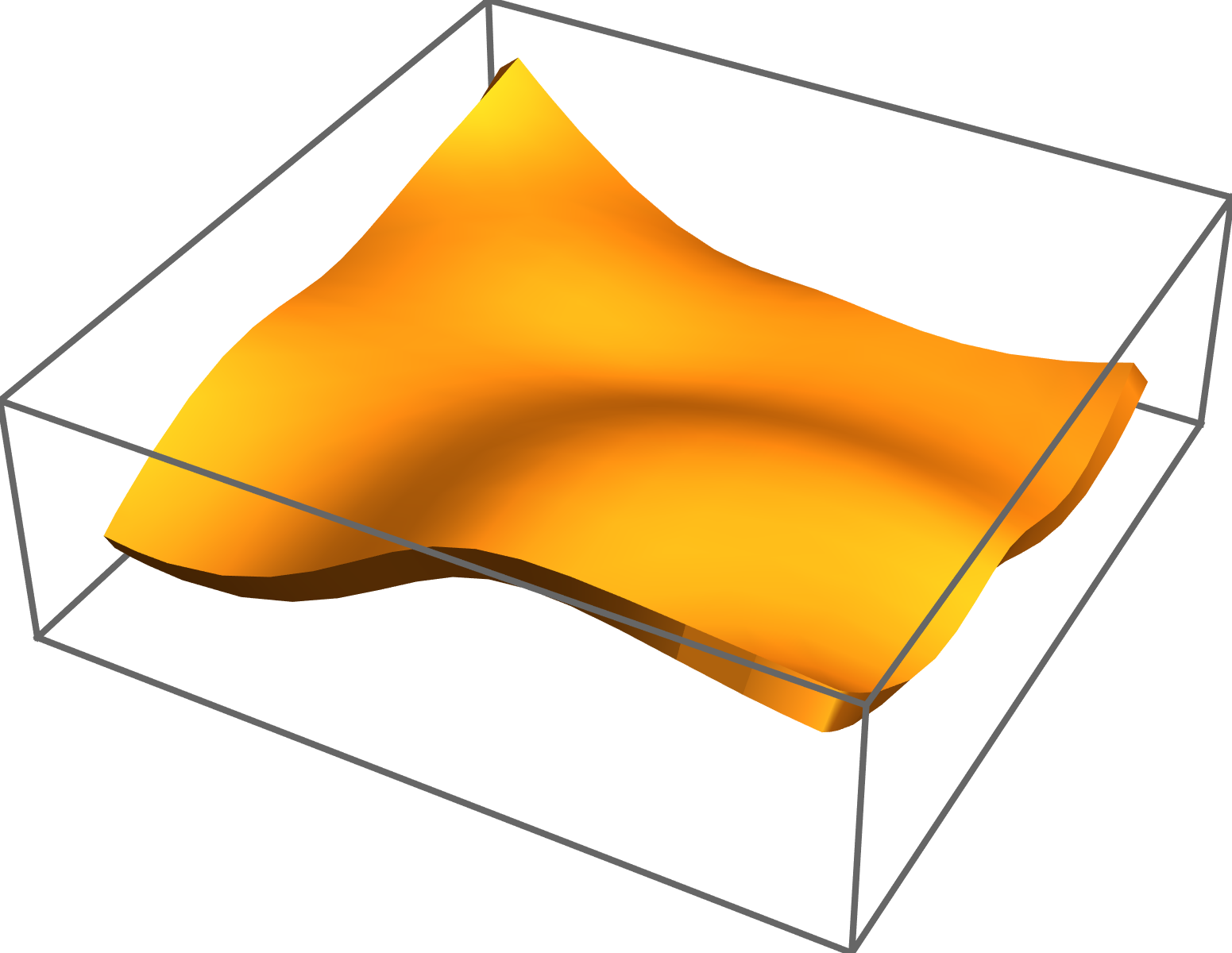}

\subsection{Getting to Yes: Meshes and Manual Discretization}
\label{sec:discretize}

Mathematica uses mesh-based geometric regions to make three-dimensional objects 3D printable. This is similar to how complex numerical expressions are converted to a decimal
\begin{mmaCell}[moredefined={Rasterize,RegularPolygon}]{Code}
N[Sqrt[5 + Sqrt[3]]]

\end{mmaCell}
\begin{mmaCell}{Output}
2.59462

\end{mmaCell}
or how a two-dimensional scene is represented as a collection of pixels.
\begin{mmaCell}[moredefined={Rasterize,RegularPolygon}]{Code}
Rasterize[RegularPolygon[7], ImageSize -> 100]

\end{mmaCell}
\mmaCellGraphics[ig={height=1in},pole2=vc,yoffset=.5ex]{Output}{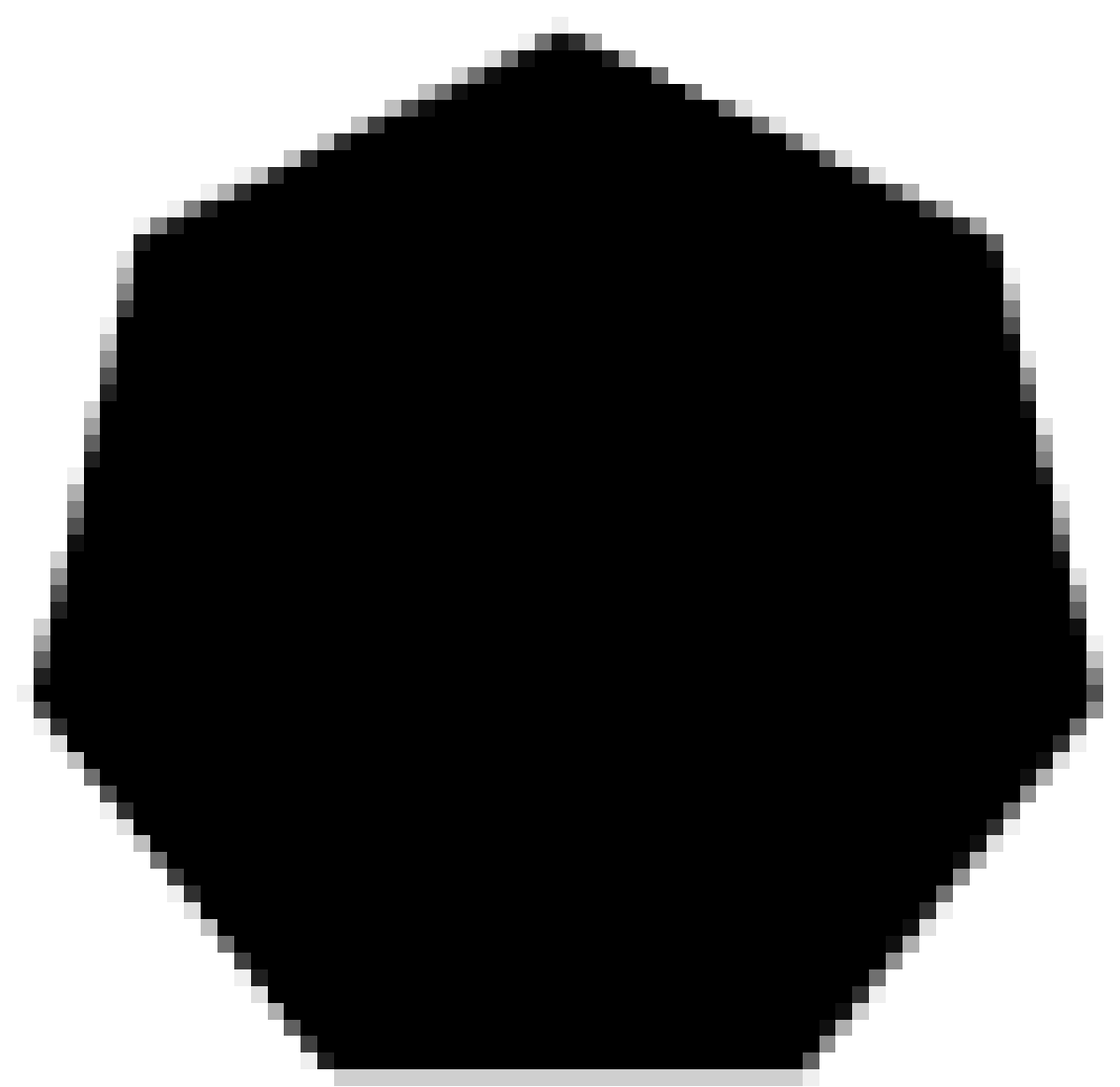}

The process of approximating an ideal object by an approximation is called discretization (or triangulation).  This replaces the smooth object by a discrete object made up of cells: points, line segments, triangles, tetrahedra, and higher-dimensional simplices. There are two main ways to describe a discretized object: either by the cells that comprise it or by the cells of its boundary.  Compare the following two discretizations of the circle:

\begin{mmaCell}[moredefined={DiscretizeRegion,BoundaryDiscretizeRegion}]{Code}
{DiscretizeRegion[Disk[]], 
 BoundaryDiscretizeRegion[Disk[]]}
 
\end{mmaCell}
\mmaCellGraphics[ig={height=1in},pole2=vc,yoffset=.5ex]{Output}{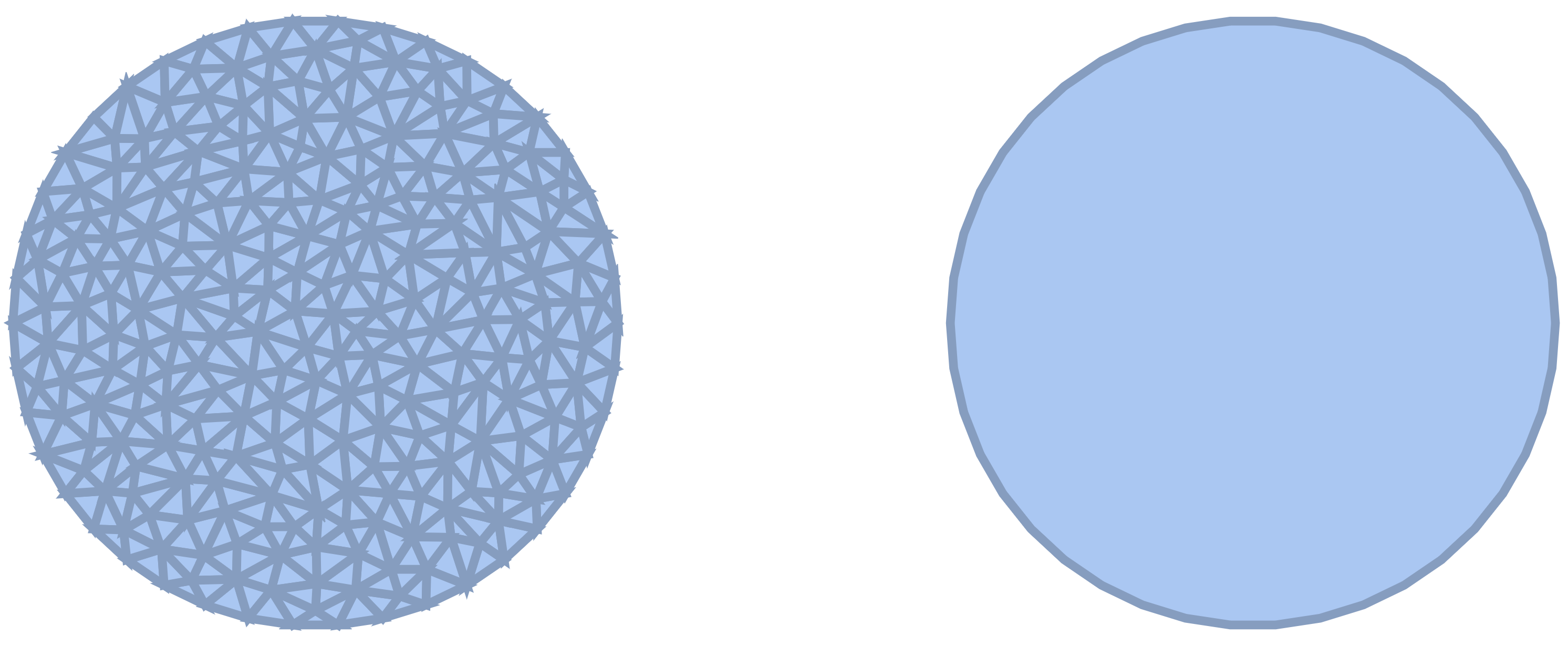}

In the first example, the circle is approximated as a collection of two-dimensional triangles. In the second example, the circle is approximated as the interior of a collection of one-dimensional line segments. In the same way, a three-dimensional form can be represented as a collection of three-dimensional tetrahedra or as the interior of a collection of two-dimensional triangles. 

When you discretize a graphics object, use {\tt DiscretizeGraphics} or {\tt Boundary\textrm{-} DiscretizeGraphics}. When you discretize a shape that is defined as a region, use {\tt DiscretizeRegion} or {\tt BoundaryDiscretizeRegion}.  The resulting objects are {\tt MeshRegion} objects or {\tt BoundaryMeshRegion} objects, respectively.

Mathematica has a large number of built-in polyhedra, accessible using the command {\tt PolyhedronData}; their meshes can be accessed directly from the command by specifying {\tt "MeshRegion"} or {\tt "BoundaryMeshRegion"}.

\begin{mmaCell}[moredefined={PolyhedronData}]{Code}
PolyhedronData["Dodecahedron", "BoundaryMeshRegion"]

\end{mmaCell}
\mmaCellGraphics[ig={height=1in},pole2=vc,yoffset=.5ex]{Output}{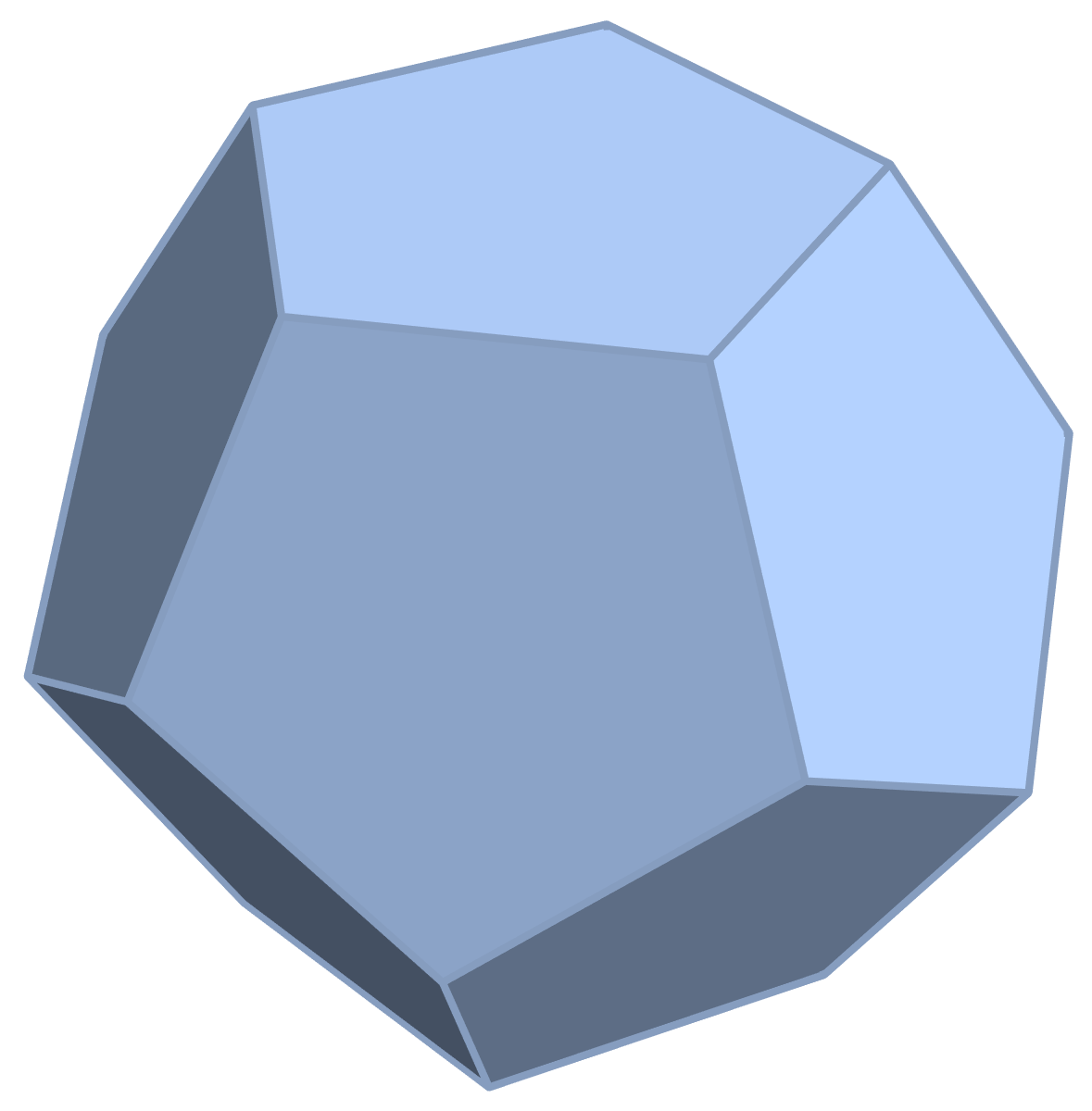}

Since boundary representations are one lower dimension than cell representations, they are computationally simpler so more of Mathematica's functionality applies.  However, not all shapes have a nice boundary representation; in that case you will have to work with the shape's cell representation. Wolfram developers have been continually improving these commands; great strides have been made in the past few years.

When you use {\tt Export}, (I think) Mathematica tries its best to discretize each piece of your model, then combines everything together.  This process sometimes goes awry, so it may make sense to discretize your objects first, use a {\tt Show} command to combine them, and export the result. Do note that it may be beneficial to discretize every single object individually instead of trying to discretize a collection of objects. Compare the following.
\begin{mmaCell}[moredefined={RandomReal,DiscretizeGraphics,Sphere}]{Code}
balls = Table[
   Sphere[{i, j, 0}, RandomReal[]], {i, 10}, {j, 10}];
{DiscretizeGraphics[balls], 
   Show[Map[DiscretizeGraphics, balls]]}
   
\end{mmaCell}
\mmaCellGraphics[ig={height=1.25in},pole2=vc,yoffset=.5ex]{Output}{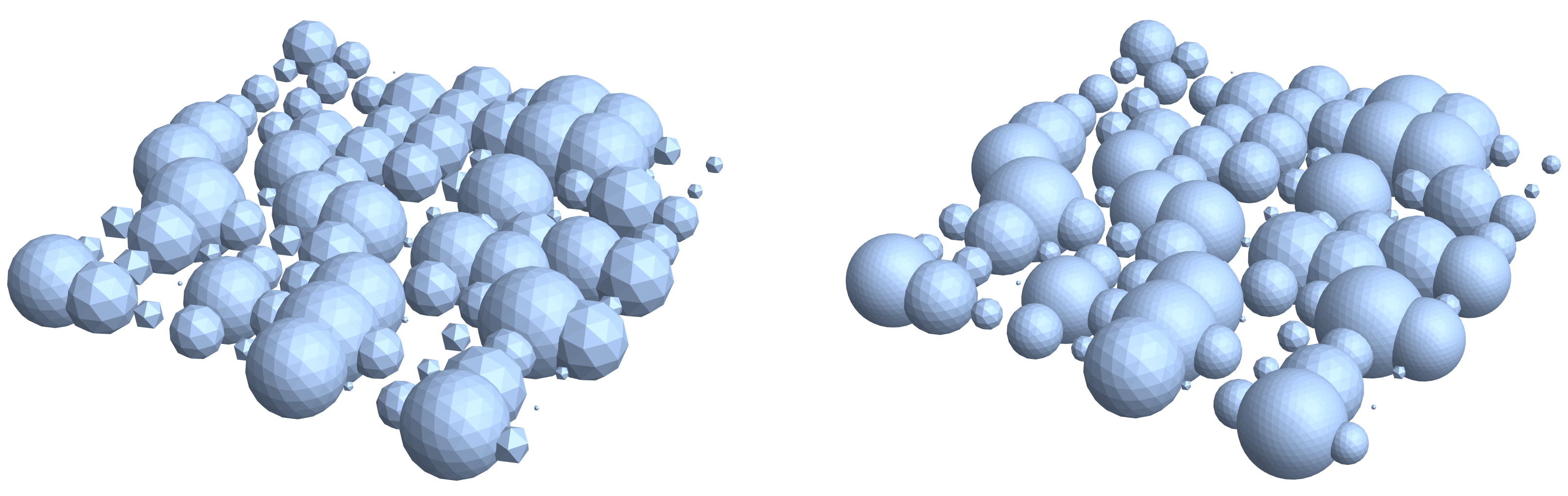}
The resolution of the second code block is better because {\tt DiscretizeGraphics} is being applied to each sphere individually. 

Furthermore, because {\tt Export} appears to try much harder than the various {\tt Discretize} commands, sometimes it is even beneficial to export and re-import each piece of the model and assemble those pieces together.

It is possible to specify the quality of a discretization of a 3D model by applying the {\tt MaxCellMeasure} option which specifies the maximum size for each highest-dimension cell. Compare:

\begin{mmaCell}[moredefined={BoundaryDiscretizeGraphics,Ball,MaxCellMeasure}]{Code}
{BoundaryDiscretizeGraphics[Ball[], MaxCellMeasure->.01],	
 BoundaryDiscretizeGraphics[Ball[], MaxCellMeasure->.001],
 BoundaryDiscretizeGraphics[Ball[], MaxCellMeasure->.0001]}
 
\end{mmaCell}
\mmaCellGraphics[ig={height=1.25in},pole2=vc,yoffset=.5ex]{Output}{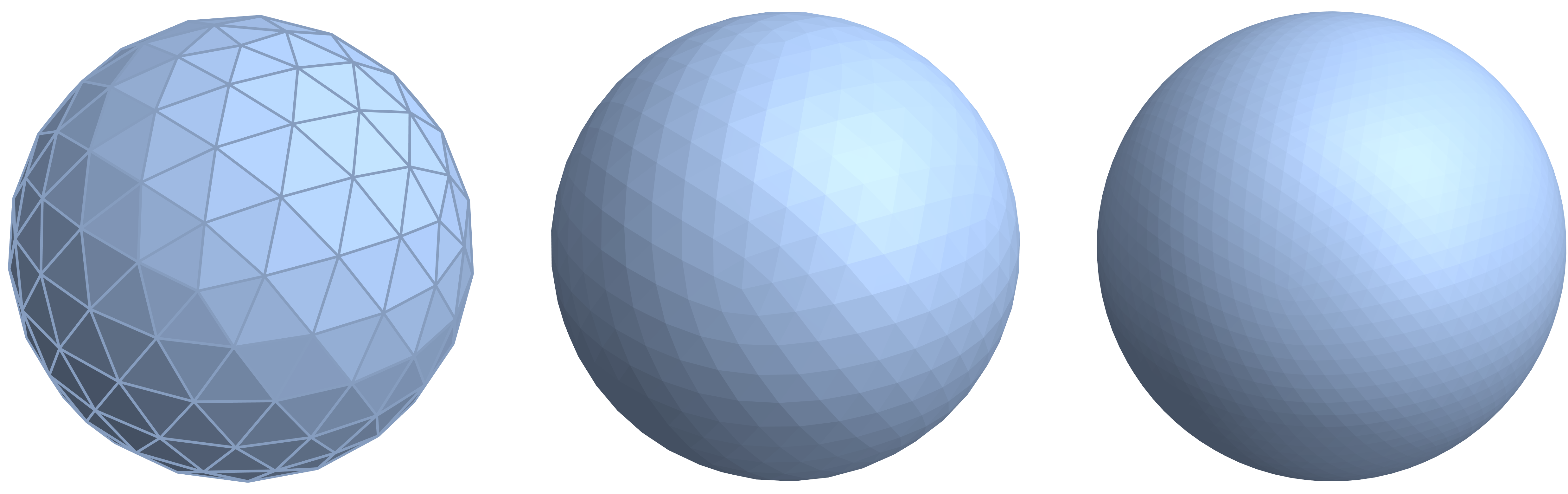}

As of Mathematica version 13, the various region objects can be imported into a {\tt Graphics3D} scene as primitives, which makes it possible to finalize each individual piece of the model and then color them independently for export as a WRL file, which is how the snowman is improved in Section~\ref{sec:snowman}. 

In my experience, applying some combination of these techniques allow the vast majority of models to be exported by Mathematica. If you are still having issues, see Section~\ref{sec:RegionPlot} for another tool.

\subsection{Extracting Mesh Information}

There are two useful commands to situate your object in 3D space. {\tt RegionBounds} gives the smallest and largest $x$, $y$, and $z$ values over the whole object. (Do not confuse {\tt RegionBounds} with {\tt RegionBoundary} which finds a region's geometric boundary.) And {\tt RegionCentroid} gives the centroid of the object. We observe that the following discretized ellipse is $3.8$ units wide and is centered at $(0,0,0)$.

\begin{mmaCell}[moredefined={BoundaryDiscretizeGraphics,Ellipsoid,MaxCellMeasure,RegionBounds,RegionCentroid}]{Code}
mesh = BoundaryDiscretizeGraphics[
  Ellipsoid[{0, 0, 0}, {2, 1, 1}], MaxCellMeasure -> 0.2];
RegionBounds[mesh]
Map[#[[2]] - #[[1]] &, 
Chop[RegionCentroid[mesh]]

\end{mmaCell}

\begin{mmaCell}{Output}
\{\{-1.90211, 1.90211\}, \{-1., 1.\}, \{-1., 1.\}\}
\{3.80423, 2., 2.\}
\{0, 0, 0\}

\end{mmaCell}

As we build complexity into our models, it becomes useful to work directly with the internal structure of the mesh. Use {\tt MeshCoordinates} to extract the coordinates of the mesh and {\tt MeshPrimitives} to extract higher-dimensional pieces of the object, like its edges or faces.  The second input is the dimension of the parts you want to extract.

\begin{mmaCell}[moredefined={MeshCoordinates,MeshPrimitives}]{Code}
coords = MeshCoordinates[mesh];
edges = MeshPrimitives[mesh, 1]; 
faces = MeshPrimitives[mesh, 2]; 

\end{mmaCell}

These commands allow you to efficiently make a wireframe version of your object or algorithmically color the faces of your object.

\begin{mmaCell}[moredefined={Sphere,Tube}]{Code}
thickness = 0.06;
Graphics3D[{
  Map[Sphere[#, thickness] &, coords],
  Map[Tube[#, thickness] &, edges]
}, Boxed -> False, SphericalRegion -> True]

\end{mmaCell}
\mmaCellGraphics[ig={height=1in},pole2=vc,yoffset=.5ex]{Output}{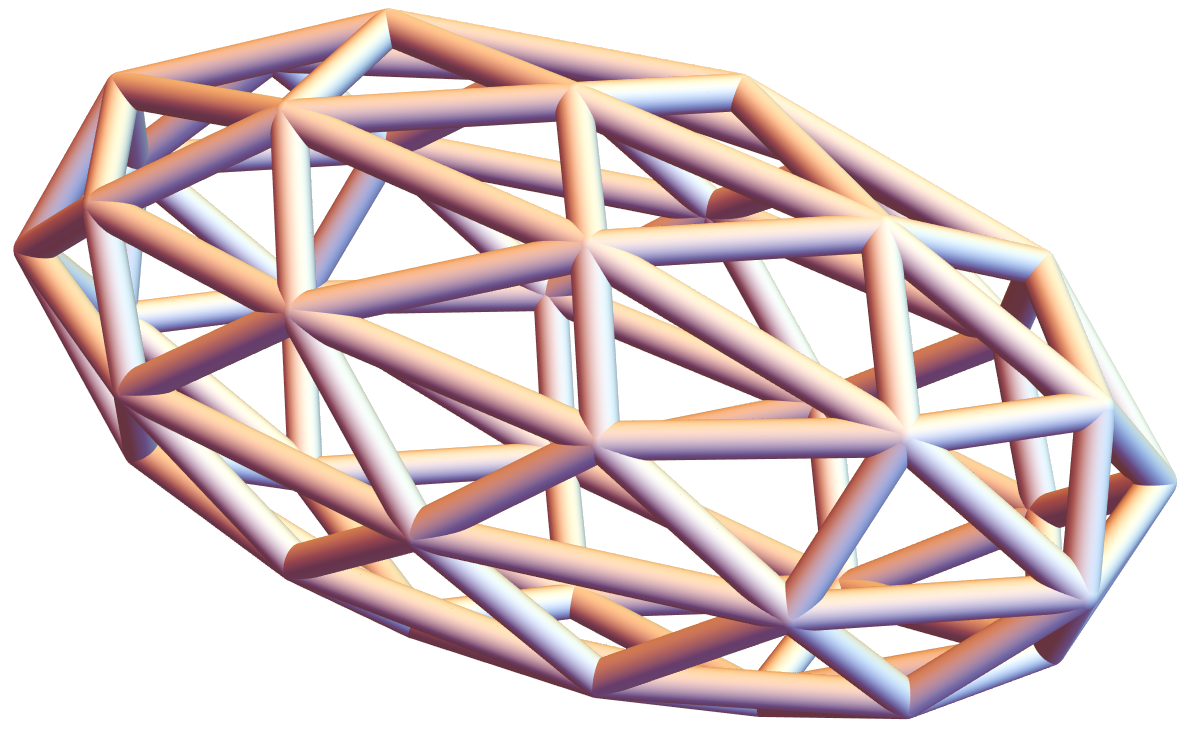}

\begin{mmaCell}[moredefined={RegionCentroid}]{Code}
colors = Map[Hue[(RegionCentroid[#][[3]]+1)/2] &, faces];
Graphics3D[MapThread[{#2, #1} &, {faces, colors}], 
  Boxed -> False, SphericalRegion -> True]
  
\end{mmaCell}
\mmaCellGraphics[ig={height=1in},pole2=vc,yoffset=.5ex]{Output}{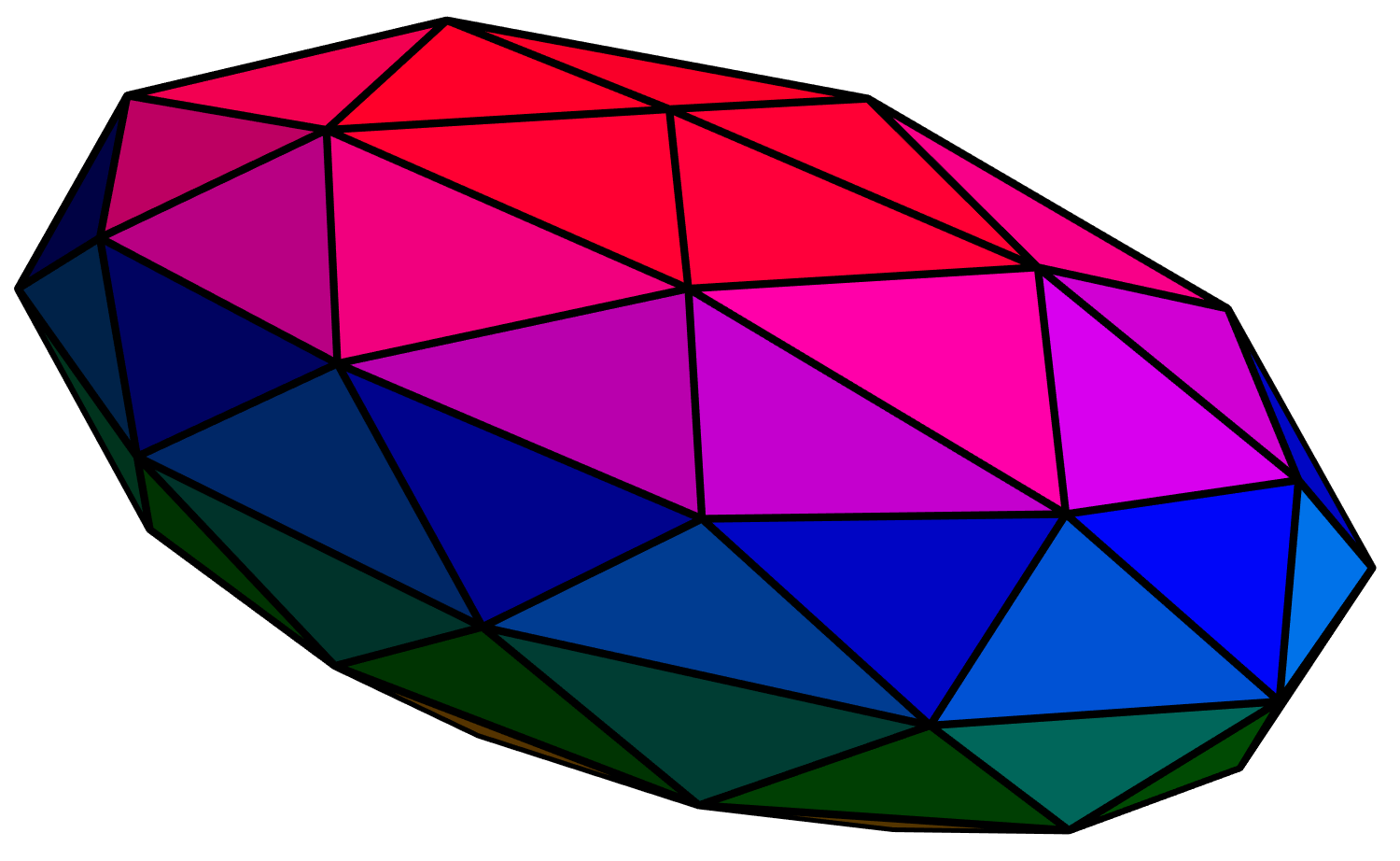}

Note that the wireframe object needs to include spheres at each vertex so that no gaps occur where two cylinders meet. Compare the 3D models generated with spheres and without:
\begin{center}
\includegraphics[height=1in]{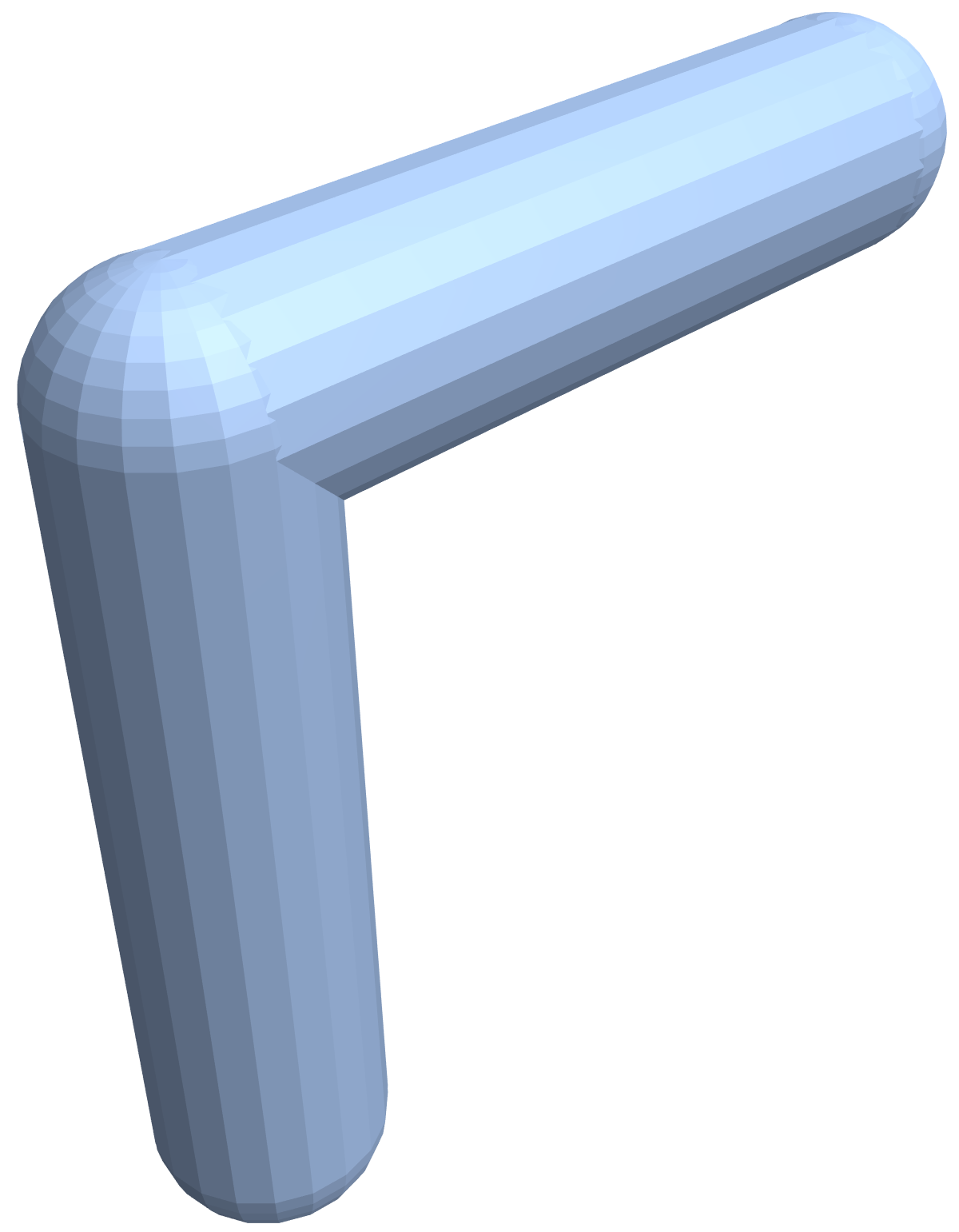}\qquad\qquad\qquad\includegraphics[height=1in]{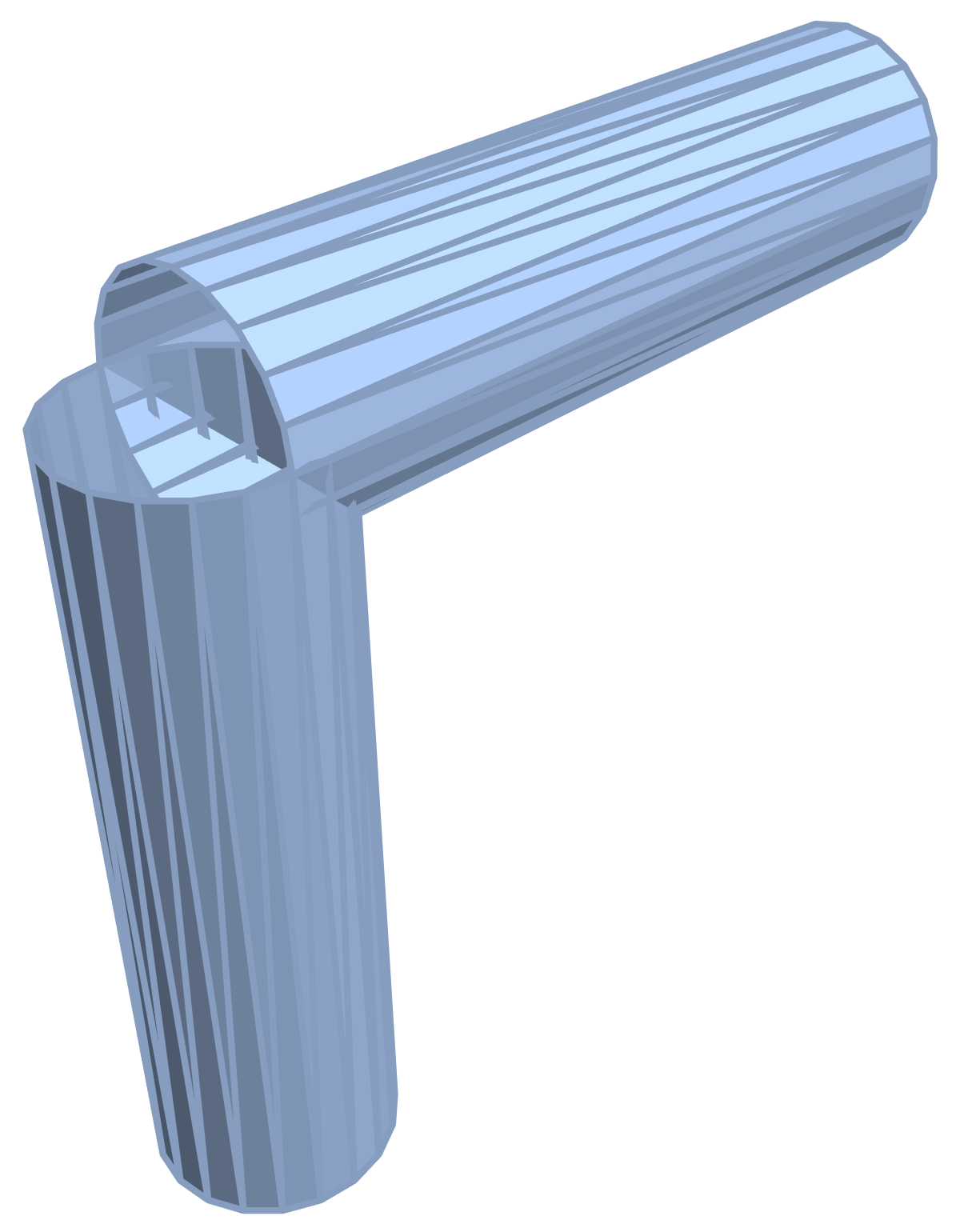}
\end{center}
The model without spheres only consists of the shell of each cylindrical edge. Because it is not solid, it would not be 3D printable.

\subsection{Region Operations}

Once you have discretized your models (or pieces of your models) using the discretizing commands or by exporting and reimporting, you have {\tt RegionMesh} or {\tt BoundaryRegionMesh} objects and you can apply various region operations to them.

It is important to be able to change the size and orientation of your model. For instance, it may be easiest to work in one coordinate system to create your mathematical form and in a different coordinate system when printing. In addition, you will want to specify the dimensions of the piece in millimeters if you are uploading your model to a 3D printing service.

To change the size of your model use {\tt RegionResize}. If your second input is one number, Mathematica will scale the object to have first coordinate equal to that input, respecting the existing box ratios.  If you instead specify a list of parameters, it will scale the object to fit into the given box, not respecting the existing box ratios.

If we want our model to be 100mm wide, we would write
\begin{mmaCell}[moredefined={RegionResize,PolyhedronData}]{Code}
mesh = PolyhedronData["Dodecahedron", "BoundaryMeshRegion"]
RegionResize[mesh, 100]

\end{mmaCell}
\mmaCellGraphics[ig={height=1in},pole2=vc,yoffset=.5ex]{Output}{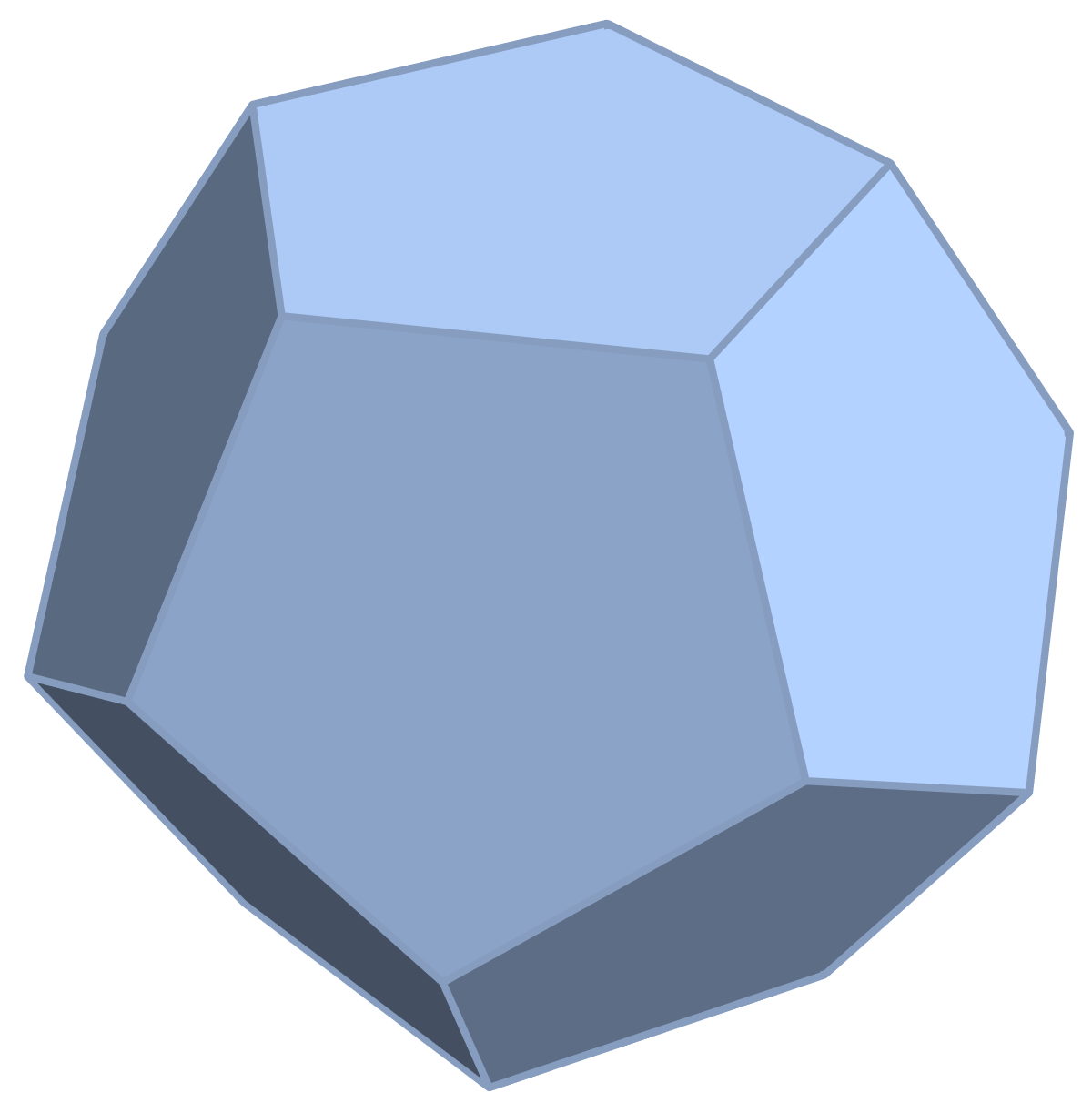}
or if we want to skew it to fit in a box that is 80mm $\times$ 50mm $\times$ 30mm, we would write
\begin{mmaCell}[moredefined={RegionResize,PolyhedronData}]{Code}
mesh = PolyhedronData["Dodecahedron", "BoundaryMeshRegion"]
RegionResize[mesh, {80, 50, 30}]

\end{mmaCell}
\mmaCellGraphics[ig={height=0.75in},pole2=vc,yoffset=.5ex]{Output}{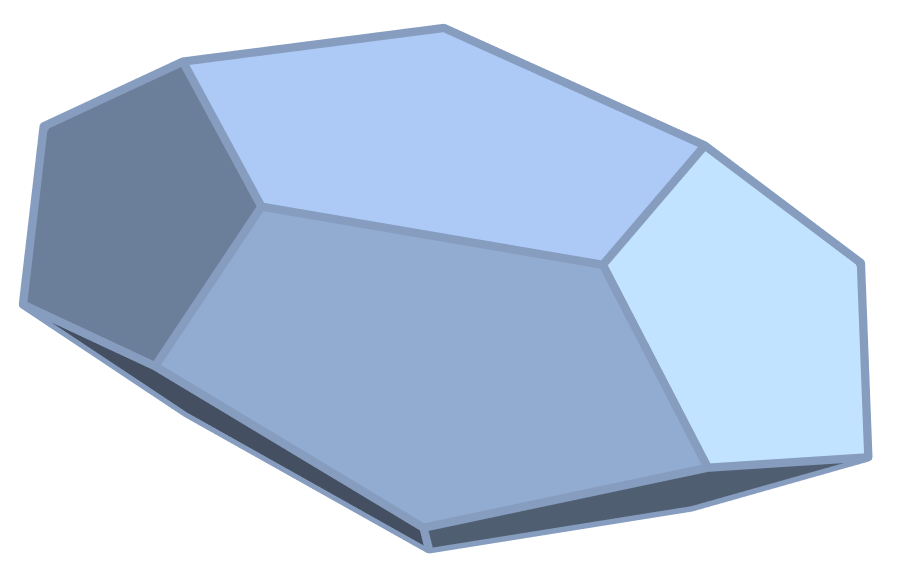}

Use {\tt TransformedRegion} and the various {\tt Transform} commands to apply transformations to a region. Suppose we want to apply transformations to the basic cube:
\begin{mmaCell}[moredefined={BoundaryDiscretizeGraphics}]{Code}
cube = BoundaryDiscretizeGraphics[Cuboid[]];

\end{mmaCell}
To rotate the object by a given angle $\theta$ around an axis of rotation (with direction $v$ and passing through point $p$), input these values into {\tt RotationTransform}.
\begin{mmaCell}[moredefined={TransformedRegion,RotationTransform}]{Code}
rotatedcube = TransformedRegion[cube, 
  RotationTransform[Pi/4, {0, 0, 1}, {1/2, 1/2, 0}]];
Graphics3D[{Green, cube, Red, rotatedcube}, 
  Boxed -> False]
  
\end{mmaCell}
\mmaCellGraphics[ig={height=1in},pole2=vc,yoffset=.5ex]{Output}{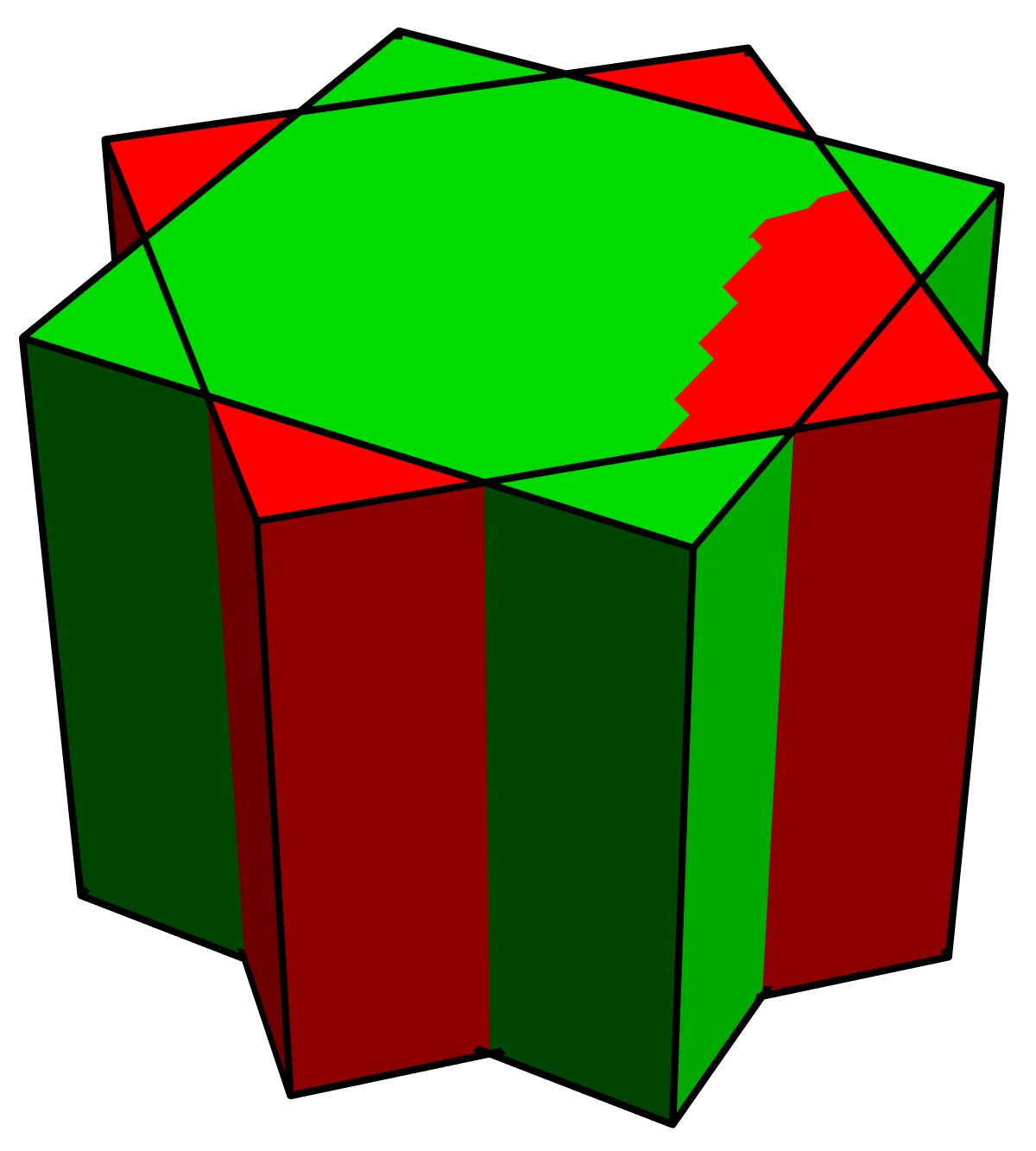}
To translate the object along a vector $v$, input that vector into {\tt Translation\textrm{-} Transform}. 
\begin{mmaCell}[moredefined={TransformedRegion,TranslationTransform}]{Code}
translatedcube = TransformedRegion[cube, 
  TranslationTransform[{0.5, 0.5, 0.5}]];
Graphics3D[{Green, cube, Red, translatedcube},
  Boxed -> False]
  
\end{mmaCell}
\mmaCellGraphics[ig={height=1.5in},pole2=vc,yoffset=.5ex]{Output}{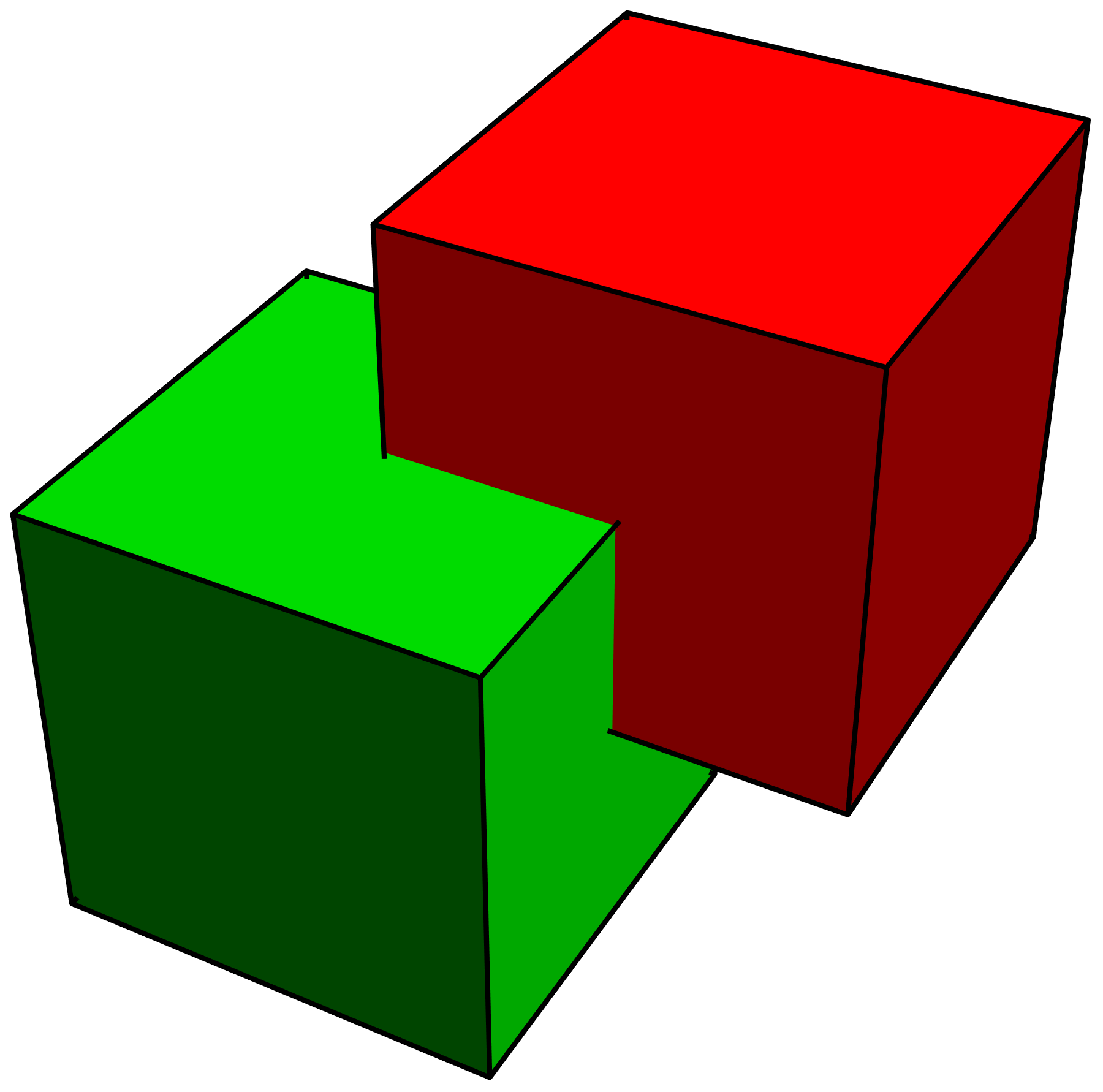}

Combining a few of these techniques is useful to orient your object. For example, if you want to move the center of your object to the origin, apply a translation equal to the negative of its centroid:
\begin{mmaCell}[moredefined={TransformedRegion,TranslationTransform,RegionCentroid}]{Code}
TransformedRegion[mesh, 
  TranslationTransform[-RegionCentroid[mesh]]],
  
\end{mmaCell}

To intersect, merge, or subtract two regions, use {\tt RegionIntersection}, {\tt Region\textrm{-} Union}, and {\tt RegionDifference}.  These commands work well in 2D. On the other hand, it has been hit or miss whether these commands work on 3D MeshRegion objects. Starting in version 12.1 and further improved in version 13 of Mathematica, these boolean operations work more often for 3D boundary representations.
\begin{mmaCell}[moredefined={BoundaryDiscretizeGraphics,Ball,RegionResize,PolyhedronData,RegionUnion,RegionIntersection,RegionDifference}]{Code}
ball = BoundaryDiscretizeGraphics[
  Graphics3D[Ball[{0, 0, 0}, 1]]];
cube = BoundaryDiscretizeGraphics[
  Graphics3D[Cuboid[{0.1, -0.5, 0.1}, {1, 0.5, 1}]]];
dodec = RegionResize[
  PolyhedronData["Dodecahedron", "BoundaryMeshRegion"], 2];
{RegionUnion[ball, cube], 
 RegionUnion[dodec, cube], 
 RegionIntersection[ball, cube], 
 RegionDifference[ball, cube],
 RegionDifference[dodec, cube]}
 
\end{mmaCell}
\mmaCellGraphics[ig={height=0.9in},pole2=vc,yoffset=.5ex]{Output}{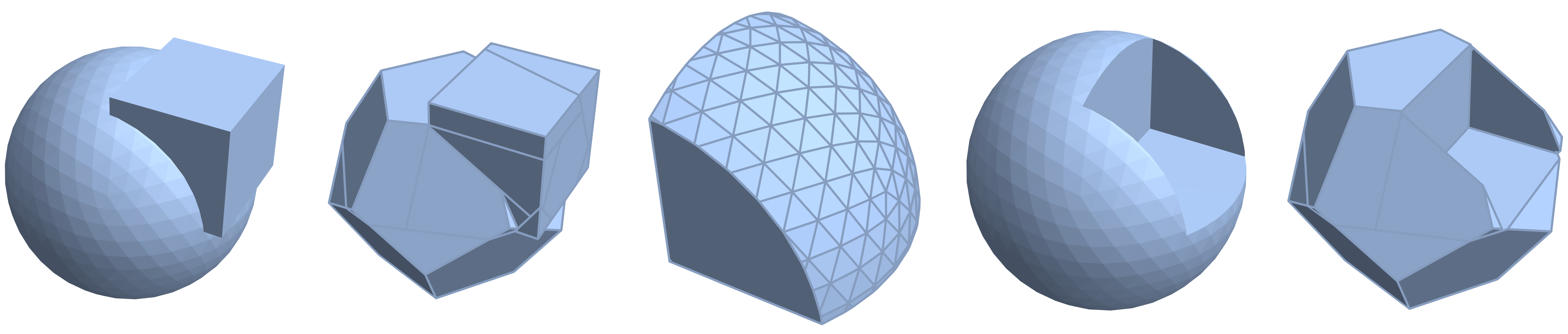}

\subsection{RegionPlot3D}
\label{sec:RegionPlot}

Another command that can achieve high resolution models is {\tt RegionPlot3D}. It can be very computationally intensive but it is sometimes the only command that works. This command takes as input a boolean expression that defines your desired region and the domain that you want sampled. {\tt RegionPlot3D} gives all the points in the specified domain that satisfy your given requirements.  

For example, if you want to find the intersection of two unit spheres centered at $(0,0,0)$ and $(1,0,0)$, you specify the inequality that defines each sphere and use {\tt \&\&} (the syntax for `and') to indicate that you want their intersection, where both inequalities are true.
\begin{mmaCell}[moredefined={RegionPlot3D}]{Code}
RegionPlot3D[
  (x^2 + y^2 + z^2 <= 1) && ((x - 1)^2 + y^2 + z^2 <= 1),
  {x, -1, 2}, {y, -1.5, 1.5}, {z, -1.5, 1.5}]
  
\end{mmaCell}
\mmaCellGraphics[ig={height=1.5in},pole2=vc,yoffset=.5ex]{Output}{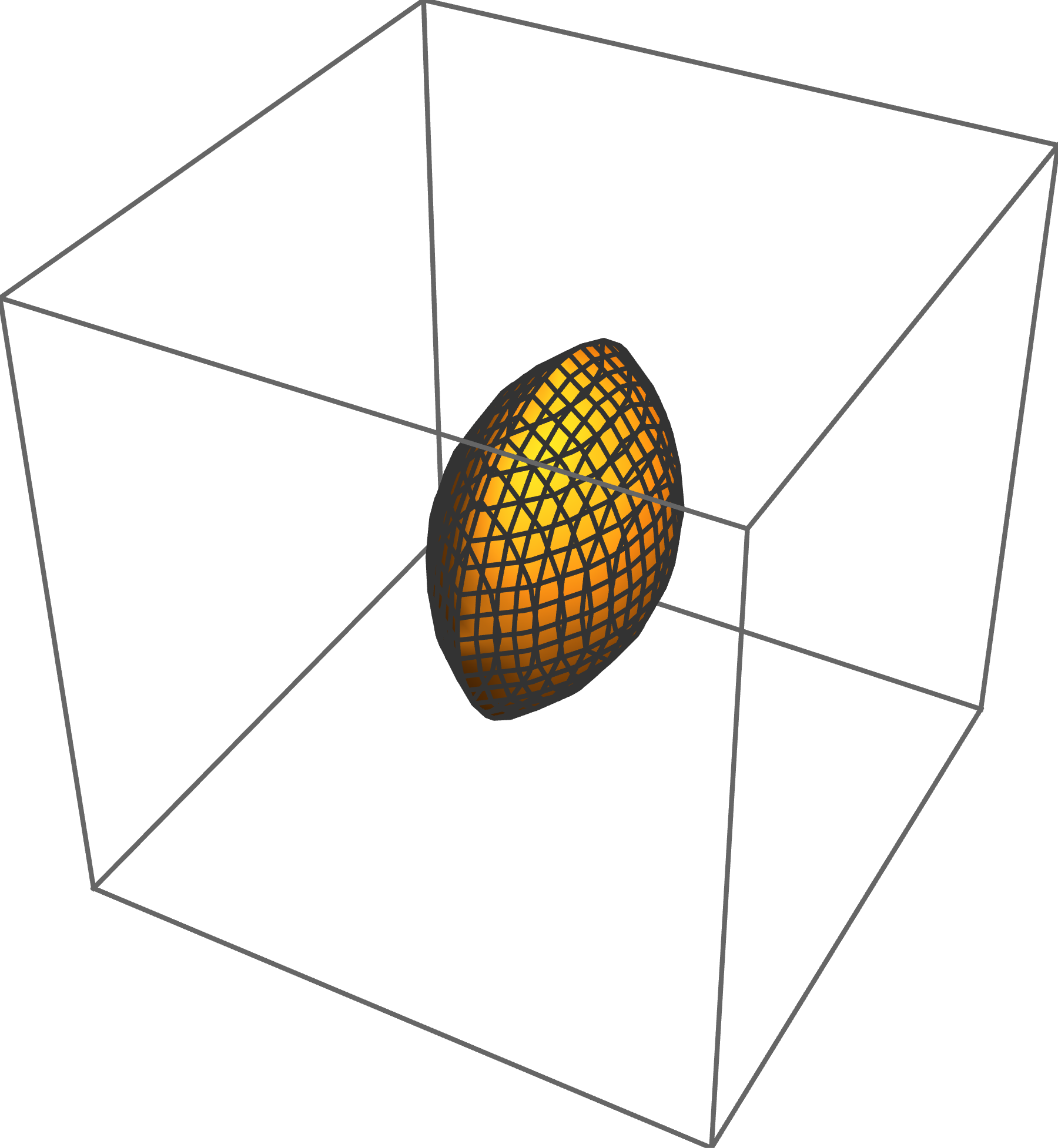}
or use {\tt ||} (the syntax for `or') to indicate that you want their union, where at least one of the two inequalities is true.
\begin{mmaCell}[moredefined={RegionPlot3D}]{Code}
RegionPlot3D[
  (x^2 + y^2 + z^2 <= 1) || ((x - 1)^2 + y^2 + z^2 <= 1),
  {x, -1, 2}, {y, -1.5, 1.5}, {z, -1.5, 1.5}]
  
\end{mmaCell}
\mmaCellGraphics[ig={height=1.5in},pole2=vc,yoffset=.5ex]{Output}{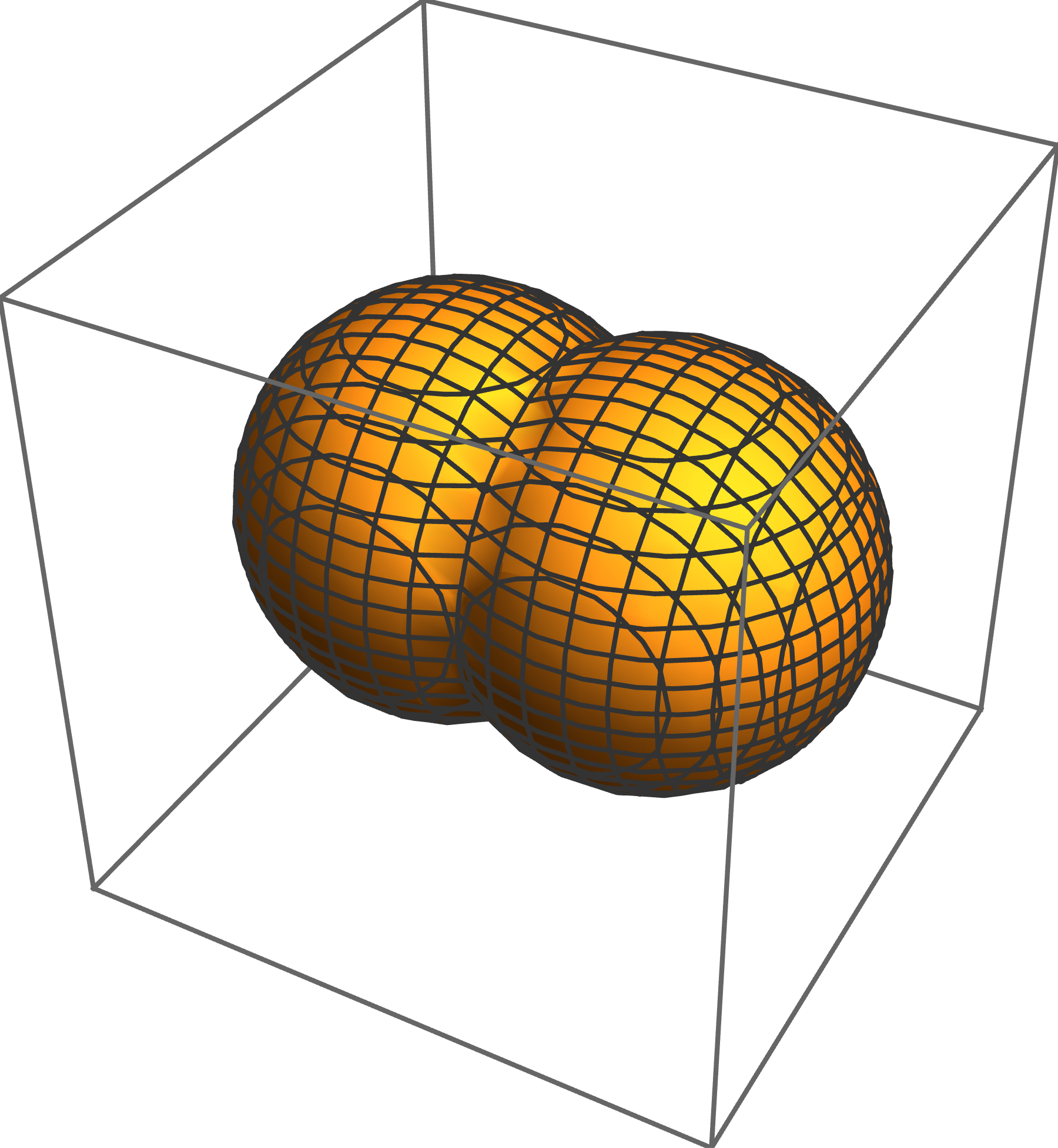}

{\tt RegionPlot3D} is a nice way to create 3D models of surfaces defined by an equation $z=f(x,y)$ that have a flat base for display. We do that by specifying the region as the set of points that lie below the curve. 
\begin{mmaCell}[moredefined={RegionPlot3D}]{Code}
block = RegionPlot3D[Sin[x] + Sin[y] > z, 
  {x, 0 - Pi/2, 6 Pi - Pi/2}, {y, 0 - Pi/2, 6 Pi - Pi/2}, 
  {z, -3, 3}, BoxRatios -> Automatic];
Import[Export[NotebookDirectory[] <> "Bumpy.stl", block]]

\end{mmaCell}
\mmaCellGraphics[ig={height=1in},pole2=vc,yoffset=.5ex]{Output}{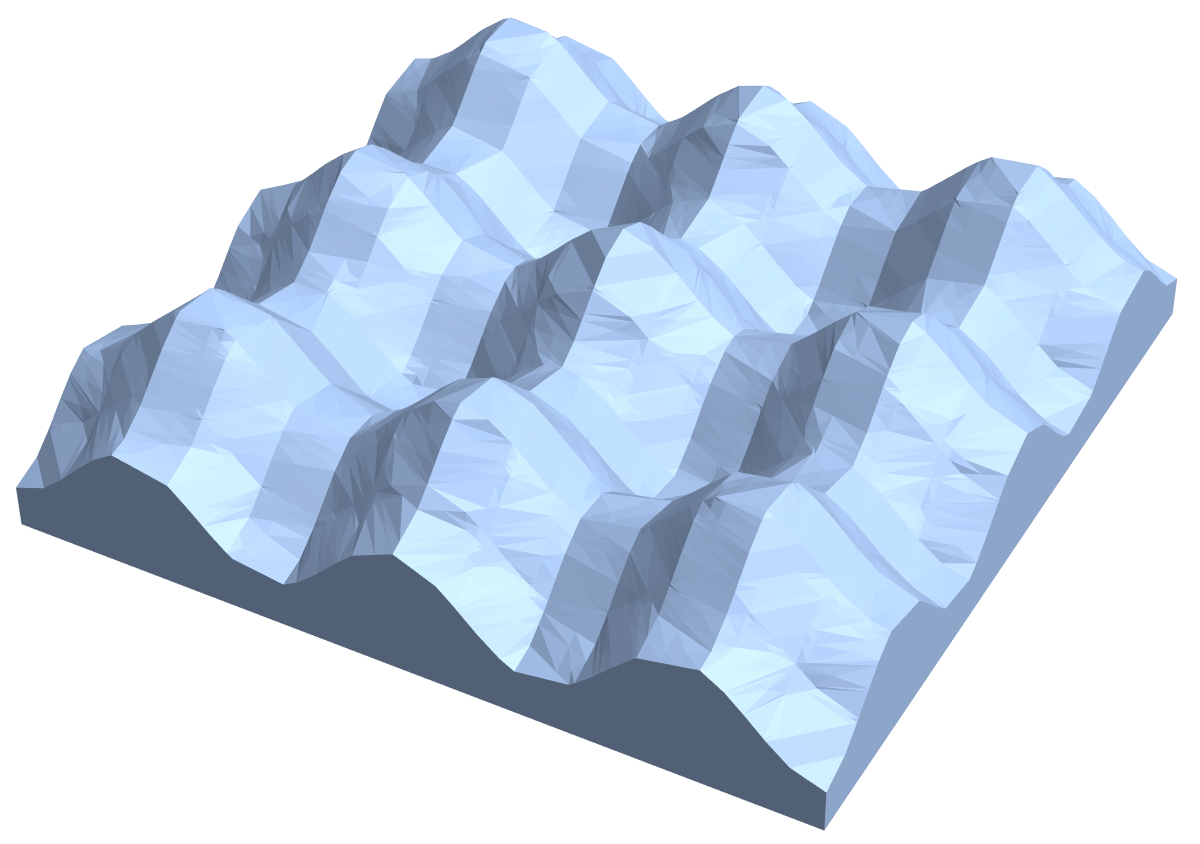}
I assign the option {\tt BoxRatios} to {\tt Automatic} so that the axes are scaled to ensure the unit distance in each direction is the same. 

We see that the quality of the surface could be improved. We can specify the number of points that Mathematica should sample in each direction using the {\tt PlotPoints} option.
\begin{mmaCell}[moredefined={RegionPlot3D}]{Code}
blockHD = RegionPlot3D[Sin[x] + Sin[y] > z, 
  {x, 0 - Pi/2, 6 Pi - Pi/2}, {y, 0 - Pi/2, 6 Pi - Pi/2}, 
  {z, -3, 3}, PlotPoints -> {100, 100, 100}, 
  BoxRatios -> Automatic];
Import[Export[
  NotebookDirectory[] <> "Bumpy.HD.stl", blockHD]]
  
\end{mmaCell}
\mmaCellGraphics[ig={height=1in},pole2=vc,yoffset=.5ex]{Output}{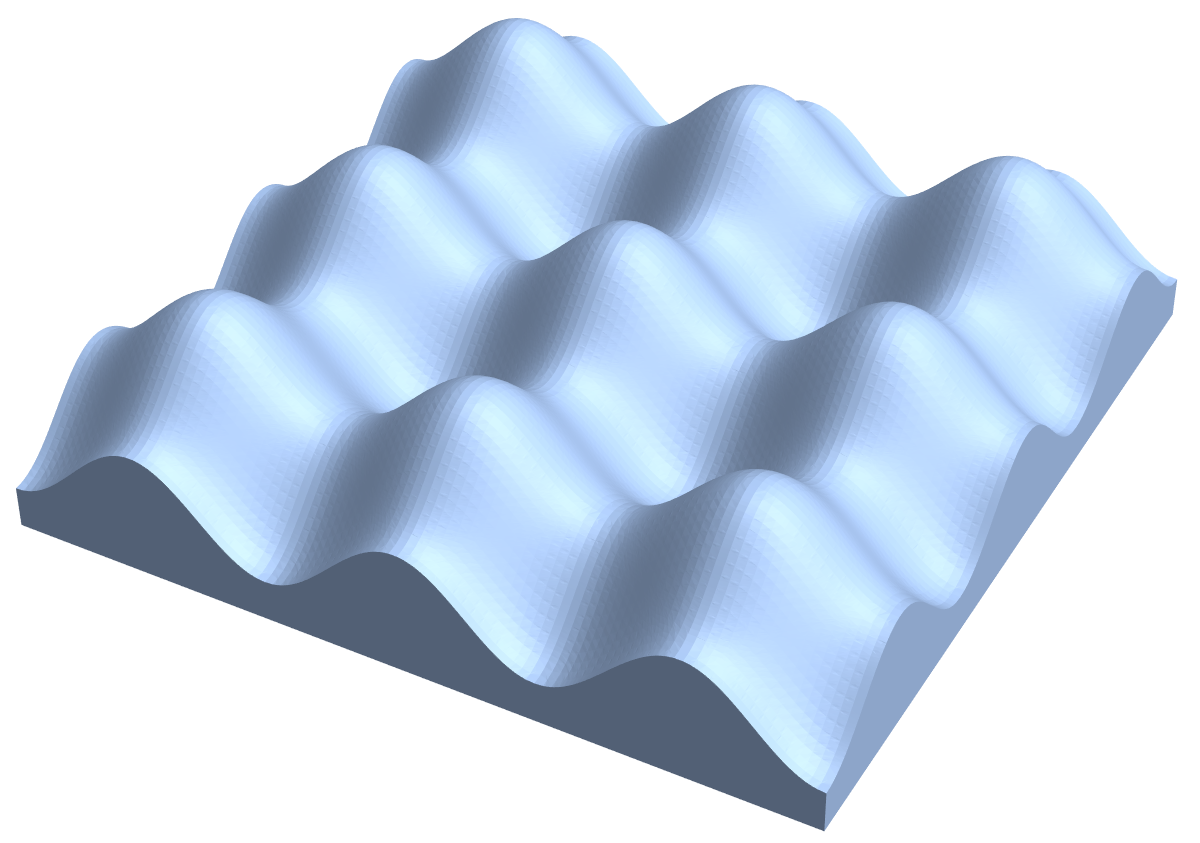}
The more points you specify, the better the result. However, this becomes computationally intensive very quickly as you increase the number of points and it rarely produces as precise of a model as we had using the earlier methods.

\subsection{High Quality Snowman}
\label{sec:snowman}

Below is code that creates a high quality, full color 3D model of the snowman from Section~\ref{sec:basic3D}. It took multiple iterations to make it work---I ran into errors to export to a {\tt .wrl} file until I converted every mesh to its polygonal representation, hence the use of {\tt MeshPrimitives}. 
\begin{mmaCell}[moredefined={BoundaryDiscretizeRegion,BoundaryDiscretizeGraphics,MaxCellMeasure,Ball,MeshPrimitives,Cylinder,ConeRegionResize,PolyhedronData,RegionUnion,RegionIntersection,RegionDifference}]{Code}
snowball1 = MeshPrimitives[BoundaryDiscretizeRegion[
  Ball[{0, 0, 0}, 1], MaxCellMeasure -> .0005], 2];
snowball2 = MeshPrimitives[BoundaryDiscretizeRegion[
  Ball[{0, 0, 1.3}, 0.8], MaxCellMeasure -> .0005], 2];
snowball3 = MeshPrimitives[BoundaryDiscretizeRegion[
  Ball[{0, 0, 2.4}, 0.6], MaxCellMeasure -> .0005], 2];
tophat1 = MeshPrimitives[BoundaryDiscretizeRegion[
  Cylinder[{{0, 0, 2.8}, {0, 0, 2.9}}, 0.7], 
  MaxCellMeasure -> 0.0005], 2];
tophat2 = MeshPrimitives[BoundaryDiscretizeRegion[
  Cylinder[{{0, 0, 2.9}, {0, 0, 3.5}}, 0.5], 
  MaxCellMeasure -> 0.0005], 2];
carrot = MeshPrimitives[BoundaryDiscretizeRegion[
  Cone[{{0, -0.55, 2.4}, {0, -0.9, 2.4}}, 0.1], 
  MaxCellMeasure -> .00005], 2];
button1 = MeshPrimitives[BoundaryDiscretizeGraphics[
  Cuboid[{0, -0.75, 1.3} - 0.1 {1, 1, 1}, 
         {0, -0.75, 1.3} + 0.1 {1, 1, 1}]], 2]; 
button2 = MeshPrimitives[BoundaryDiscretizeGraphics[
  Cuboid[{0, -0.7, 1.6} - 0.1 {1, 1, 1}, 
         {0, -0.7, 1.6} + 0.1 {1, 1, 1}]], 2]; 
button3 = MeshPrimitives[BoundaryDiscretizeGraphics[
  Cuboid[{0, -0.7, 1.0} - 0.1 {1, 1, 1}, 
         {0, -0.7, 1.0} + 0.1 {1, 1, 1}]], 2];
snowmanHD = Graphics3D[{EdgeForm[None],
  White, snowball1, snowball2, snowball3, 
  Black, tophat1, tophat2, button1, button2, button3,
  Orange, carrot}, 
  Boxed -> False, SphericalRegion -> True]
Export[NotebookDirectory[] <> "snowmanHD.wrl", snowmanHD]

\end{mmaCell}
You can interact with the final model on Sketchfab at \url{https://skfb.ly/ovIVI}. 
\begin{center}
	\includegraphics[height=3in]{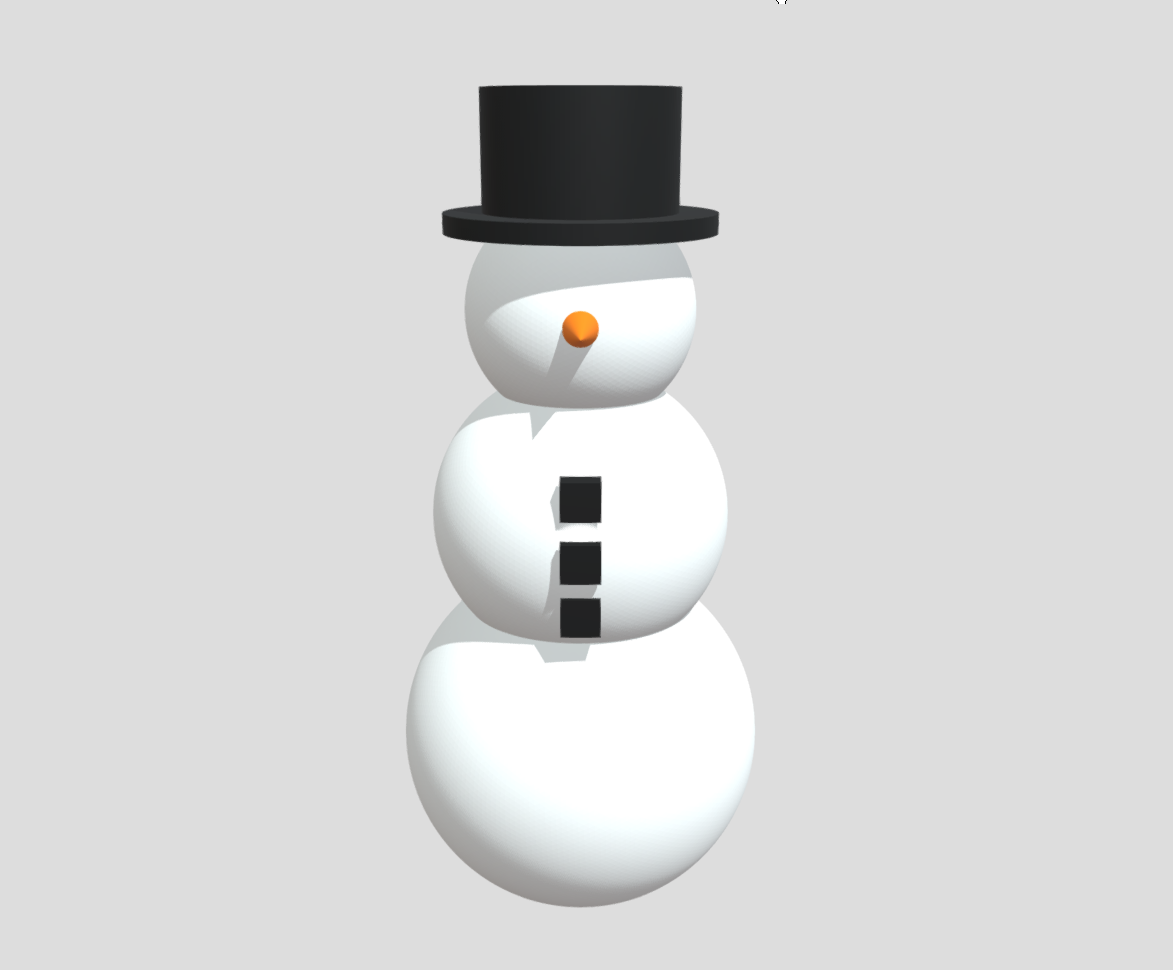}
\end{center}

\bibliographystyle{amsalpha}

\end{document}